\documentclass{amsart}

\usepackage{amsfonts, amsmath, amssymb, amsthm, amsxtra, latexsym}

\theoremstyle{plain}

\newtheorem{theorem}{Theorem}[section]
\newtheorem{proposition}[theorem]{Proposition}
\newtheorem{lemma}[theorem]{Lemma}
\newtheorem{corollary}[theorem]{Corollary}

\theoremstyle{definition}

\newtheorem{definition}[theorem]{Definition}
\newtheorem{example}[theorem]{Example}
\newtheorem{remark}[theorem]{Remark}

\numberwithin{equation}{section}

\newcommand {\C}{{{\mathbb C}}}
\newcommand {\N}{{{\mathbb N}}}
\newcommand {\R}{{{\mathbb R}}}
\newcommand {\T}{{{\mathbb T}}}
\newcommand {\Z}{{{\mathbb Z}}}
\newcommand {\Cn}{{\C^n}}
\newcommand {\Nn}{{\N^n}}
\newcommand {\Rn}{{\R^n}}
\newcommand {\Tn}{{\T^n}}
\newcommand {\Zn}{{\Z^n}}

\newcommand {\choice}{{\mathcal C}}
\newcommand {\gen}{{\mathcal G}}
\newcommand {\env}{{\mathfrak U}}

\newcommand {\testX}{{\mathcal T}}
\newcommand {\testXdual}{\testX^\prime}

\newcommand {\test}{{\mathcal D}}
\newcommand {\dist}{{{\test}^\prime}}
\newcommand {\testcOmega}{C_c^\infty(\Omega)}
\newcommand {\testcmOmega}{C_c^m(\Omega)}
\newcommand {\testOmega}{C^\infty(\Omega)}
\newcommand {\testzeroOmega}{C(\Omega)}
\newcommand {\testmOmega}{C^m(\Omega)}
\newcommand {\distOmega}{C_c^\infty(\Omega)^\prime}
\newcommand {\distmOmega}{C_c^m(\Omega)^\prime}

\newcommand {\compdistmOmega}{C^m(\Omega)^\prime}
\newcommand {\schw}{{\mathcal S}}
\newcommand {\tempdist}{{{\mathcal S}^\prime}}

\newcommand {\Fourier}{{\mathcal F}}
\newcommand {\Fourierinverse}{{\Fourier^{-1}}}

\newcommand {\typ}{T}
\newcommand {\Fouriertyp}{\widehat\typ}

\newcommand {\poly}{{\mathcal P}}
\newcommand {\polyd}{\poly_d}
\newcommand {\galg}{{\mathcal E}_{\text{\textup{im}}}}

\newcommand {\weights}{{\mathcal W}}
\newcommand {\weightsqa}{\weights_{\text{\textup{qa}}}}
\newcommand {\weightshol}{\weights_{\text{\textup{hol}}}}

\newcommand {\wtilde}{{\widetilde{w}}}
\newcommand {\wprime}{w^\prime}

\newcommand {\op}{{\mathcal O}}
\newcommand {\opw}{{\op_w}}
\newcommand {\opwtilde}{\op_\wtilde}
\newcommand {\opwdual}{{\op_w^\prime}}

\newcommand {\parta}{\partial^\alpha}
\newcommand {\pa}{p_\alpha}
\newcommand {\domain}{{\C^n_w}}
\newcommand {\emap}{\Psi}
\newcommand {\emapP}{\emap_P}
\newcommand {\charfunc}{1}

\DeclareMathOperator {\supp}{supp\,}

\DeclareMathOperator {\ima}{Im}
\DeclareMathOperator {\ord}{ord}
\DeclareMathOperator {\Span}{Span\,}
\DeclareMathOperator {\End}{End\,}
\DeclareMathOperator {\aff}{Aff}

\begin{document}


\title{Subspaces with equal closure}

\translator{}

\dedicatory{Dedicated to Gabri\"elle ten Have, Lianne, Sabine, and Rein de Graaf}

\author[Marcel~de~Jeu]{Marcel~de~Jeu}

\begin{abstract}

We take a new and unifying approach toward polynomial and trigonometric approximation in
topological vector spaces used in analysis on $\Rn$. The idea is to show in considerable
generality that in such a space a module, which is generated over the polynomials or
trigonometric functions by some set, necessarily has the same closure as the module which is
generated by this same set, but now over the compactly supported smooth functions. The
particular properties of the ambient space or generating set are, to a large degree, irrelevant
for these subspaces to have equal closure. This translation---which goes in fact beyond
modules---allows us, by what is now essentially a straightforward check of a few properties, to
replace many classical results in various spaces by more general statements of a hitherto
unknown type. Even in the case of modules with one generator the resulting theorems on, e.g.,
completeness of polynomials are then significantly stronger than the classical statements. This
extra precision stems from the use of quasi-analytic methods (in several variables) rather than holomorphic methods, combined with the classification of quasi-analytic weights. In one dimension this classification, which then involves the logarithmic integral, states that two well-known families of weights are essentially equal.

As a side result
we also obtain an integral criterion for the determinacy of multidimensional measures which
is less stringent than the classical version.

The approach can be formulated for Lie groups and this interpretation then shows that many
classical approximation theorems are ``actually'' theorems on the unitary dual of $\Rn$, thus
inviting to a change of paradigm. In this interpretation polynomials correspond to the universal
enveloping algebra of $\Rn$ and trigonometric functions correspond to the group algebra.

It should be emphasized that the point of view, combined with the use of quasi-analytic methods,
yields a rather general and precise ready-to-use tool, which can very easily be applied in new
situations of interest which are not covered by this paper.

\end{abstract}

\date{}




\address{M.F.E.~de~Jeu\\
         Mathematical~Institute\\
         Leiden~University\\
         P.O.~Box 9512\\
         2300~RA~Leiden\\
         The~Netherlands}


\email{mdejeu@math.leidenuniv.nl}

\urladdr{}

\maketitle

\vfill
\eject
\tableofcontents

\vfill
\eject
\section{Introduction and overview}\label{sec:introduction}

In this paper we develop a method to approach a certain class of approximation problems in an
arbitrary number of variables. This class includes many instances of polynomial and
trigonometric approximation as well as ``mixed'' versions, the meaning of which will be
explained below. The present section is intended to motivate this method and contains an
explanation and overview of the main ideas and results.

The reader whose main interest lies in a possible application in a concrete situation is advised
to read the present section, Section \ref{subsec:defprel} up to and including Definition
\ref{def:qaweights}, Definition \ref{def:admissiblespaces}, Theorem~\ref{thm:admissiblespaces},
and then proceed to Sections~\ref{sec:general} and~\ref{sec:applications}, reading other parts
if and when necessary.

As to the motivation of the method, let us concentrate on polynomials for the moment and make
some general observations on the existing literature on polynomial approximation.

As a starting point, consider the following vector spaces and subspaces:
\begin{enumerate}
\item $C(K)$---the continuous functions on a compact subset $K$ of $\Rn$---with the
polynomials $\poly$ in $n$ variables as subspace;
\item $C^\infty(\Omega)$---the smooth functions on an open subset $\Omega$ of
$\Rn$---with $\poly$ as subspace;
\item $\distOmega$---the distributions on an open subset $\Omega$ of $\Rn$---with
$\poly$ as subspace;
\item $L_2(\Rn,dx)$ with $\poly\,\exp(-\Vert x\Vert^2)$ as subspace.
\end{enumerate}
When the ambient spaces are supplied with the usual topologies (the weak topology in the case of
distributions), each indicated subspace is sequentially dense, as is asserted by classical
theorems. Note that the subspaces are all modules over $\poly$, even though the ambient spaces
need not be. More precisely: the subspaces in question are cyclic modules over $\poly$,
generated as they are by either the constant $1$ or the Gaussian $\exp(-\Vert x\Vert^2)$.
Evidently, these classical density results can---although this is not the classical
formulation---all be rephrased as ascertaining the density of certain cyclic modules over
$\poly$. It is not difficult to find more examples of this nature. Notwithstanding this formal
similarity, the classical proofs of such results can differ substantially from case to case.

Continuing in this terminology, we furthermore note that most results on the closure of a module
over $\poly$ are positive results for the
\emph{density} of the module under consideration. Methods of a general nature to quickly
describe a (proper) closure do not appear to be known. As an example, the space $\poly\cdot
x\exp(-x^4)$ is not dense in the Schwartz space $\schw(\R)$ on the real line since any function
in its closure vanishes at $0$. Its closure consists in fact precisely of all Schwartz functions
vanishing at $0$, but a quick method yielding this description does not appear to be available.

Also, not much seems to be known about
\emph{non}cyclic modules. As a variation on the previous example, since the generators of
$\poly\cdot x\exp(-x^4)+ \poly\cdot (x-1)\exp(-x^2)$ do not have any common zero, there is now
no obvious obstruction for this space to be dense in $\schw(\R)$. It
\emph{is} dense, but again an appropriate method to conclude this quickly appears to be absent.

In addition, we note that evidently much more is known in one variable than in
several. The analysis in one variable is sometimes very refined, but this
situation is not at all paralleled in higher dimensions.

Last but not least, one is tempted to say that the situations in which results are known about the closure
of a module over $\poly$ have the character of a list. If one
encounters a situation which seems to be different from the known cases---e.g.,
the question whether $\poly\cdot\exp(-\Vert x\Vert^2 +\sin
\sqrt{\Vert x\Vert^2+1})$ is dense in $\schw(\Rn)$, which it in fact is---then the existing literature is sometimes of little help. Of course, existing proofs can often be adapted, but it may require some thought to actually do this in sufficient detail.

The above observations motivate the present paper. We develop a \emph{method} to approach the
problem of describing the closure of modules over the polynomials in a large class of situations
in topological vector spaces used in analysis on $\Rn$. This method, which extends naturally to
trigonometric approximation and ``mixed'' versions (to be explained below), is uniform in nature
and is valid in arbitrary dimension with, so it appears, quite some precision. It has a wide
scope of application and is easy to apply. We will establish the method in this paper and
give some applications. In addition, we will propose an interpretation in terms of Lie groups and
indicate the possibility of generalization from $\Rn$ to subsets of the unitary dual of certain
classes of Lie groups.

Let us explain the method. The observation which lies at the basis of the
approach is the following characterization of continuous maps:

\begin{lemma}\label{lem:fundamentallemma}
If $X$ and $Y$ are topological spaces, then a map $\phi:X\mapsto Y$ is
continuous if and only if, for all subsets $A$ and $B$ of $X$ having equal
closure in $X$, the images $\phi(A)$ and $\phi(B)$ again have equal closure in
$Y$.
\end{lemma}

\begin{proof}
If $\phi$ is continuous and $\overline A=\overline B$, then $\phi(A)\subset
\phi(\overline{A})=\phi(\overline{B})\subset
\overline{\phi(B)}$, the latter inclusion holding by continuity. Hence $\overline{\phi(A)}
\subset\overline{\phi(B)}$. By symmetry we conclude
$\overline{\phi(A)}=\overline{\phi(B)}$.
Conversely, let $\phi$ preserve the property of having equal closure. Since for
all $A\subset X$, $\overline A$ and $A$ have equal closure in $X$, we have
$\phi(\overline A)\subset\overline{\phi(\overline A)}=\overline{\phi(A)}$ for
all $A\subset X$, showing that $\phi$ is continuous.
\end{proof}

With this result we will prove a uniform ``translation'' theorem along the
following lines. We will establish topological vector spaces $\opw$, depending
on weights $w$ (to be discussed later on) and such that:
\begin{enumerate}
\item each $\opw$ consists of smooth functions on $\Rn$;
\item $\poly$ and the compactly supported smooth functions $\test$ are
subspaces of each $\opw$;
\item $\overline\poly=\overline\test$ in each $\opw$;
\item it is common to have continuous maps from a space $\opw$
to a topological vector space $V$ of the form $f\mapsto f\cdot g$ ($f\in\opw$),
where the dot indicates the appropriate type of action on a generator $g\in V$.
\end{enumerate}
The weight $w$ under (4) has to be chosen judiciously, taking the properties of $g$ and the
topology of $V$ into account. Such an appropriate choice need not always be possible, but in
those cases where it \emph{is} possible the previous lemma asserts that the property of having
equal closure is preserved, i.e., $\overline{\poly\cdot g}=\overline{\test\cdot g}$ in $V$. If
$V$ is in addition locally convex, and if $\gen\subset V$ is a set of generators to each of
which the above principle applies, then it is easy to see that
\begin{equation}\label{eq:translation}
\overline{\sum_{g\in\gen}\poly\cdot g}=\overline{\sum_{g\in\gen}\test\cdot g}
\end{equation}
as a consequence of the Hahn--Banach lemma.

It will become apparent in the subsequent sections that the equality
$\overline\poly=\overline\test$ in each $\opw$ is a nontrivial matter, which rests on a
combination of abstract vector-valued integration, Fourier analysis, and quasi-analytic classes
in several variables. The point of the above translation is therefore that this nontrivial
analysis for the spaces $\opw$ is thus put to \emph{immediate} use in ``any'' other situation to
validate a reformulation. Indeed, the problem of describing the closure of a module over $\poly$
has, according to
\eqref{eq:translation}, now been translated to the apparently equivalent problem of describing the
closure of the corresponding module over $\test$. The latter formulation of the
problem is in general much better accessible, since the properties of $\test$ in
relations to various spaces and their duals are usually well-catalogued. As an
example, it is immediate from the results in this paper that
$\poly\cdot\exp(-\sqrt{\Vert x\Vert^2+1}/\log \sqrt{\Vert x\Vert^2 +2})$ has the
same closure in $L_p(\Rn,dx)$ ($1\leq p<\infty$) as $\test\cdot
\exp(-\sqrt{\Vert x\Vert^2+1}/\log \sqrt{\Vert x\Vert^2 +2})$. Since the latter
space is simply $\test$, which is dense by general principles, the former space is therefore
also dense. Such an immediate translation is valid in all examples so far and the reader may
verify that this simplifies matters considerably, sometimes even reducing the situation to a
triviality.

Aside, let us note that one cannot expect to use the above method to prove that
\eqref{eq:translation} holds in all possible situations, since this is simply not true.
Examples where \eqref{eq:translation} is a meaningful but false statement are
provided by some of the spaces of the classical continuous Bernstein problem,
and by some $L_2$-spaces associated with nondeterminate moment problems, both
with $\gen=\{1\}$. There are also situations (the two aforementioned classes of
spaces again provide examples for this with $\gen=\{1\}$), where
$\eqref{eq:translation}$ is meaningful and holds, but where our method does not
allow us to conclude its validity. One could argue that these cases are subtle
or at least more delicate than many practical situations one is likely to
encounter, but this is admittedly a matter of taste. At any rate, we hope to
convince the reader that the above circle of ideas has a rather large scope of
application and that therefore, when used with appropriate caution,
\eqref{eq:translation} can be regarded as a rule of thumb with theoretical foundation.

Returning to \eqref{eq:translation} we note from an inspection of its ``proof'' that an extension
of \eqref{eq:translation} presents itself, since the only relevant fact for the validity of the ``proof'' is that $\poly$
is contained in the spaces $\opw$ with the same closure as $\test$. Any other
subspace with this property will do just as well. As it turns out, the span of
$\exp i(\lambda,x)$ for $\lambda\in\Rn$---we denote this span by $\galg$---also
has these properties together with $\poly$; in fact the span of the exponentials
corresponding to the spectral parameters in a somewhere dense subset of $\Rn$
(i.e., a subset such that its closure has nonempty interior) is already dense
in $\opw$. We thus obtain variations of
\eqref{eq:translation} when replacing $\poly$ by $\galg$
(or such a proper subspace of $\galg$), for all $g$ in a subset of $\gen$ of our choice. These
are the ``mixed'' versions of polynomial and trigonometric approximation alluded to above. One
obtains a whole family of (in general) different subspaces, all with equal closure. This type of
general result appears to be new.

It is also possible, and as we will see in fact natural, to relax the conditions on the weight
$w$ and develop versions of the spaces $\opw$ which contain $\galg$ but not $\poly$, still
satisfying the crucial property $\overline{\galg}=\overline\test$. This evidently yields results
on trigonometric approximation, valid in a wider setting than the polynomial context.

We will now discuss these spaces $\opw$ and their subspaces. A moment's thought shows that, if
one aims at a large scope of application, then the spaces $\opw$ should be small, implying that
$\poly$ and $\test$ should both be dense in the polynomial case and that $\galg$ and $\test$
should both be dense in the trigonometric case. The spaces should consist of smooth functions to
make the relevant maps sending $f$ to $f\cdot g$ well-defined in as many cases as possible, and
finally the topology should be fairly strong to enhance the continuity of these maps.

It turns out that such spaces can already be found in the existing literature, but their
potential use seems to have gone unnoticed. They have been introduced and studied in
\cite{Zapata,Nachbin2}, the latter reference aiming at simplifications of the former. We will
analyze only a special case of these spaces, which is nevertheless sufficient to obtain
significant extensions of the results in [loc.cit.] in the general situation (see Section
\ref{subsec:Bernsteinspaces}).

These spaces $\opw$ can be described as corresponding to the multidimensional differentiable
analogue of the classical one-dimensional continuous Bernstein problem. For the sake of
simplicity, let us restrict ourselves to the real line and assume that we are given a weight
$w:\R\mapsto\R_{\geq 0}$. We do not need any regularity hypothesis on $w$, but we do assume that
$w$ is bounded. The associated space $\opw$ then consists of all $f\in C^\infty(\R)$ such that
$\lim_{\vert x\vert\rightarrow\infty}f^{(m)}(x)\cdot w(x)=0$ for $m=0,1,2,\ldots.$ The seminorms
$p_m$ defined by $p_m(f)=\Vert f^{(m)}\cdot w\Vert_\infty$ ($m=0,1,2,\ldots$) turn $\opw$ into a
locally convex (not necessarily Hausdorff) topological vector space. It is easy to show that
$\test$ is dense in $\opw$ for all $w$. The subspace $\galg$ is not considered in
\cite{Zapata,Nachbin2}, but we show that $\galg$ is dense if $\lim_{\vert
x\vert\rightarrow\infty}w(x)=0$, so that we then have $\overline\galg=\overline\test$ as
required. If $w$ is rapidly decreasing then $\poly\subset\opw$ and one may ask what additional
conditions on $w$ are sufficient to ensure that $\poly$ is dense in $\opw$, i.e., that
$\overline\poly=\overline\test$ ($=\overline\galg$). It is known \cite{Zapata,Nachbin2} that
this holds whenever $w$ is a so-called quasi-analytic weight, i.e., when
\begin{equation}\label{eq:introseries}
\sum_{m=1}^\infty \Vert x^m\cdot w(x)\Vert_\infty^{-1/m}=\infty.
\end{equation}
Here the term ``quasi-analytic'' refers to an application of the Denjoy--Carleman theorem in the
proofs which is validated by this divergence.

We will give an independent proof of the density of $\poly$ in $\opw$ for quasi-analytic
weights. Whereas the approach in \cite{Zapata,Nachbin2} in the multidimensional case consists of
a reduction to the one-dimensional case, we use a direct approach in arbitrary dimension which
brings out the natural connection between $\poly$ and $\galg$. This direct method not only
yields an extra trigonometric approximation result (namely that the span of the exponentials
corresponding to spectral parameters in a somewhere dense subset of $\Rn$ is already dense in
$\opw$ if $w$ is quasi-analytic), but it also shows what the natural interpretation in terms of
Lie groups is, thus indicating how the setup can be generalized to certain classes of Lie
groups.

From the properties of the spaces $\opw$ as mentioned above it is clearly desirable to have more
information on quasi-analytic weights. Although they have been studied in abstract in the
literature as above, nontrivial examples of such weights appear to be sparse, if not absent.
Such weights are necessarily rapidly decreasing and it is easily verified that weights on $\R$
of type $C\exp(-\epsilon\vert x\vert)$ and their minorants are quasi-analytic. For such weights
however the above results do not provide new insights, since one can in this case prove the
same---in fact even stronger---closure results by holomorphic methods rather than quasi-analytic
ones.

In this paper we therefore make a systematic study of weights, and of quasi-analytic weights in
particular. We are then not only able to give concrete nontrivial examples of quasi-analytic
weights (see Section~\ref{subsec:examplesqa}), but we will actually classify these weights in
arbitrary dimension. On the real line the equivalent characterization is as follows: a weight on
$\R$ satisfies
\eqref{eq:introseries} if and only if there exists a continuous even strictly positive weight $\wtilde$
majorizing $w$, such that $s\mapsto -\log\wtilde(e^s)$ is convex on
$(-\infty,\infty)$ and such that
\begin{equation*}
\int_0^\infty \frac{\log\wtilde(t)}{1+t^2}=-\infty.
\end{equation*}
This classification has two interesting aspects. First, it has long been known that such
weights $\wtilde$ have good approximation properties in a number of situations, but it has not
previously been noticed that they are essentially the same as the quasi-analytic weights with
defining property \eqref{eq:introseries} as in \cite{Nachbin1}. Second, since in particular
such a weight $\wtilde$ is itself a quasi-analytic weight, the results in the present paper can
be regarded as an \emph{explanation} of these good classical approximation properties in view of
the possibility to push the strong results on equal closure in the corresponding space
$\opwtilde$ forward to a variety of other situations by Lemma~\ref{lem:fundamentallemma}.

This classification of quasi-analytic weights has another application. As the reader may verify
from \eqref{eq:introseries}, a probability measure $\mu$ on $\R$ is determinate by Carleman's
theorem whenever there exists a strictly positive quasi-analytic weight $w$ such that $w^{-1}\in
L_1(\R,\mu)$. Since the latter condition can be tested by an integral, the conclusion of
determinacy can therefore be reached directly from the measure, which is evidently simpler than
a computation of the moment sequence itself. A similar result for determinate measures holds in
arbitrary dimension and, in view of the explicit examples of quasi-analytic weights in
Subsection~\ref{subsec:examplesqa}, this establishes a new class of multidimensional determinate
measures. In Subsection \ref{subsec:locallycompact} we will conclude the determinacy from a density
result, but it is also possible to prove this determinacy in a more direct fashion as an
application of an extended Carleman theorem in arbitrary dimension. We will report on this
separately in \cite{deJeu}.

Returning to the spaces $\opw$ we mention that the analogues of the closure results on the real
line hold in arbitrary dimension. As to the proofs, the main idea for the proof of the density
of $\poly$ in $\opw$ for quasi-analytic weights is standard, but the underlying technique
involves Fourier analysis, quasi-analytic classes in several variables and vector-valued
integration in a combination which seems to be new. Furthermore, it is interesting to note that,
whereas most applications of the Fourier transform in approximation theory use its injectivity,
we, on the contrary, use a result on its \emph{range}. This is not a coincidence and this
fundamental role of harmonic analysis becomes more transparent when reformulating the method for
Lie groups.

Turning to the applications, let us first make some comments of a general nature.

Let us assume that $V$ is a topological vector space, that $g\in V$, and that we wish to conclude
that $\overline{\poly\cdot g}=\overline{\test\cdot g}$ as above by finding a quasi-analytic
weight $w$ such that the map $f\mapsto f\cdot g$ is continuous from $\opw$ into $V$. The first
requirement is obviously that $\opw\cdot g\subset V$. Since $\opw$ contains functions which
together with their derivatives grow at infinity at the rate of $w^{-1}$, this inclusion can,
roughly speaking, be expected to hold when $g$ and its relevant derivatives decay at infinity at
the same rate as $w$. In the same fashion, in the case of trigonometric approximation one
intuitively expects $\opw\cdot g\subset V$ to hold for some weight $w$ tending to zero at
infinity when $g$ tends to zero at infinity, together with its relevant derivatives. These
imprecise statements are of course only useful as an intuitive orientation in practical
situations. An interesting case occurs when $V$ is a $C^\infty(\Rn)$-module. Then the first
requirement $\opw\cdot g\subset V$ is trivially satisfied for all $g\in V$ and all $w$ and
chances increase that
\eqref{eq:translation} holds for all generators.

Let us now consider the second requirement, i.e., the continuity of the map $f\mapsto f\cdot g$
from $\opw$ into $V$, assuming that it is well-defined. Surprisingly enough, one can prove in
considerable generality that the hypothesis of continuity is not necessary to conclude
\eqref{eq:translation}. This interesting phenomenon occurs, e.g., for all $V$ that can be regarded as a subspace of a
space of distributions, under mild conditions on the topology of $V$ which are satisfied by many
common spaces. All the spaces we have considered in this section are within the scope of this
result. The underlying reason for the redundancy of the hypothesis of continuity is the fact
that in this situation various spaces related to $V$ are Fr\'echet spaces; since it is in
addition possible to ``improve'' the spaces $\opw$ into Fr\'echet spaces, the closed graph
theorem comes within range.

This result on the redundancy of the hypothesis of continuity for a large class of spaces
simplifies the application of our method. Indeed, the only hypothesis to be verified in all such
situations is now an inclusion of the type $\opw\cdot g\subset V$, which involves no more than
the definition of $V$ \emph{as a set}. The precise nature of $V$ and $g$ are to a large degree
immaterial. For the detailed formulation of the very general results in this direction we refer
to Section~\ref{sec:general}. Here we give the following typical example, which is a special
case of Corollary~\ref{cor:nogrowth}:

\begin{theorem}\label{thm:automatic}
Let $X\subset\Rn$ be nonempty and let $V$ be an LF-space of functions on $X$
with a defining sequence $V^1\subset V^2\subset\ldots$. Suppose that convergence
of a sequence
in $V$ implies pointwise convergence on $X$ and that the spaces $V^l$ are all
$C^\infty(\Rn)$-modules.

Let $\gen\subset V$ be a nonempty set of generators and for $g\in\gen$ let
$L_g$ denote either $\test$, $\poly$, or $\Span_\C \{\, \exp i(\lambda,x)\mid
\lambda\in E_g\,\}$ where $E_g$ is a somewhere dense subset of $\Rn$ or, more generally, a uniqueness set for holomorphic functions on $\Cn$.

Then the closure in $V$ of the subspace $L=\sum_{g\in\gen}L_g\cdot g$ is
independent of the choice of the $L_g$.
\end{theorem}

The results of this type provide a theoretical foundation for the aforementioned suggestion of
regarding \eqref{eq:translation} as a rule of thumb.

Closely related to the previous paragraphs and relevant for applications is the notion of an
admissible space for an element $g$ in a topological vector space $V$, which we will introduce
for $V$ in the category of spaces for which the continuity hypothesis is redundant. Roughly
speaking, an admissible space for $g$ is any space $L$ of functions of a common nature for which
there exists a weight $w$ such that $\overline L=\overline\test$ in $\opw$ and such that
$\opw\cdot g\subset V$. The size of the collection of admissible spaces for $g$ determines the
amount of freedom in varying the spaces acting on $g$ in
\eqref{eq:translation}. Stated intuitively (and disregarding the possibly
relevant derivatives), the connection is that this freedom increases together
with the maximum rate of growth of smooth functions $f$ such that $f\cdot g\in
V$. If $V$ is a $C^\infty(\Rn)$-module, then there is evidently maximal freedom for arbitrary $g\in V$ and one obtains results as
Theorem~\ref{thm:automatic}, where the precise nature of $V$ is to a large
degree irrelevant.

We now turn to the specific applications in various common spaces as they are encountered in
analysis on $\Rn$. In all but a few examples we are able to give necessary and sufficient
conditions for modules to be dense. The annihilator of nondense modules can be described and in
a number of cases this is sufficient to determine explicitly the closure of the module under
consideration. This being said, let us indicate very briefly some of the additional features,
leaving most of the material to Section~\ref{sec:applications}.

The first application is $\schw$. We show that a closed translation-invariant
subspace of $\schw$ is dense if and only if the Fourier transforms of its
elements have no common zero. Consideration of mixed modules yields (under
Fourier transform) multidimensional results of a type that has previously been
considered on the real line by Mandelbrojt \cite{MandelbrojtRice}.

The second application treats various spaces of test functions and their duals.
Classical results such as the sequential density of $\poly=\poly\cdot 1$ in
$\distOmega$ are now an immediate consequence of the sequential density of
$\test=\test\cdot 1$.

The third application is concerned with the general type of space for the
differentiable multidimensional Bernstein problem. We extend the results in
\cite{Zapata,Nachbin2}. A noteworthy particular case is the classical continuous
Bernstein problem on the real line; the classical sufficient condition for the
density of the polynomials in that situation is seen to be part of a more
general picture as described in the theorem in this paper.

Our fourth application is concerned with various spaces associated with a locally compact subset
$X$ of $\Rn$, including $L_p$-spaces and spaces of continuous functions. Specialization of these
results yields, e.g., the following: if $\mu$ is a Borel measure on $X$ (or the completion of a
Borel measure) such that $w^{-1}\in L_1(X,\mu)$ for a measurable quasi-analytic weight $w$ which
is strictly positive on $X$, then $\poly$ is dense in $L_p(X,\mu)$ for all $1\leq p<\infty$. The
special case where $X=\Rn$, $p=2$, and $w(x)=\exp (-\epsilon\Vert x\Vert)$ for some $\epsilon>0$
yields the classical Hamburger theorem. However, we now also have that for general $X$ the
polynomials are dense for all finite $p$ if $w^{-1}\in L_1(X,\mu)$ for a measurable weight $w$
which is strictly positive on $X$ and for large $\Vert x\Vert$ equal to, e.g.,
$\exp(-\epsilon\Vert x\Vert/\log a\Vert x\Vert)$ for some $\epsilon,\,a>0$. This is a much less
stringent condition than Hamburger's, and it will become apparent from the examples of
quasi-analytic weights that still more lenient explicit conditions can be given. As alluded to
above, for $X=\Rn$ this density in $L_p$-spaces in turn implies that such measures are
determined by their moments. The resulting criterion for determinacy also extends the known case
and is stated as our fifth and final application.

We also propose an interpretation of the spaces $\opw$ in terms of Lie groups. In this
interpretation, all the results under consideration in this paper are results situated on the
unitary dual of $\Rn$; this implies a change of paradigm. We will argue that the space $\galg$
describes the action of the group algebra of $\Rn$ on the smooth vectors in the various
irreducible unitary representations of $\Rn$. Indeed, in this interpretation an element
$x\in\Rn$ acts on the smooth vectors in the irreducible representation with spectral parameter
$\lambda$ as multiplication by $\exp i(\lambda,x)$. Correspondingly, the function
$\lambda\mapsto \exp i(\lambda,x)$ with $x\in\Rn$ fixed should at an abstract level be thought
of as assigning to each irreducible unitary representation an endomorphism of its space of
smooth vectors. Similarly, the action of an element of the universal enveloping algebra of $\Rn$
on all these spaces of smooth vectors is described globally on the unitary dual as
multiplication by a polynomial, thus identifying $\poly$. Note that in this interpretation the
canonical choice for the notation of the variable in many classical approximation results is
``actually'' $\lambda$ and not $x$. For further details we refer to
Section~\ref{sec:interpretation}, where we also indicate how the setup might be generalized to
subsets of the unitary dual of certain classes of Lie groups.

This paper is organized as follows:

In Section~\ref{sec:basics} we establish the basic notations and recall some
auxiliary results.

Section~\ref{sec:weights} is primarily concerned with a systematic study of weights, and notably
quasi-analytic weights in several variables. The definition of these weights is shown to be
related to frames in an essential way. We classify such weights in a practical fashion and give
explicit examples (including notably radial weights) in terms of elementary functions. It should
be emphasized that it is precisely this classification in terms of divergent integrals from
which our method lends its practical precision. We also show that, e.g., quasi-analytic weights
can often be assumed to be of class $C^\infty$, even though no regularity at all was imposed
from the outset. Although this possibility of regularization is not essential for the remainder
of the paper, it seems less than satisfactory to let it go unnoticed. This section also contains
the definition of the spaces $\opw$. An additional and important topic in this section is the
``ubiquitous'' possibility of finding spaces $\opw$ that are Fr\'echet.

Section~\ref{sec:density} contains the proof of the density results in $\opw$.
As was already mentioned above, there is---as far as the polynomials are
concerned---an overlap of our results with \cite{Zapata,Nachbin2}. Since in
[loc.\ cit.] $\galg$ is not considered in its natural role and since our method of
proof appears to be more direct than the reduction to the continuous case in one
dimension as in [loc.\ cit.], and at the same time indicates what the possible
generalization to Lie groups is, we feel that this section provides a more
informative analysis of the spaces $\opw$ than [loc.\ cit.]. This section also
introduces the notion of admissible spaces for the various types of weights,
i.e., the various subspaces of $\opw$ of a common nature that are, depending on $w$, dense in
$\opw$ together with $\test$.

Section~\ref{sec:general} may be regarded as the formalization of our method. It contains the
general theorem on subspaces with equal closure, as well as three results on the redundancy of
the hypothesis of continuity. The section concludes with the discussion of a partial filtration
of topological vector spaces that can be related to the various types of weights. This leads
naturally to the aforementioned notion of admissible spaces for elements of topological vector
spaces. The partial filtration is also relevant for the regularity of associated Fourier
transforms.

We treat the applications in Section~\ref{sec:applications} and conclude with
the interpretation and possible generalizations to subsets of the unitary dual
of certain classes of Lie groups in Section~\ref{sec:interpretation}. This
generalization is a feasible program for at least the unitary dual of $\Tn$,
yielding results situated on $\Zn$. We make some brief comments on this case.

In Appendix~\ref{app:closurecont} we elaborate on Lemma
\ref{lem:fundamentallemma}, also considering the somewhat more delicate case of
sequential closure in topological vector spaces.

Finally, in Appendix~\ref{app:DCtheorem} we prove the theorem on quasi-analytic
classes in several variables that lies at the heart of the density of $\poly$ in
$\opw$ for a quasi-analytic weight $w$.

\section*{Acknowledgments}
It is a pleasure to thank Christian Berg, Jaap Korevaar, and Jan Wiegerinck for
helpful remarks. During the preparation of this paper the author was partially
supported by a PIONIER grant of the Netherlands Organisation for Scientific Research (NWO).

\section{Notations and auxiliary results}\label{sec:basics}

Throughout this paper we will be concerned with $\Rn$; $n$ is fixed unless
explicitly stated otherwise. We let $\{e_1,\ldots,e_n\}$ be the standard basis
of $\Rn$. In the case $n=1$ we write $e=e_1=1$. We supply $\Rn$ with the standard
inner product $(\,.\,,\,.\,)$  and corresponding norm $\Vert\,.\,\Vert$. The inner product is extended bilinearly to $\Cn$. The transpose of a
linear transformation of $\Rn$ is defined with respect to this standard inner
product.

$\poly$ will denote the complex-valued polynomials on $\Rn$; $\poly_d$ consists of the
polynomials of degree at most $d$. If $\lambda\in\Cn$, then $e_{\lambda}:\Rn\mapsto\C$ is given
by $e_{\lambda}(x)=\exp (\lambda,x)$ ($x\in\Rn$). The complex span of the $e_{i\lambda}$ for
$\lambda\in\Rn$ is denoted by $\galg$. We will refer to $\lambda$ as the spectral parameter for
$e_{i\lambda}$, suppressing the factor $i$.

If $A\subset\Rn$, then $\charfunc_A$ denotes the characteristic (indicator) function of $A$.

All topological vector spaces in this paper are assumed to be complex, with the
only and obvious exception of $\Rn$. The pairing between a topological vector
space $V$ and its continuous dual $V^\prime$ is denoted as in $\langle
v^\prime,v\rangle$.
In our convention a topological vector space need not be Hausdorff. We mention
explicitly that for a linear subspace $L$ of a locally convex space $V$, the
closure $\overline L$ of $L$ can be described in the well-known fashion as
$\overline{L}=\{\,v\in V\mid
\langle v^\prime,v\rangle=0\,\,\forall \,v^\prime\in L^\perp\}$, where the
annihilator $L^\perp$ is the subspace of $V^\prime$ consisting of all elements
vanishing on $L$. For this to be valid $V$ need not be Hausdorff (see \cite{Treves} or
\cite[p.~II.45--46, p.~II.64]{BourbakiTVS}).
A Fr\'echet space is in our terminology a locally convex metrizable complete
topological vector space. For LF-spaces, i.e., inductive limits of Fr\'echet
spaces, we refer to \cite{Treves,Horvath2,RudinFA}.

If $X$ is a topological space and $A\subset X$ then the sequential closure
${\overline {A}}^s$ is the set of all $x\in X$ for which there exists a sequence
$\{a_m\}_{m=1}^\infty\subset A$ such that $\lim_{m\rightarrow\infty}a_m=x$.
Obviously, ${\overline {A}}^s\subset \overline A$ and if all points in $X$ have a
countable neighborhood base then equality holds for all $A\subset X$. If $X$ is
a topological vector space and $L$ is a subspace, then ${\overline L}^s$ is
evidently again a subspace.

If $\Omega$ is a nonempty subset of $\Rn$ and $m\in\{0,1,2,\dots,\infty\}$, we
let $\testmOmega$ denote all functions of class $C^m$ on $\Omega$, endowed with
the Fr\'echet topology of uniform convergence of all derivatives of order at
most $m$ on compact subsets of $\Omega$. The elements in $\testmOmega$ with
compact support form a subspace $\testcmOmega$ which is endowed with the usual
LF-topology. The dual spaces $\distmOmega$ (resp., $\compdistmOmega$) are
identified with the subspaces of $\distOmega$ consisting of distributions of
order at most $m$ (resp., the compactly supported distributions of order at most
$m$), as a result of the density (in fact sequential density) of $\testcOmega$
in $\testcmOmega$ (resp., the density of $\testcOmega$ in $\testmOmega$). Both
$\distmOmega$ and $\compdistmOmega$ are modules over $\testmOmega$. When
$\distmOmega$ and $\compdistmOmega$ carry their weak topologies, $\testcOmega$
is in both cases a sequentially dense subspace. We refer to
\cite{Treves,Horvath2} for details. We write $\test$ for $C_c^\infty(\Rn)$ and
$\dist$ for the distributions on $\Rn$.

The space $\schw$ of rapidly decreasing functions on $\Rn$ and its dual space of
tempered distributions $\tempdist$ are as usual. Any statement related to a
topology on the spaces $\testcmOmega$, $\testmOmega$, or $\schw$ always refers to
the usual topologies as just described.

We adhere to the usual conventions concerning multi-indices. Thus, if
$\alpha=(\alpha_1,\ldots,\alpha_n)\in\Nn$ (where $\N$ includes $0$), then $\parta$ denotes the
corresponding partial derivative of order $\vert\alpha\vert$. If necessary we
will employ self-evident notations as $\parta_x$ to indicate the variable on
which the operator acts. For $x=(x_1,\ldots,x_n)\in\Cn$ and $\alpha\in\Nn$, we
let $x^\alpha=x_1^{\alpha_1}\cdots x_n^{\alpha_n}$ and $\vert
x\vert^\alpha=\vert x_1\vert^{\alpha_1}\ldots\vert x_n\vert^{\alpha_n}$.

We use the following normalization of the Fourier transform
$\Fourier$:
\begin{equation*}
\Fourier(f)(\lambda)={\widehat f}(\lambda)=(2\pi)^{-n/2}\int_{\Rn} f(x)
e^{-i(\lambda,x)}\,dx\quad(\lambda\in\Rn,\,f\in L_1(\Rn,dx)).
\end{equation*}
For $f\in L_1(\Rn,dx)$ such that $\Fourier(f)\in L_1(\Rn,d\lambda)$ we then have the inversion formula
\begin{equation*}
f(x)=(2\pi)^{-n/2}\int_{\Rn} \Fourier(f)(\lambda)
e^{i(\lambda,x)}\,d\lambda\quad(\text{a.e.}).
\end{equation*}
$\Fourier$ extends to tempered distributions
by
\begin{equation*}
\langle \Fourier(\typ),f\rangle=\langle
\typ,\Fourier(f)\rangle\quad(T\in\tempdist,\,f\in\schw).
\end{equation*}

On occasion we will work in the extended nonnegative real numbers $[0,\infty]$
with the usual conventions, including $0\cdot\infty=0$ and $0^0=\infty^0=1$ .

Let $R>0$. We will encounter nonincreasing functions
$w:[R,\infty)\mapsto\R_{>0}$ such that $s\mapsto -\log w(e^s)$ is convex on
$[\log R,\infty)$. If $w\in C^{(2)}([R,\infty))$, then there exists $\rho\in
C^{(1)}([R,\infty))$ with $\rho,\,\rho^\prime\geq 0$ on $[R,\infty)$ and such
that
\begin{equation}\label{eq:convexrelation}
w(t)=w(R)\exp\left(-\int_R^t \frac{\rho(s)}{s}\,ds\right)\quad(t\geq R).
\end{equation}
Conversely, if $\rho\in C^{(1)}([R,\infty))$ is nonnegative and nondecreasing
and a constant $w(R)>0$ is given, then \eqref{eq:convexrelation} defines a
function $w$ with the aforementioned properties. This type of equivalence is
well-known (cf.~\cite{Mergelyan}).

\section{Weights and associated spaces $\opw$}\label{sec:weights}

This section contains the definition of weights, the associated function spaces
$\opw$, and a first investigation of the relevant properties. The emphasis is on
the classification of quasi-analytic weights. For these weights the criterion in
Definition~\ref{def:qaweights}, which is suitable for applications of
Denjoy--Carleman-type theorems, can be shown to be related to the more practical
(and at the same time familiar) criteria as in Theorem~\ref{thm:classonedim}.
This allows us to determine explicit examples of quasi-analytic weights in
Section~\ref{subsec:examplesqa}. Another main result is Theorem
\ref{thm:improvements}, where the ubiquity of complete spaces $\opw$ for all
types of weights under consideration is noted. A difference with most of the
literature, apart from the arbitrary number of variables, is that we impose no
regularity condition on weights.

\subsection{Definitions and preliminaries}\label{subsec:defprel}

An arbitrary function $w:\Rn\mapsto\R_{\geq 0}$ is called a
\emph{weight} if it is bounded. The set of all weights is denoted by $\weights$.
The support $\sup w$ of a weight $w$ is the set $\{\,x\in\Rn\mid w(x)\neq 0\,\}$.

Given a weight $w$, we define the vector space $\op_w$ by
\begin{equation*}
\opw=\{\,f\in C^\infty(\Rn)\mid\lim_{\Vert x\Vert\rightarrow\infty}\parta
f(x)\,\cdot\,w(x)=0\,\,\forall\alpha\in\Nn\,\}.
\end{equation*}
Then $\opw$ is a locally convex topological vector
space when supplied with the topology determined by the seminorms
\begin{equation*}
\pa(f)=\Vert\parta f\,\cdot\,w\Vert_\infty\quad(f\in\opw,\,\alpha\in\Nn).
\end{equation*}

\begin{remark}\label{rem:Hausdorff}
It is straightforward to check that $\opw$ is Hausdorff if and only if the
complement of the zero locus of $w$ is dense in $\Rn$. By the countability of
the defining family of seminorms, this condition on the zero locus of $w$ is
also the necessary and sufficient condition for metrizability of $\opw$.
\end{remark}

We note that if $w(x_0)\neq 0$ for some $x_0\in\Rn$, then convergence in $\opw$
implies convergence of all derivatives at $x_0$. In the same vein, if $w$ is a
weight such that $\inf_{x\in K}w(x)> 0$ for every nonempty compact subset $K$
of $\Rn$, then $w$ is said to be
\emph{separated from zero on compact sets}. Convergence in the associated space
$\opw$ then implies uniform convergence of all derivatives on all compact sets,
i.e., there is a continuous injection $\opw\hookrightarrow C^\infty(\Rn)$.

We distinguish several classes of weights; their interrelations are summarized
in Proposition~\ref{prop:invariance}. To start with, for $0\leq d<\infty$, we let
$\weights_d$ denote the set of all weights $w$ such that $\lim_{\Vert
x\Vert\rightarrow\infty}\Vert x\Vert^d w(x)=0$. The
\emph{rapidly decreasing weights} $\weights_\infty$ are defined as
$\weights_\infty=\bigcap_{d=0}^\infty
\weights_d$; evidently, $w\in\weights_\infty$ if and only if $\poly\subset\opw$. If
$w$ is a weight and if there exist $C\geq 0$ and $\epsilon>0$ such that
$w(x)\leq C\exp(-\epsilon\Vert x\Vert)$ for all $x\in\Rn$, then $w$ is \emph{a
holomorphic weight}. The nomenclature refers to the fact that the Fourier
transform of an element of $\opwdual$ is then a holomorphic function, as we will
see in Section~\ref{sec:density}. The weight itself might be highly irregular.
The set of holomorphic weights is denoted by $\weightshol$. For
$w\in\weightshol$, define $0<\tau_w\leq\infty$ as the supremum of the set of all
$\epsilon>0$ for which $w(x)\leq C_\epsilon\exp(-\epsilon\Vert x\Vert)$ for all
$x$ and some $C_\epsilon\geq 0$. Then $\tau_w$ determines a tubular open
neighborhood of $\Rn$ in $\Cn$ for which the corresponding $e_{i\lambda}$ are
still in $\opw$.

We now come to the definition of the main class of weights that will concern us:
the \emph{quasi-analytic weights}. Here again the nomenclature refers to a
property of Fourier transforms and not to the weight itself. The definition is
as follows:

With $v\in\Rn$ we associate $\tilde v:\Rn\mapsto\R$ defined by $\tilde
v(x)=(v,x)$. If $w$ is a weight and $a\geq 0$, define the nonnegative possibly
infinite quantity
\begin{equation*}
M_w(v,a)=\Vert \vert \tilde v\vert^a\cdot w\Vert_\infty.
\end{equation*}

\begin{definition}\label{def:qaweights}
If $\{v_1,\ldots,v_n\}$ is a basis of $\Rn$ such that
\begin{equation}\label{eq:divergentseries}
\sum_{m=1}^\infty M_w(v_j,m)^{-1/m}=\infty
\end{equation}
for all $1\leq j\leq n$, then $w$ is called \emph{quasi-analytic with respect to
$\{v_1,\ldots,v_n\}$}. A weight is \emph{standard quasi-analytic} if it is
quasi-analytic with respect to the standard basis of $\Rn$. A weight is
\emph{quasi-analytic} if it is quasi-analytic with respect to some basis. The set of
quasi-analytic weights is denoted by $\weightsqa$.
\end{definition}

We will see later on that $\weightsqa\subset\weights_\infty$. A small
computation shows that holomorphic weights are quasi-analytic with respect to
any basis. It is a nontrivial fact that a weight on $\Rn$ ($n\geq 2$) can be
quasi-analytic with respect to one basis and not with respect to any other,
apart from scaling of the original basis. The notion of quasi-analytic weights
is therefore linked to frames in an essential way; we return to this point in
Subsection~\ref{subsec:qaw}.

In one variable, a condition of type \eqref{eq:divergentseries} already occurs
in Bernstein's work on weighted polynomial approximation \cite{Bernstein}. It is
a recurring theme in the work of Nachbin (cf.~\cite{Nachbin1}), who in the case
of one variable appears to have introduced the terminology. In several variables
a proper subset of the quasi-analytic weights makes it appearance in
\cite{Zapata, Nachbin2}, although not under this name and not in the generality
or central role in which we will study them.

We start with some elementary facts, which we state as lemmas for the ease of
reference. The verification is left to the reader. The extended nonnegative
real numbers $[0,\infty]$ form the context in Lemma~\ref{lem:elementarytwo}.

\begin{lemma}\label{lem:elementaryone}
\begin{enumerate}
\item Let $\{a_m\}_{m=1}^\infty\subset\R_{>0}$ and suppose $\sum_{m=1}^\infty
a_m^{-1}=\infty$. If $c\geq 0$, then $\sum_{m=1}^\infty (c+ a_m)^{-1}=\infty$.
\item Let $\{a_m\}_{m=1}^\infty\subset\R_{\geq 0}$ be nonincreasing. If $k$ and $l$
are strictly positive integers, then $\sum_{m=1}^\infty a_{km}=\infty$ if and
only if $\sum_{m=1}^\infty a_{lm}=\infty$.
\end{enumerate}
\end{lemma}

\begin{lemma}\label{lem:elementarytwo} Let $S$ be a nonempty set and suppose $g,w:S\mapsto\R_{\geq 0}$.
\begin{enumerate}
\item If $a,\,b\geq 0$ then $\Vert g^{\lambda a + (1-\lambda)b}\cdot w\Vert_\infty\leq \Vert g^a\cdot w\Vert_\infty^\lambda\,\Vert g^b\cdot
w\Vert_\infty^{1-\lambda}$ for $0\leq \lambda\leq 1$.
\item If $\Vert w\Vert_\infty<\infty$ and $\Vert g^a\cdot w\Vert_\infty=\infty$ for some $a>0$, then $\Vert g^b\cdot w\Vert_\infty=\infty$ for all $b\geq a$.
\item If $\Vert w\Vert_\infty\leq 1$, then $\Vert g^a\cdot w\Vert_\infty^{1/a}\leq \Vert g^b\cdot w\Vert_\infty^{1/b}$ whenever $0<a\leq b$.
\end{enumerate}
\end{lemma}

\begin{remark}\label{rem:convexity}
The first part of Lemma~\ref{lem:elementarytwo} implies that
$\{M_w(v,m)\}_{m=0}^\infty$ is an (obviously nonnegative) logarithmically
convex sequence for all weights $w$ and all $v\in\Rn$. We will use this
repeatedly.
\end{remark}

Given the fact that a weight can be quasi-analytic with respect to essentially
one basis, it is natural to consider the action of the general linear group. For
our purposes it is actually better to consider the action of the affine group
$\aff$, where the action on functions is defined by $Af(x)=f(A^{-1}x)$ for
$f:\Rn\mapsto\C$ and $A\in\aff$.

\begin{lemma}\label{lem:invarianceweightsqa}
Let $w$ be a quasi-analytic weight with respect to $\{v_1,\ldots,v_n\}$. Then
$w$ is rapidly decreasing. If $A\in\aff$ is linear, then $Aw$ is
quasi-analytic with respect to $\{(A^{-1})^t v_1,\ldots,(A^{-1})^t v_n\}$.
If $T\in\aff$ is a translation, then $T w$ is quasi-analytic with respect to
$\{v_1,\ldots,v_n\}$.
\end{lemma}

\begin{proof}
To see that $w$ is rapidly decreasing, first note that---in view of
\eqref{eq:divergentseries}---the second part of Lemma~\ref{lem:elementarytwo} implies that
$\Vert\tilde v_j^m\cdot w\Vert_\infty<\infty$ for all $1\leq j\leq n$ and all
$m\geq 1$. Trivially this also holds when $m=0$. Next choose a disjoint
decomposition $\Rn=\bigcup_{j=1}^n S_j$ such that for $j=1,\ldots,n$ we have
$\max_{i=1,\ldots,n}\vert
\tilde v_i(x)\vert=\vert\tilde v_j(x)\vert$ for all $x\in
S_j$. Then for $\alpha\in\Nn$ we have
\begin{equation*}
\Vert \prod_{i=1}^n
\tilde v_i^{\alpha_i}\cdot w\Vert_\infty\leq \sum_{j=1}^n
\Vert \charfunc_{S_j}\prod_{i=1}^n \tilde v_i^{\alpha_i}\cdot w\Vert_\infty\leq \sum_{j=1}^n \Vert \tilde v_j^{\vert\alpha\vert}\cdot
w\Vert_\infty<\infty,
\end{equation*}
by the finiteness noted above. Hence $w$ is rapidly decreasing.

The statement concerning the action of $A$ is immediate from the definitions.

Turning to the translation part, let $T=T_y$ where $T_y$ is the translation
$T_y(x)=x+y$ over $y$. We may assume $\Vert w\Vert_\infty=1$. Using the third
part of Lemma~\ref{lem:elementarytwo} in the next to last step we have
\begin{align*}
\Vert \widetilde v_j^m\cdot T_y w \Vert_\infty&=\Vert (v_j,x)^m\cdot w(x-y)\Vert_\infty\\
&=\Vert(v_j,x+y)^m\cdot w(x)\Vert_\infty\\&\leq \sum_{i=0}^m \binom {m}{i}
\Vert (v_j,x)^i\cdot w(x)\Vert_\infty\,\vert
(v_j,y)\vert^{m-i}\\ &=
 \vert (v_j,y)\vert^{m} + \sum_{i=1}^m \binom {m}{i}
\Vert (v_j,x)^i\cdot w(x)\Vert_\infty\, \vert (v_j,y)\vert^{m-i}\\
&= \vert (v_j,y)\vert^{m} +
 \sum_{i=1}^m \binom {m}{i}  \left\{\Vert (v_j,x)^i\cdot
 w(x)\Vert_\infty^{1/i}\right\}^i \vert(v_j,y)\vert^{m-i}\\ &\leq\vert (v_j,y)\vert^{m} +
 \sum_{i=1}^m \binom {m}{i} \left\{\Vert (v_j,x)^m\cdot
 w(x)\Vert_\infty^{1/m}\right\}^i\vert (v_j,y)\vert^{m-i}\\
 &=\left\{\vert (v_j,y)\vert+\Vert (v_j,x)^m\cdot w(x)\Vert^{1/m} \right\}^m, 
\end{align*}
so that we have the estimate $\Vert\tilde
v_j^m
\cdot T_y w\Vert_\infty^{-1/m}\geq (\vert (v_j,y)\vert +\Vert
\tilde v_j^m\cdot w\Vert_\infty^{1/m})^{-1}$ for all $m\geq 1$.

Now if $\Vert \tilde v_j^m\cdot w\Vert_\infty=0$
for some $m\geq 1$, then $\Vert \tilde v_j^m\cdot w\Vert_\infty=0$ for all
$m\geq 1$ and it follows trivially from the estimate that $\sum_{m=1}^\infty
M_{T_yw}(v_j,m)^{-1/m}=\infty$. In the case $\Vert
\tilde v_j^m\cdot w\Vert_\infty\neq 0$ for all $m\geq 1$, the necessary divergence is
implied by the estimate and the first part of Lemma~\ref{lem:elementaryone}.
\end{proof}

\begin{proposition}\label{prop:invariance}
We have
$\weightshol\subset\weightsqa\subset\weights_\infty\subset\weights_{d}\subset\weights_{d^\prime}\subset\weights_0\subset\weights$
for $d\geq d^\prime\geq 0$. Each of the sets of weights in this chain of
inclusions is a positive convex cone, with the exception of $\weightsqa$ which
is closed under multiplication by nonnegative scalars, but \emph{not} under
addition.

The affine group $\aff$ leaves $\weightshol$, $\weightsqa$, $\weights_d\,(0\leq
d\leq\infty)$, $\weights$, $\test$, $\schw$, $\galg$, and $\poly$ invariant. The
orbit of a quasi-analytic weight under the general linear group contains a
standard quasi-analytic weight. If $A\in\aff$, then $f\mapsto Af$ implements an
isomorphism between the topological vector spaces $\opw$ and $\op_{Aw}$.
\end{proposition}

\begin{proof}
The fact that $\weightsqa$ is not closed under addition will be illustrated in
the discussion following Proposition~\ref{prop:sequences}. It has already been
noted that $\weightshol\subset\weightsqa$ and the rest is either a routine
verification or contained in Lemma~\ref{lem:invarianceweightsqa}.
\end{proof}

\begin{remark}
Let us note explicitly that $\test$, $\schw$, $\galg$, and $\poly$ are subspaces
of $\opw$ for all $w\in\weightsqa$ and that Proposition~\ref{prop:invariance}
implies that in order to prove that these subspaces are dense in $\opw$ for all
$w\in\weightsqa$, it is sufficient to do so for standard quasi-analytic weights.
\end{remark}

The definition of the various classes of weights involves only the behavior at
infinity; this is a special case of the following lemma:

\begin{proposition}\label{prop:atinf}
Let $w$ and $\wtilde$ be weights and suppose there exist $C\geq 0$ and $R\geq 0$
such that $w(x)\leq C \wtilde(x)$ for all $x$ with $\Vert x
\Vert\geq R$. Then, if $\wtilde\in\weights_d$ for some $0\leq
d\leq\infty$ \textup{(}resp., $\wtilde\in\weightshol$\textup{)}, we also have $w\in\weights_d$
\textup{(}resp., $w\in\weightshol$ and $\tau_w\geq \tau_{\wtilde}$\textup{)}. If $\wtilde$ is
quasi-analytic with respect to the basis $\{v_1,\ldots,v_n\}$, then so is $w$.
\end{proposition}

\begin{proof}
This is immediate from the definitions, except when $\wtilde\in\weightsqa$. In that
case Lemma~\ref{lem:invarianceweightsqa} allows us to assume that $\wtilde$ is
standard quasi-analytic. Fix $1\leq j\leq n$. If there exists $r\geq 0$ such
that $w(x)=0$ whenever $\vert x_j\vert\geq r$, then $M_{w}(e_j,m)\leq r^m
\Vert w\Vert_\infty$ for all $m$, making it obvious that
$\sum_{m=1}^\infty M_{w}(e_j,m)^{-1/m}=\infty$. If $w$ does not have such
bounded support in the direction of the $j$th coordinate, let $C$ and $R$ be as
in the hypotheses of the proposition. Choose and fix $y\in\Rn$ with $\vert y_j\vert\geq R+1$ and
$w(y)\neq 0$. A comparison of geometric progressions shows that there exists $N$
with the property that, for all $m\geq N$ and for all $x\in\Rn$ such that $\vert
x_j\vert\leq R$, we have
\begin{equation*}
\vert x_j^m\cdot w(x)\vert\leq\vert y_j^m\cdot w(y)\vert.
\end{equation*}
Then the necessary divergence follows from the observation that, for $m\geq N$,
\begin{align*}
M_{w}(e_j,m)&=\sup \{\,\vert x_j^m \cdot w(x)\vert \mid x\in\Rn\,\}\\&=\sup
\{\, \vert x_j^m\cdot w(x)\vert\mid x\in\Rn,\,\vert x_j\vert\geq R\,\}\\ &\leq \sup
\{\,\vert x_j^m\cdot w(x)\vert\mid x\in\Rn,\,\Vert x\Vert\geq R\,\}\\&\leq C\sup
\{\,\vert x_j^m\cdot \wtilde (x)\vert\mid x\in\Rn,\,\Vert x\Vert\geq R\,\}\\&\leq C
M_{\wtilde}(e_j,m).
\end{align*}
\end{proof}

\subsection{Ostrowski's construction and weights on $\R$}\label{subsec:Ostrowski}
We now discuss a well-known type of construction, going back to Ostrowski \cite{Ostrowski}, by which we assign an even
weight on $\R$ to a nonnegative logarithmically convex sequence. This in turn
yields a method which allows us to pass from a given rapidly decreasing weight on $\R$
to an even rapidly decreasing majorant which has a certain minimal degree of
regularity, while relevant other properties remain unchanged in the process.
These two constructions \eqref{eq:sequencetoweight} and
\eqref{eq:weighttoweight} are important auxiliary tools in the study of weights.

Ostrowski's construction is as follows:

Suppose $\{a(m)\}_{m=0}^\infty\subset\R_{\geq 0}$ is logarithmically convex,
i.e., $a(m)^2\leq a(m-1)a(m+1)$ for $m\geq 1$. Define the associated weight
$\wtilde$ on $\R$ by
\begin{equation}\label{eq:sequencetoweight}
\wtilde(t)= \inf_{m\geq 0} \frac{a(m)}{\vert t\vert^m}\quad(t\in\R).
\end{equation}
In case $t=0$ we interpret this expression (and similar ones in the sequel) as
$a(0)$.

If $\{a(m)\}_{m=0}^\infty$ is not strictly positive (equivalently: $a(m)=0$ for
all $m\geq 1$), then $\wtilde(0)=a(0)$ and $\wtilde(t)=0$ for $\vert t\vert>0$,
implying that $M_\wtilde(e,m)=a(m)$ for all integral $m\geq 0$.

If the sequence is strictly positive, then as a consequence of the logarithmic
convexity we have the following local expression:

\begin{equation}\label{eq:sequencetoweightexplicit}
\wtilde(t)=
\begin{cases}
a(0) &\text{if $0\leq \vert t\vert \leq a(1)/a(0)$},\\ a(m)/\vert
t\vert^m&\text{if $a(m)/a(m-1)\leq\vert t\vert\leq a(m+1)/a(m)$}\quad(m\geq
1),\\ 0 &\text{if $\vert t\vert>\sup_{m\geq 1}a(m)/a(m-1)$}.
\end{cases}
\end{equation}
Some of the transition points can coincide, but there is no ambiguity. The value of $\sup_{m\geq 1}a(m)/a(m-1)$ can be finite or infinite; if it is finite, then the value
of $\wtilde$ at $\sup_{m\geq 1}a(m)/a(m-1)$ is left open.

For each integral $m\geq 0$ the graph of $t\mapsto\vert t\vert^m\cdot
\wtilde(t)$ is evidently rather simple and from it one concludes that again
$M_\wtilde(e,m)=a(m)$ for all integral $m\geq 0$. We also note that $\wtilde$ is
continuous, even, and strictly positive on the interval $(-\sup_{m\geq
1}a(m)/a(m-1),\sup_{m\geq 1}a(m)/a(m-1))$. Furthermore, $\wtilde$ is nonincreasing on $[0,\infty)$ and $s\mapsto
-\log \wtilde(e^s)$ is piecewise linear, nondecreasing, and convex on $(-\infty, \log\sup_{m\geq
1}a(m)/a(m-1))$.

We record that for all nonnegative logarithmically convex sequences
$\{a(m)\}_{m=0}^\infty$ the series $\sum_{m=1}^\infty a(m)^{-1/m}$ is divergent
if and only if the associated weight $\wtilde$ in
\eqref{eq:sequencetoweight} is quasi-analytic.

For the second construction, i.e., passing from a given rapidly decreasing
weight $w$ on $\R$ to a more regular even rapidly decreasing majorant, we apply
the first construction to the nonnegative sequence $\{M_w(e,m)\}_{m=0}^\infty$
(which by Remark~\ref{rem:convexity} is logarithmically convex) and obtain
\begin{equation}\label{eq:weighttoweight}
\wtilde(t)= \inf_{m\geq 0} \frac{M_w(e,m)}{\vert t\vert^m}\quad(t\in\R).
\end{equation}
It follows from the definitions that $\wtilde$ is an even rapidly decreasing
majorant of $w$. From the above discussion of \eqref{eq:sequencetoweight} we see
that $M_\wtilde(e,m)=M_w(e,m)$ for all integral $m\geq 0$. In particular,
$\wtilde$ is a quasi-analytic weight if and only if $w$ is. It is also obvious
from the discussion of \eqref{eq:sequencetoweightexplicit} that $\wtilde$ has a certain
minimal degree of regularity and a convexity property if $\supp w\nsubseteq \{0\}$.

Although we will make no specific use of it, let us mention that in case we
start with an arbitrary rapidly decreasing weight $w$ such that $\supp
w\nsubseteqq\{0\}$, then $\sup_{m\geq 1}M_w(e,m)/M_w(e,m-1)=\Delta_w$, where
$\Delta_w=\sup_{w(t)\neq 0}\vert t\vert$ is the possibly infinite radius of
the symmetric convex hull of the support of $w$. Therefore
$\Delta_\wtilde=\Delta_w$. This is evidently also true if $\supp
w\subseteq\{0\}$, so that the assignment $w\mapsto\wtilde$ satisfies $\Delta_\wtilde=\Delta_w$
for all rapidly decreasing
weights $w$. To conclude with, if $\Delta_w<\infty$ one can show that
$\wtilde(\Delta)=\lim_{s\uparrow\Delta}\sup_{\vert t\vert\geq s}w(t)=\inf_{0\leq
s<\Delta}\sup_{\vert t\vert\geq s}w(t)$, which then gives the missing value in
\eqref{eq:sequencetoweightexplicit}.

\subsection{Quasi-analytic weights}\label{subsec:qaw}

In this section we establish some basic characteristics of quasi-analytic
weights and then proceed to classify such weights.

First of all, there is an obvious and important way to construct a weight
$w^\prime$ on $\Rn$ from an even weight $w$ on $\R$ by putting
$w^\prime(x)=w(\Vert x\Vert)$ for $x\in\Rn$. Then $M_{w^\prime}(v,m)=\Vert
v\Vert^m M_w(e,m)$ for $v\in\Rn$ and all integral $m\geq 0$. This implies that
for $w^\prime$ quasi-analyticity and quasi-analyticity with respect to all bases
of $\Rn$ are both equivalent to the quasi-analyticity of $w$ on $\R$. As a
consequence, if $w$ is quasi-analytic on $\R$, then by Proposition
\ref{prop:atinf} all weights on $\Rn$ majorized at infinity by $w^\prime$ are
again quasi-analytic. To be more precise: by this proposition such minorants at
infinity are quasi-analytic with respect to any basis of $\Rn$. It will become
apparent from the classification of quasi-analytic weights on $\R$ and the
explicit examples in Section~\ref{subsec:examplesqa} that the collection of such
radial quasi-analytic weights on $\Rn$ and their minorants is in many practical
situations already an interesting one to work with.

In its essence, the situation is nevertheless fundamentally different. Whenever
$n\geq 2$, there exist quasi-analytic weights on $\Rn$ which are
\emph{not} obtainable in this way as a minorant at infinity of an essentially one-dimensional quasi-analytic weight. Indeed, we will show shortly that there
exist weights which are quasi-analytic with respect to precisely one basis, up
to scalar multiplication. These weights then clearly provide examples of such
nonminorants at infinity.

This possibility of the existence of an essentially unique basis related to a
particular element of $\weightsqa$ has the same origin as the failure of
$\weightsqa$ to be closed under addition, a fact which was already mentioned in
Proposition~\ref{prop:invariance}. To construct examples of these phenomena we
need the following proposition:

\begin{proposition}\label{prop:sequences}
For $k\geq 2$ there exist strictly positive logarithmically convex sequences
$\{a_j(m)\}_{m=0}^\infty$ for $j=1,\ldots,k$ such that
\begin{equation}
\sum_{m=1}^\infty a_j(m)^{-1/m}=\infty\quad(j=1,\ldots,k),\label{eq:divergence}
\end{equation}
but
\begin{equation}
\sum_{m=1}^\infty \left(\max (a_{j_1}(m),a_{j_2}(m))\right)^{-1/m}<\infty\quad(1\leq j_1\neq j_2\leq k).\label{eq:convergence}
\end{equation}
\end{proposition}

\begin{proof}
Consider the convex function $f:[1,\infty)\mapsto\R$ defined by $f(x)=2x\log x$.
The idea is to construct convex functions $f_j:[1,\infty)\mapsto\R$
$(j=1,\ldots,n)$ such that
\begin{align}\label{eq:sequences}
\begin{cases}
\max (f_{j_1},f_{j_2})&=f \quad (1\leq j_1\neq j_2\leq k),\\
\sum_{m=1}^\infty e^{-f_j(m)/m}&=\infty.
\end{cases}
\end{align}
Once this is accomplished, we can define $\log a_j(m)=f_j(m)$ for $j=1,\ldots,k$
and $m\geq 1$. Since it is obviously possible to define the $a_j(0)$ in
addition, preserving the strictly positive and logarithmically convex properties
of $\{a_j(m)\}_{m=1}^\infty$, this yields sequences satisfying
\eqref{eq:divergence} and \eqref{eq:convergence}.

The $f_j$ are constructed as partially ``tangentialized'' versions of $f$, as
follows. Let $N_1\geq 1$ be an arbitrary integer. Consider the tangent $T_1$ to
the graph of $f$ at $(N_1, f(N_1))$. If $T_1(x)$ is the ordinate of a point at
this tangent with abscissa $x$, then $T_1(x)/x$ trivially tends to a limit as
$x\rightarrow\infty$. It is therefore possible to select integral
$N_1^\prime>N_1$ such that $\sum_{m=N_1}^{N_1^\prime}\exp(-T_1(m)/m)>1$. Having
done this, consider the other tangent $T_1^\prime$ to the graph of $f$ which
passes through $(N_1^\prime,T_1(N_1^\prime))$ and which is tangent to the graph
of $f$ at, say, $(R_1,f(R_1))$ for some real $R_1>N_1^\prime$. Then $f$ admits
the possibility of ``tangentialization'' along the tangents $T_1$ and
$T_1^\prime$ on the interval $I_1=[N_1,N_1^\prime]\cup [N_1^\prime,R_1]$ into a
function $\widetilde f$, which is given by $\widetilde f(x)=f(x)$ for
$x\in[1,N_1]\cup[R_1,\infty)$, by $\widetilde f(x)=T_1(x)$ for
$x\in[N_1,N_1^\prime]$, and by $\widetilde f(x)=T_1^\prime(x)$ for
$x\in[N_1^\prime,R_1]$. The point of this procedure is that
$\sum_{m=N_1}^{N_1^\prime}\exp(-\widetilde f(m)/m)>1$ by construction, while
$\widetilde f$ is again convex as a consequence of our choice to work with
tangents. By this inherited convexity, we can repeat this procedure indefinitely
and determine a collection of possible tangentializations of $f$ on mutually
disjoint intervals $I_r$ ($r=1,2,\ldots$), each interval decomposable as
$I_r=[N_r,N_r^\prime]\cup [N_r^\prime,R_r]$ such that
$\sum_{m=N_r}^{N_r^\prime}\exp(-T_r(m)/m)>1$, and with the property that a
realization of this tangentialization of $f$ on any subcollection of the total
collection of all intervals $I_r$ again yields a convex function on
$[1,\infty)$. Therefore, if we let $f_j$ be defined by tangentializing $f$ on
the subcollection of all intervals $I_r$ with $r=j\mod k$, then the $f_j$ are
convex. They satisfy the first equation in
\eqref{eq:sequences} since the $I_r$ are disjoint and the second equation because $\sum_{m=N_r}^{N_r^\prime}\exp(-f_j(m)/m)=\sum_{m=N_r}^{N_r^\prime}\exp(-T_r(m)/m)>1$
whenever $r=j\mod k$.
\end{proof}

Returning to our phenomena, let us first show that $\weightsqa$ is not closed
under addition. Select sequences as in the proposition for $k=2$ and apply the
construction
\eqref{eq:sequencetoweight} to both sequences, obtaining even quasi-analytic
weights $w_{1,2}$ on $\R$ (where tildes have been omitted for short). The canonical radial procedure (which is the
identity if $\Rn=\R$) yields radial quasi-analytic weights $\wprime_{1,2}$
on $\Rn$. Now since obviously
\begin{align*}
M_{\wprime_1+\wprime_2}(v,m)&\geq
\max (M_{\wprime_1}(v,m),M_{\wprime_2}(v, m))\\
&=\Vert v\Vert^m \max (M_{w_1}(e,m),M_{w_2}(e, m))\\ &=\Vert
v\Vert^m\max (a_1(m),a_2(m))
\end{align*}
for all $v\in\Rn$ and all $m\geq 0$, \eqref{eq:convergence} shows that
$\wprime_1+\wprime_2$ is not quasi-analytic with respect to any basis.
We have therefore found two weights on $\Rn$, each quasi-analytic with respect
to all bases and such that their sum is not quasi-analytic. In particular,
$\weightsqa$ is in all dimensions not closed under addition.

In order to show the possibility of the existence of a distinguished basis for a
quasi-analytic weight on $\Rn$ ($n\geq 2$), take strictly positive
logarithmically convex sequences $\{a_j(m)\}_{m=0}^\infty$ for $j=1,\ldots,n$
satisfying
\eqref{eq:divergence} and \eqref{eq:convergence} and obtain even quasi-analytic
weights $w_j\,(j=1,\ldots,n)$ on $\R$ by
\eqref{eq:sequencetoweight}. Define $w:\R^n\mapsto\R_{\geq
0}$ as the tensor product $w(x)=\prod_{j=1}^n w_j(x_j)$. It is immediate
that $w$ is quasi-analytic with respect to the standard basis. Furthermore, if
$v=\sum_{j=1}^n \lambda_j e_j\in\Rn$, then consideration of the restriction of
$w$ to the union of each possible pair of coordinate axes shows that for $1\leq j_1\neq
j_2\leq n$ and integral $m\geq 0$:
\begin{equation*}
M_w(v,m)\geq \max\left(a_{j_1}(m),a_{j_2}(m)\right)
\left(\min(\vert\lambda_{j_1}\vert,\vert\lambda_{j_2}\vert)\right)^m
\min \left(a_{j_1}(0), a_{j_2}(0)\right)\prod_{j\neq j_1,j_2} a_j(0).
\end{equation*}
Therefore, if $w$ is quasi-analytic with respect to a basis containing $v$, then
by \eqref{eq:convergence} we must have
$\min(\vert\lambda_{j_1}\vert,\vert\lambda_{j_2}\vert)=0$ for all $j_1\neq j_2$.
This implies that the standard basis is the only possible choice, up to scaling.

This construction concludes our discussion of basic characteristics of
quasi-analytic weights and we will now start the classification. The first step
is a reduction to $\R$, for which we need a preparatory result.

\begin{proposition}[cf.\ \cite{Zapata, Nachbin1}]\label{prop:rootsareqa}
Let $w$ be a weight on $\Rn$. If $w$ is quasi-analytic with respect to the basis
$\{v_1,\ldots,v_n\}$, then $w^{\nu}$ is also quasi-analytic with respect to this
basis for all $\nu>0$.
\end{proposition}

\begin{proof}
By Lemma~\ref{lem:invarianceweightsqa} we may assume that $w$ is standard
quasi-analytic. We may obviously also assume that $\Vert w\Vert_\infty=1$. Note
that for all $1\leq j\leq n$, all strictly positive integers $p$, and all
integral $m\geq 1$ we have trivially
\begin{equation*}
M_{w^{1/p}}(e_j,m)^{1/m}=M_w(e_j,mp)^{1/mp}.
\end{equation*}
The third part of Lemma~\ref{lem:elementarytwo} and the second part of Lemma
\ref{lem:elementaryone} then imply that $w^{1/p}$ is standard quasi-analytic for
all strictly positive integral $p$. Since $w$ tends to zero at infinity, we have
$w^\nu\leq w^{1/p}$ at infinity whenever $p>1/\nu$. An appeal to Proposition
\ref{prop:atinf} finishes the proof.
\end{proof}

\begin{proposition}\label{prop:tensorprod}
The weight $w$ is a standard quasi-analytic weight on $\Rn$ if and only if there
exist even quasi-analytic weights $w_1,\ldots,w_n$ on $\R$ such that
$w(x)\leq\prod_{j=1}^n w_j(x_j)$ for all $x\in\Rn$.
\end{proposition}

\begin{proof}
The ``if''-part is obvious. As to the converse, first note that from Remark
\ref{rem:convexity} we know that $\{M_{w^{1/n}}(e_j,m)\}_{m=0}^\infty$ is a
nonnegative logarithmically convex sequence for each $j=1,\ldots,n$. We can
therefore apply construction
\eqref{eq:sequencetoweight} to each of these sequences and define
\begin{equation*}
\wtilde_j(t)=\inf_{m\geq 0} \frac{M_{w^{1/n}}(e_j,m)}{\vert t\vert^m}\quad(t\in\R).
\end{equation*}
Now since $w$ is standard quasi-analytic on $\Rn$ we conclude from Proposition
\ref{prop:rootsareqa} that $\sum_{m=1}^\infty M_{w^{1/n}}(e_j,m)^{-1/m}=\infty$
$(j=1,\ldots,n)$. As noted when discussing construction
\eqref{eq:sequencetoweight}, this implies that each $\wtilde_j$ is an (obviously
even) quasi-analytic weight on $\R$. From the definitions one sees that
$w^{1/n}(x)\leq \wtilde_j(x_j)$ for all $x\in\Rn$ and $j=1,\ldots,n$. The
$\wtilde_j$ are therefore as required.
\end{proof}

Since the orbit of a quasi-analytic weight under the affine group contains a
standard quasi-analytic weight, Proposition~\ref{prop:tensorprod} shows that the
classification of quasi-analytic weights on $\Rn$ is now reduced to
quasi-analytic weights on $\R$. The key in that case is the following lemma, in
which for the real line the defining property
\eqref{eq:divergentseries} is linked to a classical condition. This connection
seems to have gone largely unnoticed although there are similarities with
\cite[proof of Theorem 2]{Lin}, under additional regularity conditions on the
weight.

\begin{lemma}\label{lem:convexqa}
Suppose $w$ is a weight on $\R$ for which there exists $R>0$ such that
$w(t)=w(-t)>0$ for $\vert t\vert>R$ and, additionally, $s\mapsto-\log w (e^s)$ is
convex on $(\log R,\infty)$. Then $w$ is a quasi-analytic weight if and only if
\begin{equation}\label{eq:intdiverges}
\int_R^\infty\frac{\log w(t)}{1+t^2}=-\infty.
\end{equation}
\end{lemma}

\begin{proof}
The basic argument and the key inequalities that we will employ go (at least)
back to Horv\'ath \cite{Horvath1} (alternatively, see \cite[p.~98]{Koosis}) and
we refer the reader to these sources for details.

Leaving the lemma aside for the moment, let us consider the following situation.

Let $R\geq 1$ and suppose that $w:(R,\infty)\mapsto\R_{>0}$ is rapidly
decreasing such that $\psi(s)=-\log w (e^s)$ is convex and strictly increasing
on $(\log R,\infty)$.

For all real $\rho\geq 0$ (the lower bound stems from the hypothesis that $\psi$
is strictly increasing) put
\begin{equation}
M(\rho)=\sup_{t>R} t^\rho w(t).
\end{equation}
Then the convexity of $\psi$ implies that $w$ can be retrieved on $(R,\infty)$
as
\begin{equation}\label{eq:retrieved}
w(t)=\inf_{\rho\geq 0} \frac{M(\rho)}{t^\rho}\quad(t> R).
\end{equation}
Define $q:[0,\infty)\mapsto\R_{\geq 0}$ by
\begin{equation*}
q(t)= \inf_{m=0,1,2,\ldots} \frac{M(m)}{t^m}\quad(t\geq 0).
\end{equation*}
In view of \eqref{eq:retrieved}, the restriction of $q$ to $(R,\infty)$ can be
viewed as an approximation of $w$. Using that $R\geq 1$ one then derives the key
inequalities
\begin{equation}\label{eq:approximation}
w(t) \leq q(t)\leq t w(t)\quad(t>R).
\end{equation}
Since $q$ is nonincreasing on $[0,\infty)$ and $w$ is strictly positive on
$(R,\infty)$, the first of these inequalities implies that $\log q$ is bounded
from below on any compact subset of $[0,\infty)$. Evidently $\log q$ is bounded
from above by $\log M(0)$. Therefore, using
\eqref{eq:approximation}, we conclude that
\begin{equation}\label{eq:intdivergesone}
\int_R^\infty\frac{\log w(t)}{1+t^2}=-\infty
\end{equation}
if and only if
\begin{equation}\label{eq:intdivergestwo}
\int_0^\infty\frac{\log q(t)}{1+t^2}=-\infty,
\end{equation}
since the term $\log t$ is immaterial, as is for $\log q$ the integral from $0$
to $R$.

On the other hand, the first part of Lemma~\ref{lem:elementarytwo} implies that
the $M(m)$ form a logarithmically convex sequence, which is obviously strictly
positive. Therefore a standard equivalence \cite[Theorem~19.11]{RudinRCA}, relevant
in the context of quasi-analytic classes, can be used to conclude that
$\sum_{m=1}^\infty M(m)^{-1/m}=\infty$ if and only if
\eqref{eq:intdivergestwo} holds.

Both equivalences together therefore yield as conclusion that $\sum_{m=1}^\infty
M(m)^{-1/m}=\infty$ if and only if \eqref{eq:intdivergesone} holds.

Returning to the lemma proper now, suppose that $w$ is as in the lemma and that
\eqref{eq:intdiverges} holds. Introduce $\psi(s)=-\log w (e^s)$ on $(\log
R,\infty)$ and choose any $s_0\in(R,\infty)$. The convexity of $\psi$ implies
that $(\psi(s)-\psi(s_0))/(s-s_0)$ is nondecreasing as a function of $s\in
(\log R,\infty)$. If this quotient remains bounded by $M\in\R$ when
$s\rightarrow\infty$, then $w(t)\geq C t^{-M}$ for some $C>0$ and
all $t>\exp s_0$, implying that
\eqref{eq:intdiverges} is violated. Hence the quotient tends to infinity, which implies that
$w$ is rapidly decreasing on $(R,\infty)$. Since the quotient is, in particular,
eventually strictly positive, $\psi$ is eventually strictly increasing. Hence we
may assume, by enlarging $R$ if necessary, that this is actually the case on
$(\log R,\infty)$. Evidently, we can also assume $R\geq 1$.

We are then in the situation as described above and conclude from
\eqref{eq:intdivergesone} that $\sum_{m=1}^\infty M(m)^{-1/m}=\infty$, where $M(m)=\sup_{t>R} t^m w(t)$.
Now since $w$ has unbounded support, the fact that $w$ is even and the argument
in the proof of Proposition~\ref{prop:atinf} show that $\Vert t^m\cdot
w\Vert_\infty$ is equal to $M(m)$ for all sufficiently large $m$. We conclude
that $w$ is a quasi-analytic weight, as was to be shown.

Conversely, if $w$ is a quasi-analytic weight as in the lemma, let us show that
\eqref{eq:intdiverges} holds. We again introduce $\psi(s)=-\log w (e^s)$ on $(\log
R,\infty)$ and choose $s_0\in(\log R,\infty)$. This time the quotients
$(\psi(s)-\psi(s_0))/(s-s_0)$ cannot remain bounded as $s\rightarrow\infty$
since this would imply that $w\notin\weights_\infty$, contradicting Lemma
\ref{lem:invarianceweightsqa}. Hence again we may assume, by enlarging $R$ if
necessary, that $\psi$ is strictly increasing on $(\log R,\infty)$ for some
$R\geq 1$, again placing us in the situation as described above. Reversing the
arguments in the previous paragraph then shows that \eqref{eq:intdivergesone}
holds.
\end{proof}

We briefly review some well-known facts on majorants in view of the next
theorem.

As a preparation, suppose $w^\prime$ is an even weight on $\R$ such that
$s\mapsto -\log w^\prime(e^s)$ is convex on $\R$, where we use the convention
that $\log 0=-\infty$. Then, if $w^\prime(t_j)=0$ for some sequence
$\{t_j\}_{j=1}^\infty$ with $t_j\rightarrow\infty$, the convexity of $s\mapsto
-\log w^\prime(e^s)$ implies that actually $w^\prime(t)=0$ for all sufficiently
large $\vert t\vert$. We see that any even $w^\prime$ with this convexity
property is for all sufficiently large absolute values of the argument either
strictly positive or identically zero.

Next, if $w$ is a weight on $\R$, consider the set ${\overline{\weights}}_w$ of
all even weights $\wtilde$ on $\R$ that are majorants of $w$ such that $s\mapsto -\log
\wtilde(e^s)$ is convex on $\R$. Note that ${\overline{\weights}}_w$ contains at least
the sufficiently large constants, so that we may define $\overline
w(x)=\inf_{w^\prime\in {\overline{\weights}}_w} w^\prime (x)$. Then $\overline
w\in{\overline{\weights}}_w$ and by the previous paragraph $\overline w$ is for
all sufficiently large absolute values of the argument either strictly positive or
identically zero.

The quasi-analytic weights on the real line are now classified in the following
theorem:

\begin{theorem}\label{thm:classonedim}
Let $w$ be a weight \textup{(}i.e., an arbitrary nonnegative bounded function\textup{)} on $\R$.
Then the following are equivalent:
\begin{enumerate}
\item $w$ is a quasi-analytic weight, i.e.,
\begin{equation*}
\sum_{m=1}^\infty \Vert x^m \cdot w(x)\Vert_\infty^{-1/m}=\infty.
\end{equation*}
\item For all $R>0$,
\begin{equation*}
\int_R^\infty\frac{\log \overline w(t)}{1+t^2}=-\infty.
\end{equation*}
\item There exist a weight $\wtilde$ on $\R$ and $R>0$ such that
$w(t)\leq \wtilde(t)$ and $\wtilde(t)=\wtilde(-t)>0$ both hold for $\vert
t\vert>R$, $s\mapsto -\log\wtilde (e^s)$ is convex on $(\log R,\infty)$, and
\begin{equation*}
\int_R^\infty\frac{\log\wtilde(t)}{1+t^2}=-\infty.
\end{equation*}
\item There exist a strictly positive $\wtilde\in {\overline{\weights}}_w$ of class $C^\infty$ and $\epsilon>0$
such that $\wtilde$ is constant on $[0,\epsilon]$ and strictly decreasing on
$[\epsilon,\infty)$, and such that
\begin{equation*}
\int_0^\infty\frac{\log \wtilde(t)}{1+t^2}=-\infty.
\end{equation*}
\item There exists $C\geq 0$ and a nondecreasing nonnegative function $\rho:[0,\infty)\mapsto\R_{\geq 0}$ of class $C^\infty$,
which is equal to zero on $[0,\epsilon]$ for some $\epsilon>0$, such that
\begin{equation*}
w(t)\leq C \exp\left(-\int_0^{\vert t\vert}
\frac{\rho(s)}{s}\,ds\right)\quad(t\in\R)
\end{equation*}
and
\begin{equation*}
\int_{0}^\infty\frac{\rho(s)}{1+s^2}\,ds=\infty.
\end{equation*}
\end{enumerate}
\end{theorem}

\begin{proof}
We first show that (1) implies (4) in case $w$ has unbounded support, i.e., in
case there exist $t$ such that $w(t)\neq 0$ with $\vert t\vert$ arbitrarily
large. We apply a slight variation of the construction
\eqref{eq:weighttoweight} to the strictly positive logarithmically convex sequence
$\{M_w(e,m)\}_{m=0}^\infty$ and define $w^\prime$ by putting
\begin{equation*}
w^\prime(t)= 2\inf_{m\geq 0} \frac{M_w(e,m)}{\vert t\vert^m}\quad(t\in\R).
\end{equation*}
As noted when discussing \eqref{eq:sequencetoweight}, \eqref{eq:weighttoweight}
and \eqref{eq:sequencetoweightexplicit}, $w^\prime/2$ is an even quasi-analytic
weight majorizing $w$. The same then holds for $w^\prime$. Hence $w^\prime$ has
unbounded support; since $w^\prime$ is evidently nonincreasing on $[0,\infty)$,
it becomes apparent that $w^\prime$ is strictly positive. This implies that the
third possibility in \eqref{eq:sequencetoweightexplicit} does not occur and that
$w^\prime$ is therefore continuous. We also conclude from
\eqref{eq:sequencetoweightexplicit} that $s\mapsto -\log w^\prime (e^s)$ is
convex on $(-\infty,\infty)$. The graph of this map is in fact particularly
simple, since it is connected and consists of a horizontal half-line, followed
by a sequence of line segments with increasing strictly positive slope.

We therefore conclude from Lemma~\ref{lem:convexqa} that
\begin{equation*}
\int_0^\infty\frac{\log w^\prime(t)}{1+t^2}=-\infty.
\end{equation*}
Hence $w^\prime$ already has all the required properties, except that it is only
piecewise smooth. To improve this, note that the piecewise linear character of
the graph of $s\mapsto-\log w^\prime(e^s)$ shows that it is possible to smooth
this graph in a convex fashion while still remaining below the graph of
$s\mapsto-\log (w^\prime (e^s)/2)$. We choose any such convex smoothing, leaving
the original graph horizontal on $(-\infty,r)$ for some $r\in\R$, and declare
the smooth version to be the graph of $s\mapsto-\log\wtilde(e^s)$, thus defining
$\wtilde$ on $(0,\infty)$. Then $\wtilde$ is on $(0,\infty)$ a majorant of
$w^\prime/2$ by construction and evidently of class $C^\infty$. Since $\wtilde$
equals (unchanged) the constant $2M_w(e,0)$ on $(0,\exp r)$, it is evident that
$\wtilde$ can be extended to an even $C^\infty$ weight on $\R$ which then, since
$\wtilde(0)=2M_w(e,0)>M_w(e,0)=w^\prime(0)/2$ and since $w^\prime/2$ is even, is
a majorant of $w^\prime/2$ on $\R$. We already observed that $w^\prime/2\geq w$,
so $\wtilde$ is a majorant of $w$. By construction $\wtilde\leq w^\prime$, so
the relevant integral involving $\wtilde$ is divergent. Therefore $\wtilde$ has
all the required properties. This shows that (1) implies (4) if $w$ has
unbounded support.

In case $w$ has bounded support, we choose $C>0$ such that $C\exp (-\vert
t\vert)$ majorizes $w$ globally. By the previous result there is a suitable
majorant $\wtilde$ for $C\exp(-\vert t\vert)$, hence for $w$. This shows that
(1) implies (4) in all cases.

The equivalence of (4) and (5) is matter of restating, using
\eqref{eq:convexrelation} to define $\rho$ in terms of $\wtilde$ and vice versa. One reverses the order of integration in the resulting double
integral; details are left to the reader.

It is trivial that (4) implies (3).

To see that (3) implies (1), note that Lemma~\ref{lem:convexqa} asserts that
$\wtilde$ is a quasi-analytic weight. Then by Proposition~\ref{prop:atinf} the
same is true for $w$.

We have shown so far that (1), (3), (4), and (5) are equivalent. If $\wtilde$ is
as in (4) note that $\overline w\leq\wtilde$, so (2) follows.

Finally, if (2) holds, first suppose that $\overline w(t)=0$ for all
sufficiently large $\vert t\vert$. Then $w$ has bounded support and (1) holds.
In the remaining case $\overline w$ is strictly positive for all sufficiently large
absolute values of the argument, so that $\overline w$ satisfies the conditions
for $\wtilde$ in (3) if $R$ is sufficiently large. Thus (1) follows once more.
\end{proof}

\begin{remark}
The various conditions in the proposition each cover different aspects. The
first is the relevant one for the quasi-analytic property of certain transforms,
as will become apparent in Section \ref{sec:density}. The second condition
establishes the link between the notion of a quasi-analytic weight and
regularization. The requirement in (2) that $R$ is arbitrary is necessary to
avoid undue influence of an interval on which $\overline w$ vanishes
identically. The third condition aims at the practice since it provides a
positive criterion for a given weight which is easier to apply than the other
conditions in the proposition. We will encounter a negative criterion
in Proposition~\ref{prop:integralsonedim}. With respect to the fourth, let us note
that there is a continuous injection $\opwtilde\hookrightarrow\opw$ whenever
$\wtilde$ and $w$ are arbitrary weights such that $\wtilde$ majorizes $w$.
Closure results in $\opwtilde$ then transfer to $\opw$ by Lemma
\ref{lem:fundamentallemma}; the implication of (4) is therefore that, although
we made no such assumptions to start with, one can in many instances in fact
assume that a quasi-analytic weight is strictly positive on compact sets (which
is an important property) and smooth. Theorem~\ref{thm:improvements} below is
concerned with possibilities for such improvements for the various classes of
weights. The fifth condition enables one to write down whole families of
quasi-analytic weights. The collection of all smooth functions $\rho$ as in this
part is apparently sufficient to completely describe the quasi-analytic weights
on the real line. Note that the holomorphic weights are retrieved in the special
case where $\rho$ is eventually linear.
\end{remark}

The classification of quasi-analytic weights in arbitrary dimension is now a
matter of combining previous results.

\begin{theorem}\label{thm:classification}
Let $w$ be a weight \textup{(}i.e., an arbitrary nonnegative bounded function\textup{)} on $\Rn$.
Then the following are equivalent:
\begin{enumerate}
\item $w$ is a quasi-analytic weight, i.e., there exists a basis $\{v_1,\ldots,v_n\}$ of $\Rn$ such that
\begin{equation*}
\sum_{m=1}^\infty \Vert (v_j,x)^m \cdot w(x)\Vert_\infty^{-1/m}=\infty\quad(1\leq j\leq
n).
\end{equation*}
\item There exist an affine transformation $A$ and weights $w_j\,(j=1,\ldots,n)$ on
$\R$, each satisfying the five equivalent conditions in Theorem
\textup{\ref{thm:classonedim}} and such that
\begin{equation*}
w(Ax)\leq\prod_{j=1}^n w_j(x_j)
\end{equation*}
for all $x\in\Rn$ with $\Vert x\Vert$ sufficiently large.
\item There exists a quasi-analytic weight $\wtilde$ on $\Rn$ of class $C^\infty$ such
that $w(x)\leq\wtilde(x)$ and $\wtilde(x)=\wtilde(-x)$ for all $x\in\Rn$ and
such that $\inf_{x\in K}
\wtilde(x)>0$ for all nonempty compact subsets $K$ of $\Rn$.
\end{enumerate}
Furthermore, if the first statement holds for a basis $\{v_1,\ldots,v_n\}$, then
the second statement holds for any $A$ with linear component $A_0$ defined by
$A_0^t v_j=e_j\,(j=1,\ldots,n)$. Conversely, if the second statement holds and
$A_0$ is the linear component of $A$, then the first statement holds for the
basis $\{v_1,\ldots,v_n\}$ defined by $A_0^t v_j=e_j\,(j=1,\ldots,n)$. If the
first statement holds for a basis $\{v_1,\ldots,v_n\}$, then the weight
$\wtilde$ in the third statement can be chosen to be quasi-analytic with respect
to this same basis.
\end{theorem}

\begin{proof}
Assuming (1), Proposition~\ref{prop:invariance} yields a linear $A\in\aff$ such
that $A^{-1}w$ is standard quasi-analytic. By Proposition~\ref{prop:tensorprod}
$A^{-1}w$ is then majorized by the tensor product of quasi-analytic weights on
$\R$. Each of these satisfies all five equivalent conditions in Theorem
\ref{thm:classonedim}; hence (1) implies (2). Assuming (2) we note that the
tensor product is a quasi-analytic weight on $\Rn$ by Theorem
\ref{thm:classonedim} and Proposition~\ref{prop:tensorprod}. By Proposition
\ref{prop:atinf} we conclude that $A^{-1}w$ is quasi-analytic and the same is
then true for $w$ by Proposition~\ref{prop:invariance}. Hence (2) implies (1).

Assuming (1) we conclude as above that there exists a linear $A\in\aff$ such
that $A^{-1}w$ is majorized by the tensor product of quasi-analytic weights on
$\R$. Majorizing each of these in turn as in the fifth part of Theorem
\ref{thm:classonedim} we conclude that there exists a standard quasi-analytic
weight $w^\prime$ on $\Rn$ of class $C^\infty$ such that $A^{-1}w(x)\leq
w^\prime(x)$ and $w^\prime(x)=w^\prime(-x)$ for all $x\in\Rn$ and such that
$\inf_{x\in K} w^\prime(x)>0$ for all nonempty compact subsets $K$ of $\Rn$.
Therefore $Aw^\prime$ satisfies the requirements for $\wtilde$ and (1) implies
(3). It is trivial that (3) implies (1).

The additional statements on the bases follow from Lemma
\ref{lem:invarianceweightsqa}.
\end{proof}

The following proposition elaborates on quasi-analytic weights on the real line.
Note that the first part of the proposition yields a negative criterion for
quasi-analyticity of weights on $\R$. The example in the third part shows that
the necessary condition in the first part is not sufficient and that the
symmetry assumption in the second statement is necessary.

\begin{proposition}\label{prop:integralsonedim}
Let $w$ be a weight on $\R$.
\begin{enumerate}
\item If $w$ is quasi-analytic and Lebesgue measurable, then for all $R>0$ we have
\begin{equation*}
\int_R^\infty\frac{\log w(t)}{1+t^2}=\int_{-\infty}^{-R}\frac{\log w(t)}{1+t^2}=-\infty.
\end{equation*}
\item If $w$ is even and such that $s\mapsto -\log w(e^s)$ is convex on $\R$, then $w$
is quasi-analytic if and only if
\begin{equation*}
\int_R^\infty\frac{\log w(t)}{1+t^2}=-\infty
\end{equation*}
for all $R>0$.
\item There exists a strictly positive weight $w$ on $\R$ such that
$s\mapsto -\log w(e^s)$ and $s\mapsto -\log w(-e^s)$ are both convex on $\R$ and
such that
\begin{equation*}
\int_R^\infty\frac{\log w(t)}{1+t^2}=\int_{-\infty}^{-R}\frac{\log
w(t)}{1+t^2}=-\infty
\end{equation*}
for all $R>0$, but which is not quasi-analytic.
\end{enumerate}
\end{proposition}

\begin{proof}
In the first statement, the divergence of the integrals is implied by the second
part of Theorem~\ref{thm:classonedim} since $w(t)\leq \overline w(\vert t\vert)$
for all $t$. The second statement is immediate from again the second part of
Theorem~\ref{thm:classonedim} since now $w=\overline w$.

As to the third statement, Proposition~\ref{prop:sequences} provides two
strictly positive logarithmically convex sequences $\{a_j(m)\}_{m=0}^\infty$ for
$j=1,2$ such that
\begin{equation*}
\sum_{m=1}^\infty a_j(m)^{-1/m}=\infty\quad(j=1,2)
\end{equation*}
but
\begin{equation*}
\sum_{m=1}^\infty \left(\max (a_{1}(m),a_{2}(m))\right)^{-1/m}<\infty.
\end{equation*}
We apply \eqref{eq:sequencetoweight} to both sequences, obtaining even
quasi-analytic weights $w_{1,2}$ on $\R$ with the property that
$M_{w_{1,2}}(e,m)=a_{1,2}(m)$ for all integral $m\geq 0$ and such that the maps
$s\mapsto -\log w_{1,2}(e^s)$ are both convex on $\R$. The first part of the
present proposition shows that
\begin{equation*}
\int_R^\infty\frac{\log w_1(t)}{1+t^2}=\int_{-\infty}^{-R}\frac{\log w_2(t)}{1+t^2}=-\infty.
\end{equation*}
for all $R>0$. Choose $w=\charfunc_{\R_{<0}}w_1+\charfunc_{\R_{>0}}w_2 +
\charfunc_{\{0\}}$. Then, evidently, $M_w(e,m)=
\max (M_{w_1}(e,m),M_{w}(e, m))=\max (a_1(m),a_2(m))$ for all $m\geq
1$, so that $w$ is an example as required.
\end{proof}

We finally note that we have the following necessary condition for a weight to
be quasi-analytic in arbitrary dimension, again providing us with a negative criterion.

\begin{corollary}\label{cor:online}
Let $w$ be a quasi-analytic weight on $\Rn$. If $x,y\in\Rn$ are such that $y\neq
0$ and such that $t\mapsto w(x+ty)$ is Lebesgue measurable on $\R$, then for all
$R>0$ we have
\begin{equation*}
\int_R^\infty\frac{\log w(x+ty)}{1+t^2}=\int_{-\infty}^{-R}\frac{\log w(x+ty)}{1+t^2}=-\infty.
\end{equation*}
\end{corollary}

\begin{proof}
Since $\weightsqa$ is invariant under translations by Lemma
\ref{lem:invarianceweightsqa}, we may assume $x=0$. Let $\wtilde(t)=w(ty)$, then
$\wtilde$ is a weight on $\R$. Suppose $w$ is quasi-analytic with respect to the
basis $\{v_1,\ldots,v_n\}$. Select $v_j$ such that $(v_j,y)\neq 0$. Since
$M_\wtilde(e,m)=\vert (v_j,y)\vert^{-m}\sup\{\,\vert(v_j,ty)\vert^m\cdot
w(ty)\mid t\in\R\,\}\leq
\vert (v_j,y)\vert^{-m} M_w(v_j,m)$ for all $m\geq 0$, $\wtilde$ is a
quasi-analytic weight on $\R$. We now apply the first part of Proposition
\ref{prop:integralsonedim}.
\end{proof}

\subsection{Examples of quasi-analytic weights}\label{subsec:examplesqa}
On the basis of the previous results we can now determine explicit examples and
counterexamples of quasi-analytic weights on $\Rn$. We concentrate mainly on
minorants and majorants of radial weights.

\begin{proposition}
Suppose $R_0>0$, $C\geq 0$, and a nondecreasing function
$\rho:(R_0,\infty)\mapsto\R_{\geq 0}$ of class $C^1$ are such that
\begin{equation*}
\int_{R_0}^\infty\frac{\rho(s)}{s^2}\,ds=\infty.
\end{equation*}
If $w$ is a weight such that
\begin{equation*}
w(x)\leq C \exp\left(-\int_{R_0}^{\Vert x\Vert}\frac{\rho(s)}{s}\right)\,ds
\end{equation*}
whenever $\Vert x\Vert\geq R_0$, then the weight $w$ is quasi-analytic with
respect to all bases of $\Rn$.
\end{proposition}

Indeed, with the aid of the third part of Theorem~\ref{thm:classonedim} and the
discussion of \eqref{eq:convexrelation} one concludes that such $w$ is majorized
by the radial extension of an even quasi-analytic weight on $\R$; it is
therefore a quasi-analytic weight on $\Rn$.

The first part of the following result exhibits families of examples of
quasi-analytic weights in terms of elementary functions, incorporating the
holomorphic weights in a natural way. The second part shows that the first part
is in a sense sharp.

\begin{proposition}\label{prop:examplesqa}
Define repeated logarithms by $\log_0 t=t$ and, inductively, for $j\geq 1$, by
$\log_j t=\log (\log_{j-1}t)$, where $t$ is assumed to be sufficiently large for
the definition to be meaningful in the real context. For $j=0,1,2\ldots$ let
$a_j>0$ and let $p_j\in\R$ be such that $p_j=0$ for all sufficiently large $j$.
Put $j_0=\min\{j=0,1,2,\ldots\mid p_j\neq 1\}$. Let $C>0$ and suppose
$w:\Rn\mapsto\R_{\geq 0}$ is bounded. Then:
\begin{enumerate}
\item If $p_{j_0}<1$ and
\begin{equation*}
w(x)\leq C\exp\left(-\Vert x\Vert^2
\left(\prod_{j=0}^\infty\log_{j}^{p_j}a_j\Vert x\Vert\right)^{-1}\right)
\end{equation*}
for all sufficiently large $\Vert x\Vert$, then the weight $w$ is quasi-analytic
with respect to all bases of $\Rn$.
\item If $p_{j_0}>1$ and
\begin{equation*}
 w(x)\geq C\exp\left(-\Vert x\Vert^2 \left(\prod_{j=0}^\infty\log_{j}^{p_j}a_j\Vert x\Vert\right)^{-1}\right)
\end{equation*}
for all sufficiently large $\Vert x\Vert$, then $w$ is not a quasi-analytic
weight on $\Rn$.
\end{enumerate}
\end{proposition}

\begin{proof}
It is sufficient to prove this in one variable, in view of the properties of the
radial extension procedure. In that case, the second statement follows
immediately from Propositions~\ref{prop:atinf} and \ref{prop:integralsonedim}. The first statement follows from the third part of
Theorem~\ref{thm:classonedim} after some computation and asymptotics to check
that the relevant function is eventually convex.
\end{proof}

Concentrating on the radial case, a weight on $\Rn$ which for some
$C,a_0,a_1,a_2,\ldots>0$ and all sufficiently large $\Vert x\Vert$ is majorized
by whichever of the expressions
\begin{align*}
&C\exp\left(-\frac{\Vert x\Vert^{1-\nu}}{a_0}\right)\quad ,\\
&C\exp\left(-\frac{\Vert x\Vert}{a_0 \left(\log a_1\Vert
x\Vert\right)^{1+\nu}}\right)\quad ,\\&C\exp\left(-\frac{\Vert x\Vert}{a_0 \log
a_1\Vert x\Vert
\left(\log\log a_2\Vert x\Vert\right)^{1+\nu}}\right)\quad ,\\
&\ldots
\end{align*}
is in $\weightsqa$ if $\nu\leq 0$. Note that the first family of majorants is in
$\weightshol$, but that the others are not. Complementary, a weight which is for
some $\nu>0$ minorized at infinity by whichever of these expressions is not a
quasi-analytic weight.

For the sake of completeness we mention that explicit nonradial examples of
quasi-analytic weights on $\Rn$ in terms of elementary functions can be obtained
as tensor products of quasi-analytic weights on $\R$ taken from the above
proposition. All minorants at infinity of such weights are then again
quasi-analytic weights on $\Rn$.

\subsection{Improving the spaces $\opw$}\label{subsec:improving}
In this section we observe that many spaces $\opw$ are Fr\'echet and show that in a sense these
well-behaved spaces are ubiquitous. This is an important fact which we use in
Section~\ref{sec:density} to invoke vector-valued integration and in Section
\ref{sec:general} to show the redundancy of a hypothesis of continuity. We also
note that we can combine this completeness with an improvement of some
regularity properties of the weight, without leaving the class of weights under
consideration.

\begin{proposition}\label{prop:spaceisFrechet}
Let $w$ be a weight such that for all nonempty compact subsets $K$ of $\Rn$
there exists $C_K>0$ such that $K\subset\overline{\{\,x\in\Rn\mid
w(x)>C_K\,\}}$. Then $\opw$ is a Fr\'echet space.
\end{proposition}

\begin{proof}
The hypothesis implies that the complement of the zero locus of $w$ is dense,
hence $\opw$ is Hausdorff by Remark~\ref{rem:Hausdorff} and only completeness
remains to be checked. Let $\{f_n\}_{n=1}^\infty$ be a Cauchy sequence in
$\opw$. Let $K$ be an arbitrary nonempty compact subset of $\Rn$ and select
$C_K>0$ such that $K\subset\overline{A_K}$ where $A_K=\{\,x\in\Rn\mid
w(x)>C_K\,\}$. For all fixed $\alpha\in\Nn$, the restrictions of the $\parta
f_n$ to $A_K$ are then bounded on $A_K$ and are in fact a Cauchy sequence under
the supremum norm on the space of bounded functions on $A_K$. The analogous
statement is then by continuity true for the restrictions to
$\overline{A_K}\supset K$. We conclude that the $f_n$ are a Cauchy sequence in
$C^\infty(\Rn)$, the latter space being supplied with the usual Fr\'echet
topology of uniform convergence of all derivatives on all compact sets. A
routine verification shows that the limit in $C^\infty(\Rn)$ is in fact an
element of $\opw$ and that the sequence converges to this limit in the topology
of $\opw$.
\end{proof}

\begin{theorem}\label{thm:improvements}
\begin{enumerate}
\item Let $w$ and $\wtilde$ be weights such that:
\begin{enumerate}
\item $w(x)\leq\wtilde(x)$ for all $x\in\Rn$; and
\item $\wtilde$ is separated from zero on compact sets.
\end{enumerate}
Then $\opwtilde$ is a Fr\'echet space which is identified with a subspace of
$\opw$. The corresponding injection $\opwtilde\hookrightarrow\opw$ is
continuous, as is the injection $\opwtilde\hookrightarrow C^\infty(\Rn)$.
\item If $w\in\weights_d$ for some $0\leq d<\infty$ \textup{(}resp.,
$w\in\weightshol$\textup{)} there exists a radial $\wtilde\in\weights_d$ \textup{(}resp., a radial
$\wtilde\in\weightshol$ with $\tau_\wtilde=\tau_w$\textup{)} which satisfies the
requirements under \textup{(1a)} and \textup{(1b)}. If $w\in\weights_\infty$ \textup{(}resp.,
$w\in\weights\textup{)}$, then there exists $\wtilde\in\weights_\infty$ \textup{(}resp.,
$w\in\weights$\textup{)} which is of class $C^\infty$, radial, nonincreasing on each ray
emanating from the origin, and satisfying the requirements under \textup{(1a)} and \textup{(1b)}.
If $w\in\weightsqa$ is quasi-analytic with respect to the basis
$\{v_1,\ldots,v_n\}$ of $\Rn$, then there exists $\wtilde\in\weightsqa$ of class
$C^\infty$, satisfying the requirements under \textup{(1a)} and \textup{(1b)}, again
quasi-analytic with respect to $\{v_1,\ldots,v_n\}$ and such that
$\wtilde(x)=\wtilde(-x)$ for all $x\in\Rn$.
\end{enumerate}
\end{theorem}

\begin{remark}
It is important to note that for all weights $w$ of the various types that we have
distinguished, there exists a majorant $\wtilde$ \emph{of the same type} and
with the properties as in the first part.
\end{remark}

\begin{proof}
The nontrivial statements in the first part are clear in view of Proposition
\ref{prop:spaceisFrechet} and its proof. As to the second part, if
$w\in\weights_d\,(0\leq d<\infty)$ or $w\in\weightshol$ one can take
$\wtilde(x)=\sup_{\Vert y\Vert=\Vert x\Vert}w(y)+\exp(-\Vert x\Vert^2)$. The
case $w\in\weights$ is trivial. Consider the case $w\in\weights_\infty$. If
$\supp w$ is bounded then $C\exp(-\sqrt{\Vert x\Vert^2+1})$ for some
sufficiently large $C$ obviously suffices. Put
\begin{equation*}
w^\prime(x)=2\inf_{m=0,1,2,\ldots}\frac{\Vert \Vert y\Vert^m\cdot
w(y)\Vert_\infty}{\Vert x\Vert^m}\quad(x\in\Rn).
\end{equation*}
Since $\{\Vert \Vert y\Vert^m\cdot w(y)\Vert_\infty\}_{m=0}^\infty$ is a
strictly positive logarithmically convex sequence, the analogues of
\eqref{eq:sequencetoweightexplicit} and the discussion of \eqref{eq:weighttoweight}
apply, showing that $w^\prime$ can be modified into a majorant with all the
required properties. Finally, the third part of Theorem~\ref{thm:classification}
handles the quasi-analytic case.
\end{proof}

\section{Density in $\opw$}\label{sec:density}

The principal aim of this section is the proof of the density of various common
subspaces of the spaces $\opw$, as dependent on the properties of the weight
$w$. The principal tool is the Fourier transform on $\opwdual$. This transform
can be defined since, as it turns out, $\opwdual$ is canonically embedded in
$\tempdist$ for any weight $w$. If $w$ is in one of the various classes of
weights as in Section~\ref{sec:weights}, then there is a corresponding minimal
degree of regularity of the elements of $\Fourier(\opwdual)$, dependent of the
particular class. In the chain of inclusions in Proposition
\ref{prop:invariance} this regularity increases from right to left. To be more
precise: if $w\in\weights_d\,(0\leq d\leq\infty)$ then the Fourier transforms
are in fact functions and we show that these functions are of class $C^{[d]}$.
If $w\in\weightsqa$ they have the quasi-analytic property. If $w\in\weightshol$
the transforms are holomorphic on a tubular neighborhood of $\Rn$.

The notion of admissible spaces for a weight $w$, as motivated in Section
\ref{sec:introduction}, is introduced in Definition~\ref{def:admissiblespaces}.
Once the regularity properties for the elements of $\Fourier(\opwdual)$ have
been established (Theorem~\ref{thm:Fouriertransforms}), the density of these
admissible spaces in $\opw$ is immediate (Theorem
\ref{thm:admissiblespaces}).

The section concludes with a necessary condition for $\poly$ to be dense in
$\opw$ for an arbitrary rapidly decreasing weight $w$.

\subsection{Fourier transform on $\opwdual$}
The main result in this section is Theorem~\ref{thm:admissiblespaces},
describing the regularity of Fourier transforms of elements of $\opwdual$ when
$w\in\weights_0$. We need some preparation.

A basic fact for the spaces $\opw$ is the density of $\test$. This result is due
to Zapata \cite{Zapata} in a slightly more general context, to which we return in
Section~\ref{sec:applications}. We include the proof, since the result is vital
for the main results in this section.

\begin{proposition}\label{prop:testdense}
Let $w$ be a weight. Then $\test$ is dense in $\opw$.
\end{proposition}

\begin{proof}[Proof (Zapata)]
Choose $\theta\in C_c^\infty(\Rn)$ such that $0\leq \theta\leq 1$ and
$\theta(x)=1$ for $\Vert x\Vert\leq 1$. For $h>0$, define $\theta_h\in
C_c^\infty(\Rn)$ by $\theta_h(x)=\theta(hx)$; note that for all $\mu\in\Nn$
there is a constant $C_\mu\geq 0$ such that, for all $h>0$,
\begin{equation*} \Vert \partial^\mu
\theta_h\Vert_\infty \leq C_\mu h^{\vert\mu\vert}.
\end{equation*}

We claim that $\lim_{h\downarrow 0}
\theta_h f=f$ for all $f\in\opw$, which clearly proves the
proposition. In order to validate this claim, fix $f\in\opw$ and $\alpha\in\Nn$.
Let $\epsilon>0$ be given. By definition there exists $R>0$ such that
\begin{equation*}
\sup_{\Vert x\Vert\geq R}\vert \parta f(x)\cdot w(x)\vert<\epsilon/2.
\end{equation*}
Recalling the definition of the seminorm $\pa$, we see that there exists a
finite set of positive integers $c_\mu$ such that
\begin{equation*}
\pa(f-\theta_h f)\leq \left\Vert (1-\theta_h)\parta f\cdot w\right\Vert_\infty
+ \sum_{\mu\in\Nn, \mu\neq 0}c_\mu \left\Vert
\partial^\mu\theta_h\cdot\partial^{\alpha-\mu} f\,\cdot w\right\Vert_\infty.\label{expressionseminorm}
\end{equation*}
For $0<h<R^{-1}$ we have $\left\Vert (1-\theta_h)\parta f\cdot
w\right\Vert_\infty<\epsilon/2$, implying that for such $h$
\begin{equation*}
\pa(f-\theta_h f)< \epsilon/2
+\sum_{\mu\in\Nn, \mu\neq 0}c_\mu C_\mu h^{\vert\mu\vert} p_{\alpha-\mu}(f).
\end{equation*}
The claim follows by taking $h\neq 0$ small enough.
\end{proof}

If $w$ is a weight, then the injection $\schw\hookrightarrow\opw$ is continuous.
Since $\schw$ is evidently dense in $\opw$ by the previous proposition,
$\opwdual$ is canonically identified with a subspace of $\tempdist$. We use the
same notation for $\typ\in\opwdual$ and its image in $\tempdist$.

\begin{definition}
If $\typ\in\opwdual$, there exist a constant $C$, a nonnegative integer $L$, and
a finite set $\{\beta_1,\ldots,\beta_L\}$ of multi-indices such that, for all
$f\in\opw$,
\begin{equation}
\vert \langle \typ,f\rangle\vert\leq C \max_{k=1,\ldots,L} p_{\beta_k}(f).\label{eq:continuity}
\end{equation}
The minimum of $\max_{k=1,\ldots,L}\vert\beta_k\vert$, computed over all sets of
multi-indices for which there exists a $C$ such that
\eqref{eq:continuity} holds for all $f\in\opw$, is the \emph{order of $\typ$ in $\opwdual$},
denoted $\ord_\opw \typ$.
\end{definition}

\begin{remark}
If $\typ\in\opwdual$, then $\typ$ has also a finite order $\ord_\test
\typ$ as a tempered distribution. It is immediate from the definitions
that $\ord_\test\typ\leq\ord_\opw\typ$. Inequality can occur; we give an example
of this in one dimension for the constant weight $1$. Let $f(x)=x^{-2}\sin
x^3\,(x\in\R)$. Then $f\in C^{\infty}(\R)$, $f\in L_1(\R,dx)$, and
$f^\prime\notin L_1(\R,dx)$. Consider the tempered distribution $T_{f^\prime}$
of order zero determined by $f^\prime$. Since $\vert\langle
\typ_{f^\prime},\psi\rangle\vert\leq \Vert f\Vert_1\Vert\psi^\prime\cdot 1\Vert_\infty$
for $\psi\in\schw$, $\typ_{f^\prime}$ is continuous on $\schw$ in the topology
of $\op_1$. Let $\widetilde\typ_{f^\prime}$ denote the unique continuous
extension of $\typ_{f^\prime}$ to $\op_1$. Then $\ord_\test
\widetilde\typ_{f^\prime}=0$ and from the density of $\schw$ in $\opw$ one
concludes that $\ord_{\op_1}\widetilde\typ_{f^\prime}\leq 1$. Using that
$f^\prime\notin L_1(\R,dx)$ one sees that
$\ord_{\op_1}\widetilde\typ_{f^\prime}\neq 0$ and that therefore
$\ord_{\op_1}\widetilde\typ_{f^\prime}=1$.
\end{remark}

The results for the Fourier transform on $\opwdual$ follow mainly from the
degree of regularity of the map $\lambda\mapsto e_{i\lambda}$, which assumes
values in $\opw$ and is defined on a suitable domain. We need the following
definitions:

\begin{definition}\label{def:emaps}
Let $w$ be a weight.
\begin{enumerate}
\item If $w\notin\weightshol$, put $\domain=\Rn$. If $w\in\weightshol$, put
$\domain=\{\,\lambda\in\Cn\mid\Vert\ima\lambda\Vert<\tau_w\,\}$ where $\tau_w$
is as in Section~\ref{sec:weights}.
\item If $d$ is a nonnegative integer, $P\in\polyd$, and $w\in\weights_{d}$, define $\emapP:\domain\mapsto\opw$ by
$\emapP(\lambda)=Pe_{i \lambda}$ $(\lambda\in\domain)$.
\end{enumerate}
\end{definition}
The notation $\domain$ for the domains has been chosen in order to treat the
holomorphic and the nonholomorphic case in a uniform formulation as much as
possible.

Our next three results concern regularity of the maps $\emapP$.

\begin{proposition}\label{prop:mapiscontinuous}
Let $0\leq d<\infty$ be an integer. If $P\in\polyd$ and $w\in\weights_d$, then
$\emapP:\domain\mapsto\opw$ is continuous.
\end{proposition}

\begin{proof}
We treat the case where $w\notin\weightshol$; the case where $w\in\weightshol$
is similar. We then have to show that
\begin{equation*}
\lim_{\Vert
h\Vert\rightarrow 0}\pa\left(\emapP(\lambda+h)-\emapP(\lambda)\right)= 0
\end{equation*} for all
fixed $\lambda\in\Rn$ and $\alpha\in\Nn$. The definitions make it obvious that
\begin{equation*}
\pa\left(\emapP(\lambda+h)-\emapP(\lambda)\right)\leq\sum_{\mu+\nu=\alpha}c_{\mu,\nu}
\left\Vert \partial_x^\mu P \cdot \partial_x^\nu
(e_{i(\lambda+h)}-e_{i\lambda})\cdot w\right\Vert_\infty
\end{equation*}
for some nonnegative integers $c_{\mu,\nu}$. By a moment's thought, it is
therefore sufficient to show that
\begin{equation*}
\lim_{\Vert h\Vert\rightarrow 0}\left\Vert (1+\Vert x\Vert)^d\left\{(\lambda+h)^\nu
e_{ih}-\lambda^\nu\right\}\cdot w \right\Vert_\infty=0
\end{equation*}
for all $\nu$ occurring in the summation. This holds in fact for all
$\nu\in\Nn$. To see this, in view of
\begin{align*}
&\left\Vert (1+\Vert x\Vert)^d\left\{(\lambda+h)^\nu
e_{ih}-\lambda^\nu\right\}\cdot w\right\Vert_\infty\leq\\ &\leq \left\Vert
(1+\Vert x\Vert)^d\left\{(\lambda+h)^\nu e_{ih} -
\lambda^\nu e_{ih}\right\}\cdot w\right\Vert_\infty +
\left\Vert (1+\Vert x\Vert)^d\left\{\lambda^\nu e_{ih} - \lambda^\nu\right\}
\cdot w\right\Vert_\infty\\
&=\left\vert(\lambda+h)^\nu-\lambda^\nu\right\vert \left\Vert (1+\Vert
x\Vert)^d\cdot w\right\Vert_\infty +
\vert\lambda\vert^\nu\,\left\Vert (1+\Vert x\Vert)^d\left\{e_{ih} - 1\right\}\cdot w\right\Vert_\infty,
\end{align*}
it is sufficient to prove that
\begin{equation}\label{eq:limzero}
\lim_{\Vert h\Vert\rightarrow 0}\left\Vert (1+\Vert x\Vert)^d \left\{e_{ih} - 1\right\}\cdot w(x)\right\Vert_\infty=0.
\end{equation}
Let $\epsilon>0$. Choose $R>0$ such that $\sup_{\Vert x\Vert\geq
R}\left\vert(1+\Vert x\Vert)^d\left\{e_{ih}(x)-1\right\}\cdot
w(x)\right\vert<\epsilon$ for all $h$, which is possible since $w\in\weights_d$.
In addition, since $\vert
\exp(it)-1\vert\leq
\vert t\vert$ for $t\in\R$, we can arrange that also $\sup_{\Vert
x\Vert<R}\left\vert(1+\Vert x\Vert)^d
\left\{e_{ih}(x)-1\right\}\cdot w(x)\right\vert<\epsilon$ by taking $\Vert
h\Vert$ small enough. This proves \eqref{eq:limzero} and completes the proof
when $w\notin\weightshol$.
\end{proof}

\begin{proposition}\label{prop:emapsmooth}
Let $0\leq d<\infty$ be an integer. If $P\in\polyd$ and $w\in\weights_{d+1}$,
then $\frac{\partial\emap_P}{\partial
\lambda_j}$ exists for $j=1,\ldots,n$ on $\domain$ in the
topology of $\opw$ and is given by
\begin{equation*}
\left(\frac{\partial\emap_P}{\partial
\lambda_j}(\lambda)\right)(x)=ix_jP(x)e^{i(\lambda,x)}\quad\text{\textup{(}$\lambda\in\domain$,
$x\in\Rn$\textup{)}}.
\end{equation*}
If $w\in\weightshol$, the derivatives are to be interpreted as complex
derivatives.
\end{proposition}

\begin{proof}
We treat the case where $w\notin\weightshol$; the case where $w\in\weightshol$
is similar. Fix $1\leq j\leq n$. The definitions imply we then have to prove that
\begin{equation}
\lim_{h\rightarrow 0}\pa\left(U\,V_h\right)=0\label{eq:aim}
\end{equation}
for all fixed $\alpha\in\Nn$ and $\lambda\in\Rn$, where $U$ and $V_h$ $(0\neq
h\in\R)$ in $C^\infty(\Rn)$ are defined by
\begin{align*}
U(x)&=P(x)e^{i(\lambda,x)}&\quad&(x\in\Rn),\\
V_h(x)&=ix_j-\frac{e^{ihx_j}-1}{h}&\quad&(x\in\Rn).
\end{align*}
The definition of $\pa$ implies that is sufficient to prove
\begin{equation}\label{eq:subaim}
\lim_{h\rightarrow 0}\left\Vert \left(1+\Vert x\Vert\right)^d\partial_x^\nu V_h\cdot w\right\Vert_\infty=0
\end{equation}
for all $\nu\in\Nn$. For this, we distinguish three cases:
\begin{enumerate}
\item Assume $\vert\nu\vert=0$.  Let $\epsilon>0$ be given. From $\vert \exp(it)-1)\vert\leq \vert t\vert \,(t\in\R)$ we see that $\vert V_h(x)\vert\leq 2\Vert x\Vert$ for all $h\neq 0$ and all $x$. Since $w\in\weights_{d+1}$, we can therefore select $R>0$ such that $\sup_{\Vert x\Vert\geq R}\vert (1+\Vert x\Vert)^d V_h(x)\cdot w(x)\vert <\epsilon$ for all $h\neq 0$.
Next, since $\vert\exp(it)-1-it\vert\leq t^2/2$ for $t\in\R$ it follows easily
that we can then arrange that $\sup_{\Vert x\Vert<R}\vert (1+\Vert x\Vert)^d
V_h(x)\cdot w(x)\vert<\epsilon$ as well, by taking $\Vert h \Vert$ small enough.
This proves \eqref{eq:subaim} in the case $\vert\nu\vert=0$.
\item Assume $\vert\nu\vert=1$. Then the only possibility for $\partial_x^\nu$, which does not
lead to the trivial situation $\partial_x^\nu V_h= 0$, is
$\partial_x^\nu=(\partial/\partial x_j)$, in which case one verifies that, for
$h\neq 0$,
\begin{equation*}
\left\Vert \left(1+\Vert x\Vert\right)^d \partial_x^\nu V_h\cdot w\right\Vert_\infty\leq\vert h\vert\,\left\Vert
\left(1+\Vert x\Vert\right)^d \Vert x\Vert\cdot w\right\Vert_\infty.
\end{equation*}
This proves \eqref{eq:subaim} for $\vert\nu\vert=1$.
\item Assume $\vert\nu\vert\geq 2$. As in the previous case, the only possibility for $\partial_x^\nu$, which does not
lead to the trivial situation $\partial_x^\nu V_h= 0$, is
$\partial_x^\nu=(\partial/\partial x_j)^{\vert\nu\vert} $, in which case one
verifies that, for $h\neq 0$,
\begin{equation*}
\left\Vert \left(1+\Vert x\Vert\right)^d \partial_x^\nu V_h\cdot w\right\Vert_\infty\leq\vert
h\vert^{\vert\nu\vert-1}\left\Vert
\left(1+\Vert x\Vert\right)^d \cdot w\right\Vert_\infty.
\end{equation*}
This proves \eqref{eq:subaim} for $\vert\nu\vert\geq 2$.
\end{enumerate}
This concludes the proof for $w\notin\weightshol$.
\end{proof}

The following corollary follows by induction, using Propositions
\ref{prop:mapiscontinuous} and \ref{prop:emapsmooth}.

\begin{corollary}\label{cor:emapsmooth}
Let $d\geq 0$ be an integer. If $w\in\weights_d$, then
$\emap_1:\domain\mapsto\opw$ is strongly of class $C^{d}$, with derivatives
given by
\begin{equation*}
\left(\partial^\alpha\emap_1(\lambda)\right)(x)=i^{\vert\alpha\vert}x^\alpha e_{i\lambda}(x)\quad(\lambda\in\domain,\,x\in\Rn; \alpha\in\Nn,\,\vert\alpha\vert\leq
d).
\end{equation*}
If $w\in\weightshol$, then the derivatives exist for all $\alpha\in\Nn$ and are
to be interpreted as complex derivatives.
\end{corollary}

As a last auxiliary result, we need the following estimates:

\begin{lemma}\label{lem:estimatederivative} Let $\alpha\in\Nn$. If $w\in\weights_{\vert\alpha\vert}$ then, for all $\beta\in\Nn$ and
$\lambda\in\domain$,
\begin{equation*}
p_\beta\left(x^\alpha e^{i(\lambda,x)}\right)\leq
2^{\vert\beta\vert}\left(1+\Vert\lambda\Vert\right)^{\vert\beta\vert}
\prod_{j=1}^n (1+\alpha_j)^{\vert\beta\vert}\cdot\left\Vert\prod_{j=1}^n
(1+\vert x_j\vert^{\alpha_j})\cdot e^{-(\ima
\lambda,x)}\cdot w(x)\right\Vert_\infty.
\end{equation*}
\end{lemma}

\begin{proof}
It is evident that
\begin{equation}\label{eq:estimateseminorm}
p_\beta\left(x^\alpha
e^{i(\lambda,x)}\right)\leq\sum_{\mu+\nu=\beta}c_{\mu,\nu}\left\Vert\partial_x^\mu
x^\alpha\cdot\partial_x^\nu e^{i(\lambda,x)}\cdot w(x)\right\Vert_\infty
\end{equation}
for some nonnegative integers $c_{\mu,\nu}$ satisfying
\begin{equation}\label{eq:sumconstants}
\sum_{\mu+\nu=\beta}c_{\mu,\nu}=2^{\vert\beta\vert}.
\end{equation}
We estimate the functions occurring in the summands in
\eqref{eq:estimateseminorm}. First, in such a summand one has, for all
$x$,
\begin{align*}
\left\vert\partial_x^\mu
x^\alpha\right\vert&\leq\prod_{j=1}^n(1+\alpha_j)^{\beta_j}(1+\vert
x_j\vert^{\alpha_j})\\ &\leq \prod_{j=1}^n(1+\alpha_j)^{\vert\beta\vert}(1+\vert
x_j\vert^{\alpha_j}).
\end{align*}
In addition, in such a summand one has, for all $x$ and $\lambda$,
\begin{align*}
\vert\partial_x^\nu
e^{i(\lambda,x)}\vert&\leq\prod_{j=1}^n(1+\vert\lambda_j\vert)^{\beta_j}\cdot
e^{-(\ima\lambda,x)}\\
&\leq\left(1+\Vert\lambda\Vert\right)^{\vert\beta\vert}e^{-(\ima\lambda,x)}.
\end{align*}
The statement now follows from an application of the estimates in
\eqref{eq:estimateseminorm} and a subsequent use of \eqref{eq:sumconstants}.
\end{proof}

We are now in the position to prove the basic result for the Fourier transform
on $\opwdual$.

\begin{theorem}\label{thm:Fouriertransforms}
Let $d\in\{0,1,2,\ldots,\infty\}$, assume $w\in\weights_d$, and suppose
$\typ\in\opwdual$. Put $N=\ord_{\opwdual}\typ$. Then:
\begin{enumerate}
\item There is a natural injective map $\opwdual\rightarrow\tempdist$, defined by letting to $T\in\opwdual$
correspond the tempered distribution $f\mapsto\langle \typ,f\rangle$
\textup{(}$f\in\schw$\textup{)}, which we again denote by $T$. The order of $T$ as a tempered
distribution is at most $N$.
\item The Fourier transform of $\typ$, where $T$ is viewed as a tempered distribution, is a continuous function $\widehat \typ$. It is given for $\lambda\in\Rn$ by
\begin{equation}\label{eq:Fourierfunction}
\widehat \typ(\lambda)=(2\pi)^{-n/2}\langle \typ,e_{-i\lambda}\rangle.
\end{equation}
The map from $\opwdual$ into $C(\Rn)$ defined by $\typ\mapsto\Fouriertyp$ is
injective.
\item $\widehat \typ$ is of class $C^{d}$ on $\Rn$, with derivatives for $\lambda\in\Rn$
and $\alpha\in\Nn\,(\vert\alpha\vert\leq d)$ given by
\begin{equation}
\parta\Fouriertyp(\lambda)=(2\pi)^{-n/2}(-i)^{\vert\alpha\vert}\langle \typ,x^\alpha
e_{-i\lambda}\rangle.\label{eq:Fourierderivatives}
\end{equation}
\item There exists a constant $C>0$ such that, for all $\lambda\in\Rn$ and for all $\alpha\in\Nn$ with
$\vert\alpha\vert\leq d$,
\begin{equation*}
\left\vert\parta\Fouriertyp(\lambda)\right\vert\leq C (1+\Vert\lambda\Vert)^N \prod_{j=1}^n (1+
\alpha_j)^N \left\Vert\prod_{j=1}^n (1+\vert
x_j\vert^{\alpha_j})\cdot w(x)\right\Vert_\infty.
\end{equation*}
\item If $w\in\weightshol$, then \eqref{eq:Fourierfunction} defines a holomorphic extension
of $\widehat\typ$ to $\domain$, also denoted by $\widehat\typ$. The complex
derivatives of $\widehat \typ$ on $\domain$ are then given by
\eqref{eq:Fourierderivatives} and there exists a constant $C>0$ such that, for all $\lambda\in\domain$ and all $\alpha\in\Nn$,
\begin{equation*}
\left\vert\parta\Fouriertyp(\lambda)\right\vert\leq
C (1+\Vert\lambda\Vert)^N \prod_{j=1}^n (1+
\alpha_j)^N \left\Vert\prod_{j=1}^n (1+\vert
x_j\vert^{\alpha_j})\cdot e^{(\ima
\lambda,x)}\cdot w(x)\right\Vert_\infty .
\end{equation*}
\item If $w$ is standard quasi-analytic, then there exist a constant $C>0$ and nonnegative
sequences $\{\widetilde M(j,m)\}_{m=0}^\infty\,(j=1,\ldots,n)$ such
that
\begin{equation*}
\sum_{m=1}^\infty \frac{1}{\sqrt[m]{\widetilde M(j,m)}}=\infty\quad(j=1,\ldots,n)
\end{equation*}
and, for all $\lambda\in\Rn$ and all $\alpha\in\Nn$,
\begin{equation*}
\left\vert\parta\Fouriertyp(\lambda)\right\vert\leq C
(1+\Vert\lambda\Vert)^N \prod_{j=1}^n \widetilde M(j,\alpha_j).
\end{equation*}
Moreover, if $w$ is standard quasi-analytic and if there exists
$\lambda_0\in\Rn$ such that $\parta
\Fouriertyp(\lambda_0)=0$ for all $\alpha\in\Nn$, then $\typ=0$.
\end{enumerate}
\end{theorem}

\begin{proof}
As to (1): this was already observed following Proposition
\ref{prop:testdense}.

For (2), we first consider as preparation a weight $\wtilde\in\weights_0$
which is separated from zero on compact sets. Let $f\in\test$ be arbitrary. The
map $\lambda\mapsto f(\lambda)e_{-i\lambda}=f(\lambda)\emap_1(-\lambda)$ is then
continuous by Proposition~\ref{prop:mapiscontinuous}, evidently compactly
supported and with values in $\opwtilde$, which is a Fr\'echet space by
Proposition~\ref{prop:spaceisFrechet}. Therefore, as a consequence of
\cite[Theorem 3.27]{RudinFA}, the weak integral
\begin{equation*}
(2\pi)^{-n/2}\int_\Rn f(\lambda)e_{-i\lambda}\,d\lambda
\end{equation*}
exists in $\opwtilde$. Applying point evaluations, which are continuous since
$\wtilde$ has no zeros, the integral is identified as $\widehat f$. We now turn
to the statements in (2) and note that the function $\Fouriertyp$ as
described in the theorem is certainly continuous on $\Rn$ (in fact on $\domain$)
as a consequence of Proposition~\ref{prop:mapiscontinuous}. We also note that,
since $w\in\weights_d\subset\weights_0$, Theorem~\ref{thm:improvements} allows
us to choose a majorant $\wtilde\in\weights_0$ of $w$ which is separated from
zero on compact sets. The previous result on the weak integrals in $\opwtilde$
then imply, in view of the continuity of the injection
$\opwtilde\hookrightarrow\opw$, that
\begin{equation*}
\langle \typ,\widehat f\rangle=(2\pi)^{-n/2}\int_\Rn f(\lambda)\langle
\typ, e_{-i\lambda}\rangle\,d\lambda
\end{equation*}
for all $f\in\test$. Therefore, with $\Fourier(\typ)$ denoting the Fourier
transform of the tempered distribution $\typ$, we have
\begin{equation}\label{eq:Fourierint}
\langle \Fourier(\typ),f\rangle=\int_\Rn
f(\lambda)\Fouriertyp(\lambda)\,d\lambda
\end{equation}
for all $f\in\test$. Now note that as a consequence of
\eqref{eq:continuity} and the estimates in Lemma~\ref{lem:estimatederivative}
(applied with $\alpha=0$), $\Fouriertyp$ is of at most polynomial growth on
$\Rn$. In view of the density of $\test$ in $\schw$ and the continuity of
$\Fourier(\typ)$, the dominated convergence theorem implies that
\eqref{eq:Fourierint} actually holds for all $f\in\schw$, showing that the function $\Fouriertyp$ represents the Fourier transform of the tempered distribution $T$ as claimed.
If $\Fouriertyp=0$, then $\typ=0$ in $\tempdist$, implying $\typ=0$ in $\opwdual$
in view of (1). This concludes the proof (2).

Corollary~\ref{cor:emapsmooth} implies the degree of differentiability and the
expressions for the derivatives in (3) and (5).

In view of \eqref{eq:continuity}, (4) and the estimates in (5) follow
immediately from Lemma~\ref{lem:estimatederivative}.

We now prove the assertions on the sequences in (6) as a consequence of the
estimates in (4), as follows. For $j=1,\ldots,n$ and $m\geq 0$ put
\begin{equation*}
\widetilde M(j,m)=\left(1+m\right)^N \left\Vert\ (1+\vert x_j\vert^{m})\cdot
w^{1/n}(x)\right\Vert_\infty.
\end{equation*}
Since $w^{1/n}$ is a standard quasi-analytic weight by Lemma
\ref{prop:rootsareqa}, we have
\begin{equation}\label{eq:tempseries1}
\sum_{m=1}^\infty \Vert x_j^m \cdot
w^{1/n}\Vert_\infty^{-1/m}=\infty\quad(j=1,\ldots,n).
\end{equation}
Fix $j$. If there does not exists $x\in\Rn$ with $\vert x_j\vert>1$ such that
$w(x)>0$, then certainly
\begin{equation}\label{eq:tempseries2}
\sum_{m=1}^\infty \frac{1}{\sqrt[m]{\widetilde
M(j,m)}}=\infty\quad(j=1,\ldots,n),
\end{equation}
since the terms tend to a constant. In the remaining case where there exists
$x\in\Rn$ with $\vert x_j\vert>1$ such that $w(x)>0$, one verifies that the
quotient of the terms in the series in \eqref{eq:tempseries1} and
\eqref{eq:tempseries2} tends to $1$ as $m\rightarrow\infty$; \eqref{eq:tempseries2}
therefore follows from \eqref{eq:tempseries1}.

Finally, suppose that all derivatives of $\Fouriertyp$ vanish at some
$\lambda_0$. Then it follows from the estimates on the derivatives and Theorem
\ref{thm:DCtheorem} on quasi-analytic classes that the translate of
$\Fouriertyp$ over $\lambda_0$ is identically zero. Therefore $\Fouriertyp=0$,
implying $\typ=0$ (2).
\end{proof}

\begin{remark}
Theorem \ref{thm:Fouriertransforms}(5), when specialized to $\alpha=0$, is one-half of a
Paley--Wiener-type theorem. For a distribution $T$ with compact convex support
$K$ one can in fact derive the easy half of the geometric Paley--Wiener theorem
(cf.\ \cite[Exercise~29.5]{Treves}) from this specialization of (5), since such
$T$ is continuous on $\test$ in the induced topology of $\op_{\charfunc_K}$ by
\cite[p.~99, Theorem XXXIV.1]{Schwartz}.
\end{remark}

\subsection{Admissible spaces and their density in $\opw$}\label{subsec:admissible}

It is now easy to derive various density results from Theorem
\ref{thm:Fouriertransforms}. As may have become apparent from Section
\ref{sec:introduction}, in applications the relevant property of a dense
subspace of $\opw$ is not the density itself, but the property that it has the
same closure as any other dense subspace of $\opw$, and $\test$ in particular.
This leads to the concept of admissible spaces for weights (Definition
\ref{def:admissiblespaces}) and motivates the formulation of Theorem
\ref{thm:admissiblespaces}.

\begin{definition}\label{def:admissiblespaces} Let $w$ be a weight. An \emph{admissible space for $w$} is any
of the following subspaces of $C^\infty(\Rn)$:
\begin{enumerate}
\item If $w\notin\weights_0$: a dense subspace of $\schw$.

\item If $w\in\weights_0$ but $w\notin\weightsqa$:
\begin{enumerate}
\item a dense subspace of $\schw$; or
\item a subspace of the form $\Span\{\,e_{i\lambda}\mid \lambda\in E\,\}$, where $E\subset\Rn$ is dense.
\end{enumerate}

\item If $w\in\weightsqa$ but $w\notin\weightshol$:
\begin{enumerate}
\item a dense subspace of $\schw$; or
\item a subspace of the form $\Span\{\,e_{i\lambda}\mid \lambda\in E\}$, where
$E\subset\Rn$ is such that $\overline E$ has nonempty interior; or
\item $\poly$.
\end{enumerate}

\item If $w\in\weightshol$:
\begin{enumerate}
\item a dense subspace of $\schw$; or
\item a subspace of the form $\Span\{\,e_{i\lambda}\mid \lambda\in
E\,\}$, where $E\subset\domain$ is such that any holomorphic function on
$\domain$ vanishing on $E$ is identically zero on $\domain$; or
\item $\poly$.
\end{enumerate}

\end{enumerate}
\end{definition}

\begin{remark}
In the one-dimensional case, the possible choices for spectral subsets $E$ as
required above for $w\in\weightshol$ can alternatively be described as the
subsets of $\C_w$ having an accumulation point in $\C_w$. In arbitrary
dimension, the situation is more involved and a simple description of such
uniqueness sets is not known. $\Rn$ is a uniqueness set for holomorphic
functions on $\Cn$ and, furthermore, the boundary uniqueness theorem and the
maximum principle \cite[p.~289]{Chirka} supply nontrivial examples of
uniqueness sets for holomorphic functions on $\domain$. A somewhere dense subset
of $\Rn$ is a uniqueness set for holomorphic functions on $\domain$.
\end{remark}

\begin{theorem}\label{thm:admissiblespaces}
Let $w$ be a weight. Then all admissible spaces for $w$ are subspaces of $\opw$
and they all have the same closure in $\opw$. In fact, they are all dense.
\end{theorem}

\begin{proof}
It is obvious that the admissible spaces are subspaces. To prove that they are
all dense, we start with an arbitrary weight such that $w\notin\weights_0$ and
then prove the various incremental statements on density as the properties of
$w$ improve.

In the starting case, the density in $\opw$ of an arbitrary dense subspace of
$\schw$ follows from Lemma~\ref{lem:fundamentallemma}, the density of $\schw$ in
$\opw$, and the continuity of the injection $\schw\hookrightarrow\opw$.

In the case $w\in\weights_0$ but $w\notin\weightsqa$, suppose that $\typ\in\opwdual$
vanishes on the span of the exponentials as described in Definition
\ref{def:admissiblespaces}. By the continuity of $\Fouriertyp$ we conclude that
$\Fouriertyp=0$, hence $\typ=0$.

In the case $w\in\weightsqa$ but $w\notin\weightshol$, then as a consequence of
Proposition~\ref{prop:invariance} we may assume that $w$ is standard
quasi-analytic. In that case, if $\typ\in\opwdual$ vanishes on the span of the
exponentials as described in Definition~\ref{def:admissiblespaces}, then by
continuity there exists an open set on which $\Fouriertyp$ vanishes. In view of
Theorem \ref{thm:Fouriertransforms}(6) we conclude that $\typ=0$. Also,
if $\typ\in\opwdual$ vanishes on $\poly$, then Theorem~\ref{thm:Fouriertransforms}(3) and (6) show that $\typ=0$.

Finally, in the case $w\in\weightshol$ the remaining statement on the density of the
span of the exponentials as described follows from the holomorphy of
$\Fouriertyp$ on $\domain$ and again the injectivity of the map $\typ\mapsto\Fouriertyp$.
\end{proof}

\begin{remark}\label{rem:generalizations}
\begin{enumerate}
\item A crucial step in the proof of Theorem
\ref{thm:admissiblespaces} is the conclusion that $\typ$ vanishes on $\opw$
whenever the function $\Fouriertyp$---which could as well have been defined
directly as $\Fouriertyp=\langle T,e_{-i\lambda}\rangle$---vanishes on $\Rn$.
For this step it is not essential that $\Fouriertyp$ can be interpreted as a
Fourier transform of a tempered distribution. Indeed, the conclusion follows by
considering a weak integral $\int_\Rn f(\lambda) e_{-i\lambda}\,d\lambda$ in
$\opw$ as in the proof of Theorem~\ref{thm:Fouriertransforms} when one uses that
$\Fourier(\test)$ is dense in $\schw$ and therefore also in $\opw$. A range
result of this kind appears to be essential; the closure theorems
on exponentials and polynomials in $\opw$ to a considerable extent ultimately rest on
harmonic analysis. This connection with harmonic analysis will become even more
apparent in the discussion of possible generalizations in
Section~\ref{sec:interpretation}.
\item It is instructive to note that the various density results in $\opw$ have different origins. The density of $\test$ and $\schw$ for general $w$ can be
proved in a direct fashion. The density of the trigonometric functions for
$w\in\weights_0$ has, as noted above, its origin in harmonic analysis. For
$w\in\weightshol$ this can still be said for the additional density of $\poly$
and the span of exponentials as described. For $w\in\weightsqa$ but $w\notin\weightshol$ both the density of $\poly$ and the density of the span of
exponentials as described are relatively deep statements, based not only on a
result from harmonic analysis but on a theorem on quasi-analytic classes in
several variables as well.
\item The theorem can be generalized (cf.\ Theorem~\ref{thm:Bernsteinspaces}).
\end{enumerate}
\end{remark}

For $w\in\weights_\infty$ we have established that the condition
$w\in\weightsqa$ is sufficient for the polynomials to be dense in $\opw$. The
following proposition, which should be compared with Corollary~\ref{cor:online},
gives a necessary condition. It concludes our study of density in $\opw$.

\begin{proposition}\label{prop:hall}
Let $w\in\weights_\infty$ and suppose that $\poly$ is dense in $\opw$. If
$x,y\in\Rn$ are such that $y\neq 0$ and such that $t\mapsto w(x+ty)$ is Lebesgue
measurable on $\R$, then
\begin{equation}\label{eq:hall}
\int_{-\infty}^\infty\frac{\log w(x+ty)}{1+t^2}=-\infty.
\end{equation}
\end{proposition}

\begin{proof}
Fix $x$ and $y$ as in the hypotheses. Define $\iota:\R\mapsto\R^n$ by
$\iota(t)=x+ty$. If $f$ is a function on $\Rn$ let $\iota^*f=f\circ\iota$ be its
pull-back to $\R$. Then $\iota^*w$ is a rapidly decreasing weight on $\R$.
Consider the one-dimensional Bernstein space $C_{\iota^*w}=\{\,f\in C(\R)\mid
\lim_{\vert x\vert\rightarrow\infty}f(x)\cdot
\iota^*w(x)=0\,\}$, supplied with the topology given by the seminorm $f\mapsto\Vert f\cdot
\iota^*w\Vert_\infty$. Then $\iota^*$ maps $\opw$ continuously into
$C_{\iota^*w}$. Furthermore, $\iota^*\poly(\Rn)=\poly(\R)$ and
$\iota^*\test(\Rn)=\test(\R)$. Now if $\poly(\Rn)$ is dense in $\opw$, i.e., if
$\poly(\Rn)$ and $\test(\Rn)$ have the same closure in $\opw$, then we conclude
from this and Lemma~\ref{lem:fundamentallemma} that $\poly(\R)$ and $\test(\R)$
have the same closure in $C_{\iota^*w}$. Since $\test(\R)$ is dense in
$C_{\iota^*w}$, so is $\poly(\R)$. Therefore, the weight $\iota^*w$ must satisfy
Hall's integral condition \cite{Hall,Mergelyan}, i.e., \eqref{eq:hall} must hold.
\end{proof}

\section{Subspaces with equal closure}\label{sec:general}

We are now able to prove the basic theorems on subspaces with equal closure.
These are supplemented with a discussion of a natural partial filtration on
topological vector spaces. The latter is motivated by the question what the
amount of freedom is to obtain variations on $\eqref{eq:translation}$, as
depending on the generator, and leads to the concept of admissible spaces for
elements of a topological vector space.

For ease of formulation we limit ourselves to locally convex spaces. In view of
practice, this is not too severe a limitation and if necessary it is, on the basis
of Theorem~\ref{thm:closurecont}, easy to adjust the results to the situation of
finitely generated subspaces of spaces which are not locally convex.

\subsection{Formalization of the method}\label{subsec:formalization}

The following basic theorem is immediate from the combination of Theorems
\ref{thm:closurecont} and \ref{thm:admissiblespaces}. The admissible
spaces for weights in its formulation can be found in Definition
\ref{def:admissiblespaces}.

\begin{theorem}[Subspaces with equal closure]\label{thm:generalthm}
Let $V$ be a locally convex topological vector space and let $\gen\subset V$ be
a nonempty set of generators. Suppose that for each $g\in\gen$ a weight $w_g$
on $\R^{n_g}$ is given, together with a continuous linear map from $\op_{w_g}$
into $V$ denoted by $f\mapsto f\cdot g$ $(f\in\op_{w_g})$. For each $g\in\gen$
choose an admissible space $L_g$ for $w_g$ and consider
$L=\sum_{g\in\gen}L_g\cdot g$. Then:
\begin{enumerate}
\item The closure of $L$ is independent of the choice of the $L_g$.
\item The annihilator of $L$ is independent of the choice of the $L_g$ and
is equal to $\bigcap_{g\in\gen} \left( L_g^\prime\cdot g\right)^\perp$ for any
choice of admissible spaces $L_g^\prime$.
\item The sequential closure of $L$ is independent of the choice of the $L_g$ when $\gen$ is finite.
\end{enumerate}
\end{theorem}

In the applications we have in mind $V$ is a space of functions, equivalence
classes of functions, distributions, sequences, etc., with a subset of $\Rn$ as
the underlying point set. The $n_g$ are then all equal to $n$ and we will assume
this for the remainder of the discussion. The linear map in the above theorem
is, in those situations, the appropriate type of multiplication by smooth
functions, as is defined in $V$. As the reader may verify in examples, this
multiplication invariably results in a natural module structure over $\poly$,
$\galg$, and $\test$ for some---in general proper---subspaces of $V$.

The salient point in the above theorem is the possibility to switch from one
choice of admissible spaces to another more convenient one (and, in particular, to
the choice of $\test$ for all $g\in\gen$) without affecting the closure of the
subspace involved. As explained in Section~\ref{sec:introduction}, this is the
mechanism which makes the nontrivial analysis in the spaces $\opw$ immediately
available in other situations.

We briefly summarize the switching possibilities. If $w_g\in\weights$, then one
can choose as admissible space $L_g$ a dense subspace of $\schw$ (such as
notably $\test$). If $w_g\in\weights_0$, then in addition one can choose a space
of the form $\Span\{\,e_{i\lambda}\mid\lambda\in E\,\}$, where $E\subset\Rn$ is
dense. If $w_g\in\weightsqa$, then there is still more choice since for $L_g$
one may now in addition also choose $\poly$ or
$\Span\{\,e_{i\lambda}\mid\lambda\in E\,\}$ if $E$ is a somewhere dense subset
of $\Rn$. If $w_g\in\weightshol$, then maximum choice is possible since, in addition to
all spaces just described, there is then another type of admissible space,
related to the exponentials. This additional degree of freedom depends on more
details of the holomorphic weight and we refer to
Definition~\ref{def:admissiblespaces} for the precise description.

The choices for admissible spaces can be made independently for all $g\in\gen$
and all possible resulting subspaces $V$ have the same closure. An application
of the theorem could in many situations, e.g., be the conclusion that the modules
generated by $\gen$ over $\test$, $\galg$, and $\poly$, all three have the same
closure, and that for a description of this closure---equivalently, a
description of the annihilator---one can assume that the module is generated
over $\test$, which in most situations is a convenient choice. It is apparent
that such a conclusion is a special case of the assertions in the theorem.

As an illustration we give some straightforward examples,
returning to more elaborate situations in Section~\ref{sec:applications}.

\begin{example}
\begin{enumerate}
\item Let $\Omega\subset\Rn$ be open and nonempty and let $V=\distOmega$. Supply $\distOmega$ with its weak topology and consider
$g=1\in\distOmega$. For the holomorphic weight $w(x)=\exp(-\Vert x\Vert^2)$ the
map $f\mapsto f\cdot g$ is continuous from $\opw$ into $\distOmega$. Since
$\test\cdot1\supset\testcOmega$ and the latter space is well-known to be sequentially
dense in $\distOmega$, the same holds for $\test\cdot 1$. The third part of the
theorem then asserts that this sequential density also holds for $\poly\cdot 1$
and any space of the form $\Span\{\,e_{i\lambda}\cdot 1\mid\lambda\in E\,\}$,
where $E\subset\C^n$ is a uniqueness set for holomorphic functions on $\Cn$.

\item Let $V=\schw$ and let $g(x)=\exp(-\Vert x\Vert^2)\in\schw$.
For the holomorphic weight $w(x)=\exp(-\Vert x\Vert^2)$ the map $f\mapsto f\cdot
g$ is continuous from $\opw$ into $\schw$. Since $\test\cdot\exp(-\Vert
x\Vert^2)=\test$ is dense in $\schw$, we conclude from the theorem
that $\poly\cdot\exp(-\Vert x\Vert^2)$ is also dense, as is
$\Span\{\,e_{i\lambda}\cdot\exp(-\Vert x\Vert^2)\mid
\lambda\in E\,\}$, where $E\subset\C^n$ is a uniqueness set for holomorphic
functions on $\Cn$.
\item Let $V=L_p(\Rn,dx)$ $(1\leq p<\infty)$ and consider $g(x)\in L_p(\Rn,dx)$ given by $g(x)=\exp(-\sqrt{\Vert x\Vert^2+1}/\log \sqrt{\Vert x\Vert^2 +2})$. Let $w$ be a weight
which is equal to $\exp (-\epsilon \Vert x\Vert/2\log \Vert x\Vert)$ for large
$\Vert x\Vert$. Then $w$ is quasi-analytic by Proposition \ref{prop:examplesqa};
the associated map $f\mapsto f\cdot g$ from $\opw$ into $V=L_p(\Rn,dx)$ is then
continuous. Since $\test\cdot g=\test$ is dense in $L_p(\Rn,dx)$, we conclude
from the theorem that $\poly\cdot g$ is dense in $L_p(\Rn,dx)$, as is
$\Span\{\,e_{i\lambda}\cdot g\mid\lambda\in E\}$ if $E\subset\Rn$ is somewhere
dense. We will later consider the case of nonsmooth $g$ (see
Section~\ref{subsec:locallycompact}).
\end{enumerate}
It will become apparent in Theorem~\ref{thm:spacedist} that the continuity of
the maps of type $f\rightarrow f\cdot g$ from $\opw$ into $V$ in these examples
does not have to be verified. It is sufficient to observe that these maps can be
defined.
\end{example}

For many spaces $V$ the natural action of a space $\opw$ on $g\in V$ can be
described as the transpose of an action of $\opw$ on functions on some set $X$,
the latter being defined by pointwise multiplication. In such situations a
uniform formulation becomes possible, covering in one concept a variety of
``appropriate'' ways of multiplying objects by smooth functions. At the same
time a new phenomenon appears, namely the redundancy of the hypothesis of
continuity. This is due to the completeness of various spaces and it illustrates
the importance of Theorem~\ref{thm:improvements}.

\begin{theorem}\label{thm:spacedist}
Let $X\subset\Rn$ be nonempty and let $\testX$ be an LF-space of functions on
$X$ with a defining sequence $\testX^1\subset\testX^2\subset\ldots$. Assume that
the $\testX^l$ are all $C^\infty(\Rn)$-modules under pointwise multiplication
and suppose that convergence of a sequence in $\testX$ implies pointwise convergence on $X$.
Let $\testXdual$ be the dual of $\testX$, supplied with the weak topology.

Then the action of each element of $C^\infty(\Rn)$ on $\testX$ is automatically
continuous, and we equip $\testXdual$ with the structure of a
$C^\infty(\Rn)$-module by transposition.

Let $V\subset\testXdual$ be a locally convex topological vector space and let
$\gen\subset V$ be a nonempty set of generators. Suppose that for each
$g\in\gen$ a weight $w_g$ on $\Rn$ is given. Assume
\begin{enumerate}
\item that the topology of $V$ is weaker than the induced weak topology of $\testXdual$ on $V$ and that $\op_{w_g}\cdot g\subset V$ for all $g\in\gen$; or
\item that the topology of $V$ is stronger than the induced weak topology of $\testXdual$ on
$V$, that $V$ is an LF-space with a defining sequence $V^1\subset
V^2\subset\ldots$, and that for each $g\in\gen$ there exists $V^l$ such that
$g\in V^l$ and $\op_{w_g}\cdot g\subset V^l$.
\end{enumerate}
For each $g\in\gen$ choose an admissible space $L_g$ for $w_g$. Consider the
subspace $L=\sum_{g\in\gen}L_g\cdot g$ of $V$. Then the following hold:
\begin{enumerate}
\item The closure of $L$ in $V$ is independent of the choice of the $L_g$.
\item The annihilator of $L$ in $V^\prime$ is independent of the choice of the $L_g$ and
equal to $\bigcap_{g\in\gen} \left( L_g^\prime\cdot g\right)^\perp$ for any
choice of admissible spaces $L_g^\prime$.
\item The sequential closure of $L$ in $V$ is independent of the choice of the $L_g$ when $\gen$ is finite.
\end{enumerate}
\end{theorem}

Before giving the proof of the theorem, let us motivate its formulation by
indicating what types of spaces can occur as $\testX$ and $V$ in the theorem.

The hypotheses on the space of test functions $\testX$ are for all nonempty
open $\Omega\subset\Rn$ and all $0\leq m\leq\infty$ satisfied by
$C_c^{m}(\Omega)$ and $C^{m}(\Omega)$, the latter case corresponding to a
constant sequence defining $\testX$. Thus a number of common dual spaces can
figure as $\testXdual$ in the theorem and for each of these one has a variety of
subspaces $V\subset\testXdual$. As an example we consider
$\testXdual=\distOmega$. Then $\distOmega$ itself in the weak topology is
covered by the first possibility. Taking the defining sequence in the second
possibility equal to a fixed Fr\'echet space, we see that Fr\'echet spaces as
$\testOmega$, $L_p$-spaces (as well as their local versions), various types of
Sobolev spaces and $\schw$ (in the case $\Omega=\Rn$) are all covered by this second
possibility. Other spaces which are not Fr\'echet spaces but still LF-spaces are
also covered by this second possibility. This holds, e.g., for certain Sobolev
spaces and for $\testcOmega$. We will see in Section~\ref{subsec:locallycompact}
that several spaces associated to a locally compact subset $X$ of $\Rn$, and
notably $L_p$-spaces with $X$ as an underlying point set, are also within the range
of the theorem.

In all such spaces an inclusion of sets immediately implies equality of closure
of subspaces.

\begin{proof}
We first settle the continuity of the action of $C^\infty(\Rn)$. Since
convergence of a sequence in $\testX$ implies pointwise convergence on $X$, the same holds for each space $\testX^l$. The closed graph
theorem \cite[Theorem 2.15]{RudinFA} then readily shows that the action of each fixed element of
$C^\infty(\Rn)$ on each space $\testX^l$ is continuous as a map from $\testX^l$
into itself. Composed with the continuous injection
$\testX^l\hookrightarrow\testX$ the action is a continuous map from each space
$\testX^l$ into $\testX$. This in turn implies by general principles the
continuity of the action as a map from $\testX$ into itself, as was to be shown.
We conclude that there is indeed a $C^\infty(\Rn)$-module structure on
$\testXdual$.

As a preparation for the conclusion of the proof, we first improve the
properties of the $\op_{w_g}$. Indeed, Theorem~\ref{thm:improvements} supplies
majorants for the $w_g$ which still satisfy the hypotheses of the present
theorem (since the corresponding function spaces are contained in the spaces
$\op_{w_g}$) and which are of the same type, so that they have the same
admissible spaces as the original $w_g$. Replacing the $w_g$ by these majorants, we may therefore assume by Theorem~\ref{thm:improvements} that the $\op_{w_g}$
are Fr\'echet spaces in which convergence implies pointwise convergence on
$\Rn$.

After this preparation we note that in view of Theorem~\ref{thm:generalthm} the
conclusions of the present theorem will follow once we can prove that $f\mapsto
f\cdot g$ is for each $g\in\gen$ automatically continuous as a map from
$\op_{w_g}$ into $V$. To this end fix $g\in\gen$.

Since convergence in $\op_{w_g}$ now implies pointwise convergence and
$\op_{w_g}$ is now Fr\'echet, the closed graph theorem yields that, for each
fixed $t\in\testX$, say $t\in\testX^l$, the map $f\mapsto f\cdot t$ is
continuous from $\op_{w_g}$ into $\testX^l$ and is therefore also continuous as
a map from $\op_{w_g}$ into $\testX$. By transposition we conclude that
$f_n\cdot g\overset{w}{\rightarrow}f_\infty\cdot g$ in the weak topology of
$\testXdual$ whenever $f_n\rightarrow f_\infty$ in $\op_{w_g}$.

The required continuity in the case of the first possibility for the topology of $V$
is now immediate.

Turning to the second possibility, suppose that $g\in V^l$ and $\op_{w_g}\cdot
g\subset V^l$. By continuity of the injection $V^l\hookrightarrow V$ it is
sufficient to show that the map $f\mapsto f\cdot g$ is continuous as a map from
$\op_{w_g}$ into $V^l$. This can be done by a final application of the closed
graph theorem. Suppose then that $f_n\rightarrow f_\infty$ in $\op_{w_g}$ and
that $f_n\cdot g\rightarrow t^\prime$ in $V^l$. We have to show that
$t^\prime=f_\infty\cdot g$. Now the second convergence also holds in $V$ and,
since by assumption the topology of $V$ is stronger than the induced weak
topology of $\testXdual$, we see that $f_n\cdot g\overset{w}{\rightarrow}
t^\prime$ in the weak topology of $\testXdual$. On the other hand, we had
already observed that $f_n\cdot g\overset{w}{\rightarrow}f_\infty\cdot g$ in the
weak topology of $\testXdual$ whenever $f_n\rightarrow f_\infty$ in $\op_{w_g}$.
By the uniqueness of weak limits in $\testXdual$, we conclude that
$t^\prime=f_\infty\cdot g$, as was to be shown. The closed graph theorem
therefore does indeed apply and the proof is complete.
\end{proof}

The following theorem has an overlap with Theorem~\ref{thm:spacedist}, but it
covers some important spaces like $C(X)$, $C_c(X)$ and $C_0(X)$ for a locally
compact subset $X$ of $\Rn$, spaces which are not, for general $X$, covered by the
previous theorem. The proof is a simpler version of the proof of Theorem
\ref{thm:spacedist} and is left to the reader.

\begin{theorem}\label{thm:spacefunctions}
Let $X\subset\Rn$ be nonempty and let $V$ be an LF-space of functions on
$X$ with a defining sequence $V^1\subset V^2\subset\ldots$. Suppose that convergence of a sequence in $V$ implies pointwise convergence on $X$.
Let $\gen\subset V$ be a nonempty set of generators and suppose that for each
$g\in\gen$ a weight $w_g$ on $\Rn$ is given.
Assume that for each $g\in\gen$ there exists $V^l$ such that $g\in V^l$ and $\op_{w_g}\cdot g\subset V^l$.
For each $g\in\gen$ choose an admissible space $L_g$ for $w_g$. Consider the
subspace $L=\sum_{g\in\gen}L_g\cdot g$ of $V$.

Then the following hold:
\begin{enumerate}
\item The closure of $L$ in $V$ is independent of the choice of the $L_g$.
\item The annihilator of $L$ in $V^\prime$ is independent of the choice of the $L_g$ and
equal to $\bigcap_{g\in\gen} \left( L_g^\prime\cdot g\right)^\perp$ for any
choice of admissible spaces $L_g^\prime$.
\item The sequential closure of $L$ in $V$ is independent of the choice of the $L_g$ when $\gen$ is finite.
\end{enumerate}
\end{theorem}

\begin{corollary}\label{cor:nogrowth}
Let $V$ be as in Theorem~\textup{\ref{thm:spacedist}} or~\textup{\ref{thm:spacefunctions}} and
suppose that $V$ \textup{(}or each space $V^l$ in the case of an LF-space\textup{)} is a
$C^\infty(\Rn)$-module. Let $\gen\subset V$ be a nonempty set of generators and,
for $g\in\gen$, let $L_g$ denote an admissible space for a holomorphic weight.

Then the closure in $V$ of the subspace $L=\sum_{g\in\gen}L_g\cdot g$ is
independent of the choice of the $L_g$.
\end{corollary}

The precise nature of $V$ and the generators is apparently irrelevant. This situation occurs, e.g., for $\testcmOmega$, $\testmOmega$ $(0\leq
m\leq\infty$) and their duals, and for $C_c(X)$ and $C(X)$ and their duals for
locally compact $X$. For these spaces one can, for arbitrary sets of generators,
always switch between admissible spaces for holomorphic weights without
affecting the closure.

\subsection{A partial filtration associated with types of weights}\label{subsec:filtration}

Given a topological vector space $V$ as in Theorem~\ref{thm:spacedist} or \ref{thm:spacefunctions} we now introduce a partial filtration on $V$
which measures, for $g\in V$, both the amount of choice for admissible spaces and
the minimal regularity of the Fourier transforms of associated linear
functionals on $\schw$.

For the first possibility of Theorem~\ref{thm:spacedist} let $V_{\weights}$
denote the set of all $v\in V$ such that $\opw\cdot v\subset V$ for some
$w\in\weights$. For $0\leq d\leq\infty$ let $V_d$ denote the set of all $v\in V$
such that $\opw\cdot v\subset V$ for some $w\in\weights_d$; the sets $V_{\textup{\text{qa}}}$
and $V_{\textup{\text{hol}}}$ are defined similarly. For the second possibility of Theorem
\ref{thm:spacedist}, as well as for Theorem~\ref{thm:spacefunctions}, let
$V_{\weights}$ denote the set of all $v\in V$ such that there exists $V^l$ with
the property that $v\in V^l$ and $\opw\cdot v\subset V^l$ for some
$w\in\weights$. The spaces $V_d$, $V_{\textup{\text{qa}}}$ and $V_{\textup{\text{hol}}}$ are defined similarly.

In all cases Proposition~\ref{prop:invariance} implies that
\begin{equation}\label{eq:filtration}
0\subset V_{\textup{\text{hol}}}\subset V_{\textup{\text{qa}}}\subset V_\infty\subset V_{d}\subset
V_{d^\prime}\subset V_{0}\subset V_{\weights}\subset V\quad(d\geq d^\prime\geq
0).
\end{equation}
One verifies that in all cases $V_{\textup{\text{hol}}}$, $V_d$ $(0\leq d\leq\infty$) and
$V_{\weights}$ are vector spaces since
$\op_{w_1+w_2}\subset\op_{w_1}\cap\op_{w_2}$ for arbitrary weights $w_1$ and
$w_2$ and since the sets of weights corresponding to these subspaces of $V$ are
closed under addition by Proposition~\ref{prop:invariance}. Since this latter
property fails for $\weightsqa$ it may well be in doubt whether $V_{\textup{\text{qa}}}$ is in
general closed under addition, so that we prefer to refer to
\eqref{eq:filtration} only as a partial filtration.

\begin{definition}\label{def:admissibleforelement}
Let $V$ be a topological vector space as in Theorems~\ref{thm:spacedist} or
\ref{thm:spacefunctions}. If $g\in V_{\textup{\text{hol}}}$, where $V_{\textup{\text{hol}}}$ is as defined
above, then an \emph{admissible space for $g$} is any admissible space for a
holomorphic weight as in Definition~\ref{def:admissiblespaces}. We define
analogously the notion of an admissible space for elements of $V_{\textup{\text{qa}}}$, $V_{0}$,
and $V_\weights$.
\end{definition}

With this terminology the essence of Theorem~\ref{thm:spacedist} and
\ref{thm:spacefunctions} can be rephrased as the validity of their conclusions
for any nonempty $\gen\subset V_\weights$ and any choice of admissible spaces
for the elements of $\gen$. The largest freedom of choice for a fixed $g\in\gen$
is then determined by the leftmost set in
\eqref{eq:filtration} that still contains $g$.

We turn to Fourier transforms. In the proof of Theorems~\ref{thm:spacedist}
and~\ref{thm:spacefunctions} we start with $g\in V$ and a weight $w_g$ such that
$\op_{w_g}\cdot g\subset V$ (resp.\ $g\in V^l$ and $\op_{w_g}\cdot g\subset V^l$
in the case of LF-spaces). The proof consists of replacing $w_g$ by a majorant
$\wtilde_g$ which is separated from zero on compact sets and then showing that
the map $f\mapsto f\cdot g$ from $\op_{\wtilde_g}$ into $V$ is necessarily
continuous. Now if, e.g., $g\in V_{\textup{\text{qa}}}$, so that $w_g\in\weightsqa$, then we can
choose $\wtilde_g\in\weightsqa$ again. Now for all $\typ\in V^\prime$ the
associated map $\typ_g:\schw\mapsto\C$ defined by
$\langle\typ_g,f\rangle=\langle\typ,f\cdot g\rangle$ $(f\in\schw)$ is not only a
tempered distribution but in fact in $\op_{\wtilde_g}^\prime$ since it factors
over $\op_{\wtilde_g}$. Since $\wtilde_g\in\weightsqa$ we therefore then have
knowledge from Theorem~\ref{thm:Fouriertransforms} on the minimal degree of
regularity of the Fourier transform of $\typ_g$. The same argument is valid for
the other classes of weights and the partial filtration \eqref{eq:filtration}
thus corresponds to various minimal degrees of regularity of the Fourier
transforms $\Fourier(\typ_g)$ of all these associated tempered distributions.
This minimal degree of regularity of all $\Fourier(\typ_g)$ for $T\in V^\prime$
is determined by, again, the leftmost set in \eqref{eq:filtration} that still
contains $g$.

To be more precise, fix $g\in V_\weights$. If $g\in V_d$ $(0\leq d\leq\infty)$
then the Fourier transforms $\Fourier(\typ_g)$ $(\typ\in V^\prime)$ of the
associated tempered distributions $\langle \typ_g,f\rangle=\langle T, f \cdot
g\rangle$ ($f\in\schw$) are all functions of class $C^{[d]}$ and are given by
$\widehat{\typ}_g(\lambda)=(2\pi)^{-n/2}\langle\typ,e_{-i\lambda}\cdot g\rangle$
$(\lambda\in\Rn)$. If $g\in V_{\textup{\text{qa}}}$ the Fourier transforms are each in a
quasi-analytic class and if $g\in V_{\textup{\text{hol}}}$ they all extend to holomorphic
functions on a common tubular neighborhood of $\Rn$. The expressions for the
possible derivatives of the $\Fourier(\typ_g)$, as well as estimates for these
derivatives, follow from those in Theorem~\ref{thm:Fouriertransforms}.

\section{Applications}\label{sec:applications}

In Section~\ref{sec:general} we have given some first applications of the spaces
$\opw$ in the case of one generator. We will now consider the case of arbitrary
sets of generators in a number of familiar spaces. The outline and the type of
results will be quite similar for the various spaces, in accordance with our aim
to develop a uniform approach.

With one exception we will use Theorems~\ref{thm:spacedist}, \ref{thm:spacefunctions}, and Corollary~\ref{cor:nogrowth} as the basic results,
so that we need not verify the continuity of maps of type $f\rightarrow f\cdot
g$ from a space $\opw$ into a space of interest. A direct verification of this
continuity is of course also possible. The exception is formed by the general
spaces of Bernstein type in Subsection~\ref{subsec:Bernsteinspaces}, a class of
spaces for which one can find examples that are not covered by the
aforementioned results.

With the exception of Subsection~\ref{subsec:Bernsteinspaces} we will use notations
like, e.g., $V_{\textup{\text{qa}}}$ as in the partial filtration in Subsection
\ref{subsec:filtration}. We will also employ the notion of admissible spaces for
elements of a topological vector space as in Definition
\ref{def:admissibleforelement}.

All subsets $\gen$ of generators are assumed to be nonempty. If $f$ is a
function, $Z(f)$ denotes its zero locus. We will also employ this notation for
elements of $L_p$-spaces, in which case the zero locus is not necessarily well
defined. Since the measure of such sets is all we will be concerned with in that
situation, we will allow ourselves this imprecision. The notations concerning
test functions and distributions are as in Section~\ref{sec:basics}. If
$y\in\Rn$ and $f$ is a function on $\Rn$, then the translate $T_y f$ is defined
by $T_y f(x)=f(x-y)$.

\begin{remark}
As explained in Subsection~\ref{subsec:filtration} we have in applications a
minimal degree of regularity of the Fourier transforms of associated tempered
distributions, as well as estimates on these transforms. We will refrain from
stating these results explicitly in the present section, since they are quite
similar in all cases and not necessary for our purposes. It should be noted,
however, that in other cases this extra information may be a useful
additional means of investigation.
\end{remark}

\subsection{The space $\schw$ of rapidly decreasing functions}\label{subsec:schwartz}
The Schwartz space $\schw$ is covered by the second possibility in
Theorem~\ref{thm:spacedist} by taking $\testX=\test$.

We start by discussing the partial filtration of Subsection
\ref{subsec:filtration}. It is immediate that $\schw=\schw_d$ for all finite $d$
(consider $w(x)=(1+\Vert x\Vert^{d+1})^{-1}$), so that $\galg$ and all dense
subspaces of $\schw$ are admissible spaces for all $f\in\schw$. For $\schw_{\textup{\text{qa}}}$
and $\schw_{\textup{\text{hol}}}$ the following lemma describes a class of examples:

\begin{lemma}
Let $g_1\in C^\infty(\Rn)$ and $g_2\in\schw$. Assume that $\vert
g_2\vert\in\weightsqa$ \textup{(}resp., $\vert g_2\vert\in\weightshol$\textup{)}. Suppose that for
all $\alpha,\,\beta\in\Nn$ there exist polynomials $P_\alpha,\,Q_\beta$ and
exponents $\mu_\alpha,\,\nu_\beta\geq 0$ such that:
\begin{enumerate}
\item
$\vert \parta g_1 (x)\vert\leq \vert P_\alpha(x)\vert \vert g_2\vert^{\mu_\alpha}$ for all $x\in\Rn$ with $\Vert x\Vert$ sufficiently large.
\item
$\vert \partial^\beta g_2(x)\vert\leq \vert Q_\beta(x)\vert \vert g_2\vert^{\nu_\beta}$ for all $x\in\Rn$ with $\Vert x\Vert$ sufficiently large.
\item $\min_{\alpha,\,\beta\in\Nn} (\mu_\alpha+\nu_\beta)>0$.
\end{enumerate}
Then $g=g_1g_2\in\schw_{\textup{\text{qa}}}$ \textup{(}resp., $g=g_1g_2\in\schw_{\textup{\text{hol}}}$\textup{)}. If $g_2$ has no
zeros for sufficiently large $\Vert x\Vert$, then the condition that the
$\mu_\alpha$ and $\nu_\beta$ are nonnegative is superfluous.
\end{lemma}

The straightforward proof, which we leave to the reader, consists of verifying
that $\opw\cdot g_1g_2\subset\schw$ for the quasi-analytic (resp., holomorphic)
weight $w=\vert
g_2\vert^{1/2\min_{\alpha,\,\beta\in\Nn}(\mu_\alpha+\nu_\beta)}$.

It is obvious from the lemma that $\exp(-\epsilon\Vert x\Vert^2)\in\schw_{\textup{\text{hol}}}$
for all $\epsilon>0$. The same is true for a more exotic example as $\sin
(\exp x)\,\cdot\,\exp(-\vert x\vert^2+\sin (\exp x))$ on the real line. For the
quasi-analytic case the lemma provides the following relatively elementary (with
$\mu_\alpha=0$ and $\nu_\beta=1$ for all $\alpha$ and $\beta$) examples: if
$g_2\in\schw$ is for large $\Vert x\Vert$ equal to one of the majorants in the
first part of Proposition~\ref{prop:examplesqa} and if $g_1$ is of polynomial
growth together with its derivatives, then $g_1g_2\in\schw_{\textup{\text{qa}}}$.

With these remarks in mind we formulate the main theorem for $\schw$. It becomes
clear when one changes all admissible spaces in its formulation to $\schw$---as
is validated by Theorem~\ref{thm:spacedist}---and notes that point evaluations
are continuous. The statement on supports is proved as in Lemma
\ref{lem:testfunctionsannihilators} below.

\begin{theorem}\label{thm:schw}
Let $\gen\subset\schw$. For each $g\in\gen$ choose an admissible subspace $L_g$
and consider the corresponding subspace $L=\sum_{g\in\gen} L_g\cdot g$ of
$\schw$.

Then the annihilator of $L$ consists of all $T\in\tempdist$ such that $g\cdot
T=0$ for all $g\in\gen$. Any such $T$ has support in $\bigcap_{g\in\gen}Z(g)$.

The closure of $L$ is independent of the choice of the $L_g$ and equal to the
closed multiplicative ideal generated by $\gen$ in $\schw$. The subspace $L$ is
dense if and only if the elements of $\gen$ have no common zero.
\end{theorem}

\begin{corollary}\label{cor:schw}
Let $\gen\subset\schw$. Then the module generated by $\gen$ over $\galg$ is
dense in $\schw$ if and only if the elements of $\gen$ have no common zero. If
$\gen\subset\schw_{\textup{\text{qa}}}$ then the analogous statement holds for the module
generated over $\poly$.
\end{corollary}

\begin{remark}\label{rem:schw}
\begin{enumerate}
\item In the cyclic case we retrieve the density of $\poly\cdot \exp(-\Vert x\Vert^2)$
in $\schw$ since the Gaussian has no zeros.
\item In case $L$ is not dense the natural approach in describing
$\overline L$ is via the determination of $L^\perp$. If $\bigcap_{g\in\gen}Z(g)$
is sufficiently regular the structure theorems in \cite{Schwartz}, for
distributions with support in a prescribed set, can then be of assistance.

The easiest case occurs when $\bigcap_{g\in\gen}Z(g)$ has no accumulation point
in $\Rn$. Then $L^\perp$ is as a vector space canonically isomorphic with a
subspace of the direct product $\prod_{x\in\bigcap_{g\in\gen}Z(g)}L_x^\perp$,
where $L_x^\perp$ is the subspace of $L^\perp$ of elements supported in $\{x\}$;
if $\bigcap_{g\in\gen}Z(g)$ is finite, then this isomorphism is onto. This
description is computable in concrete cases and in general the isomorphism shows
what the type of description of $\overline L$ is: a function $\psi\in\schw$ is
in $\overline L$ if and only if $D\psi(x)=0$ for all
$x\in\bigcap_{g\in\gen}Z(g)$ and for all constant coefficient differential
operators $D$ corresponding to a (possibly infinite) base of $L_x^\perp$.

In the
one-dimensional case the answer can be described explicitly, as follows. If
$\bigcap_{g\in\gen}Z(g)$ has no accumulation point in $\R$, then $\overline L$
consists precisely of those $\psi\in\schw(\R)$ having at all
$x\in\bigcap_{g\in\gen}Z(g)$ a zero of an order which is at least the minimum order, as $g$ ranges over $\gen$, of
the zero of $g$ at $x$.
\end{enumerate}
\end{remark}

The density statement in the following corollary parallels the well-known
$L_1$-result:

\begin{corollary}\label{cor:translates}
Let $\gen\subset\schw$. Then the annihilator of $\Span \{\,T_y g\mid
g\in\gen,\,y\in\Rn\,\}$ consists precisely of the tempered distributions $T$
such that $\Fourier(g)\cdot\Fourier^{-1}(\typ)=0$ for all $g\in\gen$. For any
such $\typ$, $\Fourierinverse(\typ)$ has support in
$\bigcap_{g\in\gen}Z(\Fourier(g))$. The span of the translates of the elements
of $\gen$ is dense in $\schw$ if and only if the Fourier transforms of the
elements of $\gen$ have no common zero in $\Rn$.
\end{corollary}

This follows immediately from Corollary~\ref{cor:schw} under Fourier transform.
We note that if $g\in\gen$ is such that $\Fourier(g)\in\schw_{\textup{\text{qa}}}$, then one has
proper subspaces of $\galg$ which are admissible spaces for $\Fourier(g)$ but
correspond to spectral parameters taken from a somewhere dense set only.
Correspondingly, for $g$ such that $\Fourier(g)\in\schw_{\textup{\text{qa}}}$ one can restrict
the translations in the above corollary to those over vectors in a somewhere
dense set. By the same argument, for $\Fourier(g)\in\schw_{\textup{\text{hol}}}$ the
translations can be restricted still further in accordance with Definition
\ref{def:admissiblespaces}. These various smaller sets of translations can be
chosen independently for each of the elements of
$\gen\cap\Fourier^{-1}\left(\schw_{\textup{\text{qa}}}\cup\schw_{\textup{\text{hol}}}\right)$ and both the
description of the annihilator and the density statement in the corollary then
remain valid. As an example, $\{\,T_{y}\exp(-x^2)\mid y\in E\,\}$ is dense in
$\schw(\R)$ for any set $E\subset\R$ having an accumulation point in $\R$.

\begin{corollary}\label{cor:derivatives}
Let $\gen\subset\schw$ be such that $\Fourier(g)\in\schw_{\textup{\text{qa}}}$ for all
$g\in\gen$. Then the annihilator of the span of the derivatives of the elements
of $\gen$ consists precisely of the tempered distributions $T$ such that
$\Fourier(g)\cdot
\Fourierinverse(\typ)=0$ for all $g\in\gen$. For any such $T$,
$\Fourierinverse(\typ)$ has support in $\bigcap_{g\in\gen}Z(\Fourier(g))$. The
span of the derivatives of the elements of $\gen$ is dense in $\schw$ if and
only if the $\Fourier(g)$ have no common zero in $\Rn$.
\end{corollary}

Again the proof is immediate from Corollary~\ref{cor:schw}.

\begin{remark}\label{rem:schwartz}
If $g\in\schw$ is such that $\Fourier(g)\in\schw_{\textup{\text{qa}}}$, then the translates of $g$ over the vectors in a somewhere dense set
span a subspace of $\schw$ with the same closure as the span of the derivatives
of $g$, as is immediate from Theorem~\ref{thm:schw}. One can refine this if
$g\in\schw_{\textup{\text{hol}}}$ and one can formulate similar results for arbitrary sets of
generators. It is also clear how to obtain results on equality of the closure of
the span of the derivatives of one set of generators and the closure of the span
of the translates of a second set of generators. Using the continuity of the
inclusion of $\schw$ in various spaces, these results are subsequently
transferred to such other spaces.

As an example on the real line, suppose that $\gen_{1}\subset\schw(\R)$,
$\gen_{2}\subset\schw(\R)$, and $\Fourier(\gen_{2})\subset\schw(\R)_{\textup{\text{qa}}}$. Assume
that $\bigcap_{g_1\in{\gen_1}}Z(\Fourier(g_1))$ and
$\bigcap_{g_2\in{\gen_2}}Z(\Fourier(g_2))$ are equal and have no accumulation
point in $\R$, so that Remark~\ref{rem:schw} applies to $\{\,\Fourier(g_1)\mid
g_1\in\gen_1\,\}$ and $\{\,\Fourier(g_2)\mid g_2\in\gen_2\,\}$. Then in
$\schw(\R)$ the closure of the span of the translates of the elements of
$\gen_1$ is equal to the closure of the span of the derivatives of the elements
of $\gen_2$ if and only if for all
$x\in\bigcap_{g_1\in{\gen_1}}Z(\Fourier(g_1))$ the minimum order, as $g_1$
ranges over $\gen_1$, of the zero of $\Fourier(g_1)$ at $x$, is equal to the
minimum order, as $g_2$ ranges over $\gen_2$, of the zero of $\Fourier(g_2)$ at
$x$. By continuity the sufficiency of this condition persists in $L_p(\R,dx)$
($1\leq p\leq\infty$). This type of result is similar in spirit to
\cite{MandelbrojtRice}, where refined results for two sets of one generator in
$L_p$-spaces on the real line are derived.

\end{remark}

\subsection{The spaces $\testcmOmega$ and $\testmOmega$ $(0\leq m\leq\infty)$ and their duals}
These spaces of test functions on a nonempty open subset $\Omega$ of $\Rn$ and
their duals are covered by Theorem~\ref{thm:spacedist}. Since these spaces, as
well as the spaces in a defining sequence for the LF-spaces among them, are all
modules over $C^\infty(\Rn)$ we see that all admissible spaces for holomorphic
weights are admissible spaces for arbitrary elements, as in Corollary
\ref{cor:nogrowth}.

We start with a preparatory result.

\begin{lemma}\label{lem:testfunctionsannihilators}
Suppose $0\leq m\leq \infty$ and let $V$ denote $\testcmOmega$ or $\testmOmega$.
Then, for $T\in V^\prime$ and $g\in V$, the following are equivalent:
\begin{enumerate}
\item $\langle \typ,\test\cdot g\rangle=0$.
\item $\langle \typ,\testcOmega\cdot g\rangle=0$.
\item $\langle \typ,V\cdot g\rangle=0$.
\end{enumerate}
If $\typ$ satisfies these conditions, then $\supp \typ\subset Z(g)$.
\end{lemma}

\begin{proof}
We give the proof for $\testcmOmega$; the case $\testmOmega$ is similar. It is
trivial that (1) implies (2). Assume (2), i.e., assume that $g\cdot
T\in\distmOmega$ vanishes on $\testcOmega$. Then, by density, $g\cdot T$ is
identically zero on $\testcmOmega$. Thus (2) implies (3). Next assume (3) and
note that $g\cdot T\in\compdistmOmega$. By density of $\testcmOmega$ in
$\testmOmega$ we conclude from (3) that $g\cdot T=0$ in $\compdistmOmega$. In
particular, the product vanishes on the restriction of elements of $\test$ to
$\Omega$. Hence (3) implies (1). Finally, if $\psi\in\testcOmega$ and
$\supp\psi\cap Z(g)=\emptyset$, then $\psi$ factors as $\psi=g\cdot\chi$ for
some $\chi\in\testcmOmega$. Consequently, $\langle T,\psi\rangle=\langle T,g\cdot
\chi\rangle=0$.
\end{proof}

For $m=0$ the condition $\supp \typ\subset Z(g)$ is equivalent to the other
statements (see Lemma~\ref{lem:testfunctionsannihilatorsX}).

The equivalence of (1) and (3) in the above lemma yields the following:

\begin{corollary}\label{cor:testfunctionsideals}
Let $V$ denote $\testcmOmega$ or $\testmOmega$ and suppose $\gen\subset V$ is a generating set. Then
$\overline{\sum_{g\in\gen}\test\cdot g}$ is equal to the closed ideal in $V$
generated by $\gen$.
\end{corollary}

The following theorem becomes clear when one changes all admissible spaces in
its formulation to $\test$---as is validated by Theorem
\ref{thm:spacedist}---and subsequently uses Lemma
\ref{lem:testfunctionsannihilators}, Corollary~\ref{cor:testfunctionsideals}, and
the continuity of point evaluations.

\begin{theorem}\label{thm:testfunctions}
Let $0\leq m\leq\infty$. Endow $\testcmOmega$ with its usual LF-topology \textup{(}resp.,
$\testmOmega$ with its usual Fr\'echet topology\textup{)}. Suppose
$\gen\subset\testcmOmega$ \textup{(}resp., $\gen\subset\testmOmega$\textup{)}. For each $g\in\gen$
choose an admissible space $L_g$ for a holomorphic weight and consider the
corresponding subspace $L=\sum_{g\in\gen} L_g\cdot g$.

Then the annihilator of $L$ consists of all $T\in\distmOmega$ \textup{(}resp., all
$T\in\compdistmOmega$\textup{)} such that $g\cdot T=0$ for all $g\in\gen$. Any such $T$
has support in $\bigcap_{g\in\gen}Z(g)$.

The closure of $L$ in $\testcmOmega$ \textup{(}resp., $\testmOmega$\textup{)} is equal to the
closed ideal generated by $\gen$ in $\testcmOmega$ \textup{(}resp., $\testmOmega$\textup{)}. The subspace $L$ is dense in $\testcmOmega$ \textup{(}resp., $\testmOmega$\textup{)} if and only if the
elements of $\gen$ have no common zero.

In the case of $\testcmOmega$ the sequential closure of $L$ is independent of
the choice of the $L_g$ if $\gen$ is finite.
\end{theorem}

\begin{corollary}
Suppose $\gen\subset\testcmOmega$ \textup{(}resp., $\gen\subset\testmOmega$\textup{)}. Then the module generated by $\gen$ over $\poly$ is dense in
$\testcmOmega$ \textup{(}resp., $\testmOmega$\textup{)} if and only if the elements of $\gen$ have no common zero.
The analogous statement holds for the module generated over $\galg$.
\end{corollary}

\begin{remark}
\begin{enumerate}
\item Specializing the corollary to $\gen=\{1\}$ we retrieve that $\poly$ is dense in $\testmOmega$ because the constant function $1$
has no zeros (cf.\ \cite[p.~160]{Treves}).
\item Analogously to Remark~\ref{rem:schw} we note that if $\bigcap_{g\in\gen}Z(g)$ has no accumulation point in $\Omega$, then
$L^\perp$ is as linear space canonically isomorphic to the direct sum
$\bigoplus_{x\in\bigcap_{g\in\gen}Z(g)}L_x^\perp$ in the case of $\testmOmega$
(resp., the direct product $\prod_{x\in\bigcap_{g\in\gen}Z(g)}L_x^\perp$ in the
case of $\testcmOmega$), where $L_x^\perp$ is the subspace of $L^\perp$ of
elements supported in $\{x\}$. Again this makes it evident that in both the
cases $\testcmOmega$ and $\testmOmega$ an element $\psi$ of the respective space
is in $\overline L$ if and only if $D\psi(x)=0$ for all
$x\in\bigcap_{g\in\gen}Z(g)$ and for all constant coefficient differential
operators $D$ corresponding to a (possibly infinite) base of $L_x^\perp$. In the
one-dimensional situation one finds again that, if $\bigcap_{g\in\gen}Z(g)$ has
no accumulation point in $\Omega\subset\R$, then $\overline L$ consists precisely
of those $\psi\in\testcmOmega$ (resp., $\psi\in\testmOmega$) having in all
$x\in\bigcap_{g\in\gen}Z(g)$ a zero of an order which is at least the minimum
order, as $g$ ranges over $\gen$, of the zero of $g$ at $x$.
\item For finite $m$, the fact that the ideal $\testmOmega\cdot\gen$ is dense in $\testmOmega$ if the
generators have no common zero can also be seen to be a consequence of the results on
closed ideals in $\testmOmega$ in \cite{Whitney}.

In \cite{Nachbin3} a necessary and sufficient condition for the density of
subalgebras of $C^m(M)$ ($m\geq 1$) is given, where $C^m(M)$ is the algebra of
functions of class $C^m$ on a manifold $M$ of class $C^m$, supplied with the
topology of uniform convergence of order $m$ on compact subsets. It also follows
from this condition that the ideal $\testmOmega\cdot\gen$ is dense in
$\testmOmega$ for finite $m \geq 1$.
\item For $m=0$ more precise results are available in Theorem~\ref{thm:testfunctionsX}.
\end{enumerate}
\end{remark}

The dual theorem is as follows:

\begin{theorem}\label{thm:distributions}
Let $0\leq m\leq\infty$. Endow $\distmOmega$ \textup{(}resp., $\compdistmOmega$\textup{)} with the
weak topology. Suppose $\gen\subset\distmOmega$ \textup{(}resp.,
$\gen\subset\compdistmOmega$\textup{)}. For each $T\in\gen$ choose an admissible space
$L_T$ for a holomorphic weight and consider the corresponding subspace
$\sum_{T\in\gen} L_T\cdot T$.

Then the annihilator of $L$ consists of all $\psi\in\testcmOmega$ \textup{(}resp., all
$\psi\in\testmOmega$\textup{)} such that $\psi\cdot T=0$ for all $T\in\gen$. Any such $\psi$ vanishes on $\overline{\bigcup_{T\in\gen}\supp T}$, where closure is
taken in $\Omega$.

The closure of $L$ in $\distmOmega$ \textup{(}resp., $\compdistmOmega$\textup{)} is independent of
the choice of the $L_T$. Each element of $\overline L$ has support in
$\overline{\bigcup_{T\in\gen}\supp T}$. The subspace $L$ is dense in $\distmOmega$ \textup{(}resp.,
$\compdistmOmega$\textup{)} if and only if $\overline{\bigcup_{T\in\gen}\supp T}=\Omega$.

The sequential closure of $L$ is in both cases independent of the choice of the
$L_T$ if $\gen$ is finite.

In the case of $\distmOmega$, if the set $\gen$ of
arbitrary cardinality consists of continuous functions on $\Omega$, then the
sequential closure of $L$ is independent of the choice of the $L_T$. $L$ is in this case sequentially dense if the elements of
$\gen$ have no common zero.
\end{theorem}

\begin{proof}
When changing all admissible spaces to $\test$---as is validated by Theorem
\ref{thm:spacedist}---all statements become clear, except the additional results
on sequential closure in $\distmOmega$ in the last paragraph. As to these, note
that by Theorem~\ref{thm:testfunctions} the closure of $L$, when considered as
subspace of $\testzeroOmega$ in its Fr\'echet topology, is independent of the
choice of the $L_T$. The continuity of the canonical injection
$\testzeroOmega\hookrightarrow\distmOmega$ then implies by Theorem
\ref{thm:closurecont} that the sequential closure in $\distmOmega$ is also
independent of the choice of the $L_T$. If in addition the functions in $\gen$
have no common zero, then by Theorem~\ref{thm:testfunctions} $L$ is dense in
$\testzeroOmega$. Stated otherwise: $L$ and $\testcOmega$ have equal closure in
$\testzeroOmega$. Again, by Theorem~\ref{thm:closurecont}, when $L$ is considered
as a subspace of $\distmOmega$, it still has the same sequential closure as
$\testcOmega$. Since the latter space is well-known to be sequentially dense in
$\distmOmega$, so is $L$.
\end{proof}

\begin{corollary}
\begin{enumerate}
\item Suppose $\gen\subset\distmOmega$ \textup{(}resp.\ $\gen\subset\compdistmOmega$\textup{)}.
Then the modules generated by $\gen$ over $\poly$ and $\galg$ are both weakly dense in
$\distmOmega$ \textup{(}resp.\ $\compdistmOmega$\textup{)} if and only if
$\overline{\bigcup_{T\in\gen}\supp T}=\Omega$, where closure is taken in
$\Omega$.
\item Suppose $\gen\subset C(\Omega)$ and assume that the elements of $\gen$ have no
common zero. Then the modules generated by $\gen$ over $\poly$ and $\galg$ are
both sequentially dense in $\distmOmega$.
\end{enumerate}
\end{corollary}

\begin{remark}
\begin{enumerate}
\item Specializing to $\gen=\{1\}$ we retrieve that $\poly$ is
sequentially dense in $\distmOmega$ because the constant continuous function $1$
has no zeros (cf.\ \cite[p.~304]{Treves}).
\item For $m=0$ more precise results are available in Theorem~\ref{thm:Radonmeasures}.
\end{enumerate}
\end{remark}

\subsection{General spaces of Bernstein type}\label{subsec:Bernsteinspaces}
We will now consider a generalization of the classical Bernstein problem, both
in terms of differentiability and dimensionality. The relevant spaces, of which
our spaces $\opw$ are particular examples, were introduced and studied in
\cite{Zapata,Nachbin2}. We will generalize the results in [loc.\ cit.] and explain
the connection with the well-known sufficient condition for the polynomials to
be dense in the classical continuous one-dimensional case.

First we introduce the analogue of our spaces $\opw$ in arbitrary degree
of differentiability, as follows. Suppose $m\in\{0,1,2,\ldots\,,\infty\}$ and
$w\in\weights$. Let $C_w^m$ consist of all $f\in C^m(\Rn)$ such that
$\lim_{\Vert x\Vert\rightarrow\infty}\parta f(x)\cdot w(x)=0$ for all
$\alpha\in\Nn$, $\vert\alpha\vert\leq m$. Then $C_w^m$ becomes a locally convex
space when supplied with the seminorms $\pa\,(\vert\alpha\vert\leq m)$ defined
by $\pa(f)=\Vert \parta f\cdot w\Vert_\infty\,(f\in C_w^m)$. The spaces $\opw$
correspond to $m=\infty$.

\begin{lemma}\label{lem:Bernsteindense}
Let $0\leq m\leq\infty$ and $w\in\weights$. If $\gen\subset C_w^m$ is such that
$\bigcap_{g\in\gen}Z(g)=\emptyset$, then $\sum_{g\in\gen}\test\cdot g$ is dense
in $C_w^m$.
\end{lemma}

\begin{proof}
The injection $j: C_c^m(\Rn)\mapsto C_w^m$ is continuous and the image is dense
\cite{Zapata}. Let $j^t:(C_w^m)^\prime\mapsto C_c^m(\Rn)^\prime$ be the
injective transposed map. If $T\in (C_w^m)^\prime$ vanishes on the subspace in
the lemma, then for all $g\in\gen$ we see in particular that $g\cdot j^t(T)\in
C_c^m(\Rn)^\prime$ vanishes on $\testcOmega$. Hence $g\cdot j^t(T)=0$ and therefore $\supp
j^t(T)\subset Z(g)$ for all $g$, implying that $j^t(T)=0$. Hence $T=0$.
\end{proof}

The routine verification of the following lemma is left to the reader:

\begin{lemma}\label{lem:Bernsteincontinuity}
Let $0\leq m\leq\infty$ and $w\in\weights$ and suppose $0\leq \nu\leq 1$. Then
$C_{w^\nu}^m\subset C_w^m$. If $g\in C_{w^\nu}^m$, then the assignment $f\mapsto
f\cdot g$ maps $\op_{w^{1-\nu}}$ continuously into $C_w^m$.
\end{lemma}

\begin{theorem}\label{thm:Bernsteinspaces}
Let $0\leq m\leq\infty$ and $w\in\weights$. Then:
\begin{enumerate}
\item All admissible spaces for $w$ are dense in $C_w^m$.
\item Suppose $\gen\subset C_w^m$ is such that $\bigcap_{g\in\gen}Z(g)=\emptyset$. For
each $g\in\gen$ choose $0\leq\nu_g\leq 1$ such that $g\in C_{w^{\nu_g}}^m$. Let
$L_g$ be an admissible space for the weight $w^{1-\nu_g}$.

Then $\sum_{g\in\gen} L_g\cdot g$ is dense in $C_w^m$.
\end{enumerate}
\end{theorem}

\begin{proof}
For the first part, note that the canonical injection $C_c^m(\Rn)\hookrightarrow
C_w^m$ is continuous. The image is dense \cite{Zapata} and we conclude that
$\test$ is dense in $C_w^m$. The canonical injection $\opw\hookrightarrow C_w^m$
is also continuous, so we subsequently infer from Lemma
\ref{lem:fundamentallemma} and Theorem~\ref{thm:admissiblespaces} that all
admissible spaces for $w$ are dense in $C_w^m$.

As to the second part, for each $g\in\gen$, Lemma~\ref{lem:Bernsteincontinuity}
provides a continuous map from $\op_{w^{1-\nu_g}}$ into $C_w^m$ defined by
sending $f$ to $f\mapsto f\cdot g$. The theorem is therefore an immediate
consequence of Theorem~\ref{thm:generalthm} and Lemma~\ref{lem:Bernsteindense}.
\end{proof}

Taking all $\nu_g$ in the second part equal to, e.g., one-half, we have the following
consequence. The growth conditions below can be relaxed for any explicitly given
weight, in particular, in the second part where polynomial bounds can be replaced
by the reciprocals of appropriate quasi-analytic weights.

\begin{corollary}\label{cor:Bernsteinspaces}
Let $0\leq m\leq\infty$. Suppose $\gen\subset C_w^m$ is such that
$\bigcap_{g\in\gen}Z(g)=\emptyset$.
\begin{enumerate}
\item If $w\in\weights_0$ and if each of the derivatives of order at most $m$
of an arbitrary element of $\gen$ is bounded, then the module generated by
$\gen$ over $\galg$ is dense in $C_w^m$.
\item If $w\in\weightsqa$ and if each of the derivatives of order at most $m$
of an arbitrary element of $\gen$ is of at most polynomial growth, then the
module generated by $\gen$ over $\poly$ is dense in $C_w^m$.
\end{enumerate}
\end{corollary}

In order to establish the connection with the continuous one-dimensional
Bernstein problem, consider a continuous even strictly positive weight $w:\R\mapsto\R_{\geq 0}$ with the property
that $s\mapsto
-\log w(e^s)$ is convex on $(-\infty,\infty)$ and such that
\begin{equation*}
\int_0^\infty\frac{\log w(t)}{1+t^2}=-\infty.
\end{equation*}
As a consequence of the second part of Theorem~\ref{thm:classonedim} $w$ is a quasi-analytic weight
and we conclude from Corollary~\ref{cor:Bernsteinspaces} for $\gen=\{1\}$ that
the polynomials are dense in $C_w^m$ for all $m$; moreover, Theorem
\ref{thm:Bernsteinspaces} asserts that the same holds for the span of the
exponentials with spectral parameters in a somewhere dense subset of $\Rn$.
Specializing further by letting $m=0$ one retrieves a well-known sufficient condition for the
polynomials to be dense in the continuous one-dimensional Bernstein space
\cite{Mergelyan,Carleson,Koosis}. It becomes apparent that this classical result is an
aspect of a more general picture which is described by Theorem
\ref{thm:Bernsteinspaces}.

Let us now explain the connection with the results in \cite{Zapata,Nachbin2}. In
[loc.\ cit.] the most general version of the spaces $C_w^m$ was introduced, as
follows. If $m\in\{0,1,2,\ldots,\infty\}$, let ${\bf w}=\{\,w_\alpha\mid
\alpha\in\Nn,\,\vert\alpha\vert\leq m\,\}$ be a set of weights on $\Rn$. These weights are
assumed to be upper semicontinuous in [loc.\ cit.], but this can be relaxed: we
will continue to assume only that these weights are bounded. Consider the space
$C_{\bf w}^m$ of all $f\in C^m(\Rn)$ such that $\lim_{\Vert
x\Vert\rightarrow\infty}\parta f(x)\cdot w_\alpha(x)=0$ for all
$\alpha\in\Nn,\,\vert\alpha\vert\leq m$. The obvious choice of seminorms makes
$C_{\bf w}^m$ into a locally convex space.

In [loc.\ cit.] the emphasis is on the situation where ${\bf w}$ is a decreasing
set of weights, i.e., where for all $\alpha,\,\beta\in\Nn$
($\vert\alpha\vert,\,\vert\beta\vert\leq m$) with $\alpha_j\geq\beta_j$ ($1\leq
j\leq n$) there exists a constant $C_{\alpha,\beta}$ such that
$w_{\alpha}(x)\leq C_{\alpha,\beta} w_{\beta}(x)$ for all $x\in\Rn$. It is
easily verified that in that case there is a continuous injection
$C_{w_0}^m\hookrightarrow C_{\bf w}^m$, where $w_0\in{\bf w}$ is the weight
corresponding to $\alpha=0$. Now $C_c^m(\Rn)$ is dense in both $C_{w_0}^m$ and
$C_{\bf w}^m$ (see \cite{Zapata}), so Lemma~\ref{lem:fundamentallemma} implies
that any dense subspace of $C_{w_0}^m$ is also dense in $C_{\bf w}^m$. In this
way the density results in the above theorem for $C_{w_0}^m$ are transferred to
$C_{\bf w}^m$, provided that ${\bf w}$ is a decreasing set of weights. One thus
obtains various (not necessarily cyclic) density statements in $C_{\bf w}^m$. In
the cyclic case one concludes, e.g., that $\poly\cdot p$ is dense in $C_{\bf
w}^m$ for all $p\in\poly$ without zeros on $\Rn$, provided that $\bf w$ a decreasing set of weights such that $w_0\in\weightsqa$.
Specializing still further by letting $p=1$ yields the main result in
[loc.\ cit.], namely that $\poly$ is dense in $C_{\bf w}^m$ if ${\bf w}$ is a
decreasing set of weights such that $w_0\in\weightsqa$.

\subsection{Spaces associated with a locally compact set $X$}\label{subsec:locallycompact}
Throughout this section $X$ denotes a locally compact subset of $\Rn$. As a
consequence of \cite[I.3.3 and I.9.7]{Bourbaki1} such $X$ can alternatively be
characterized as the intersection of an open subset $\Omega$ and a closed subset
$C$ of $\Rn$. This structure equation is convenient in the proof of several
statements below. We choose and fix such $\Omega$ and $C$ and assume
$X=\Omega\cap C\neq\emptyset$.
In our terminology a Borel measure on $X$ takes finite values on compact subsets of $X$.

Let $C_0(X)$ denote the continuous functions on $X$ vanishing at infinity (with respect to $X$),
supplied with the maximum norm. $C(X)$ denotes the continuous functions on $X$ with its subspace $C_c(X)$ of functions with compact support. These two spaces carry a standard topology which we now describe.

Choose a sequence $\{K_m\}_{m=1}^\infty$ of compact subsets of $X$ such that
$X=\bigcup_{m=1}^\infty K_m$ and such that $K_m$ is contained in the interior
(in the topology of $X$) of $K_{m+1}$. Define for $m=1,2,\ldots$ the seminorms
$p_m$ on $C(X)$ by $p_m(f)=\max_{x\in K_m}\vert f(x)\vert$. Then the $p_m$
determine on $C(X)$ the Fr\'echet topology of uniform convergence on compact
subsets of $X$. For $m=1,2\dots$, let $C_{K_m}(X)$ be the subspace of $C(X)$
consisting of the continuous functions on $X$ with support in $K_m$, supplied
with the induced topology of $C(X)$. Then the $C_{K_m}$ form an increasing
sequence of Fr\'echet spaces, each inducing for $m\geq 2$ on its predecessor its given
topology. This sequence therefore provides $C_c(X)$ with an LF-topology, which
is independent of the choice of the $K_m$. The canonical injection
$C_c(X)\hookrightarrow C(X)$ is continuous with dense image.

The dual space $C_c(X)^\prime$ consists of the Radon measures on $X$. Such a
Radon measure is locally expressible as a complex-valued measure. The dual
$C(X)^\prime\subset C_c(X)^\prime$ of $C(X)$ is identified with complex Borel
measures on $X$ with compact support. We supply $C_c(X)^\prime$ and
$C(X)^\prime$ each with its weak topology.

We will be concerned with the spaces $C(X)$, $C_c(X)$, $C_0(X)$, $C(X)^\prime$, and $C_c(X)^\prime$, as well as with $L_p(X,
\mu)$ ($1\leq p<\infty$) for an arbitrary Borel measure $\mu$ on $X$ (and its  completion). For these $L_p$-spaces it is important to note that any relatively open
subset of $X$ can be written as a countable union of compact subsets, since this holds for $\Omega$. This excludes certain measure-theoretic pathologies since it not
only implies that $\mu$ is $\sigma$-finite but also that $\mu$ is automatically
regular \cite[Theorem 2.18]{RudinRCA}.

Although the continuity hypothesis in Theorem~\ref{thm:generalthm} is easily
verified in a direct fashion in the situations below, it is instructive to note
that our spaces of interest are in fact covered by Theorems \ref{thm:spacedist},
\ref{thm:spacefunctions}, and Corollary~\ref{cor:nogrowth}. Indeed,
Corollary~\ref{cor:nogrowth} applies to the $C^\infty(\Rn)$-modules $C(X)$,
$C_c(X)$, $C(X)^\prime$, and $C_c(X)^\prime$ and Theorem~\ref{thm:spacefunctions}
handles $C_0(X)$. For the remaining case of $L_p$-spaces we note that for any
Borel measure $\mu$ on $X$ the space $L_1^{\text{\textup{loc}}}(X,\mu)$ of locally integrable
functions is canonically identified with a subspace of $C_c(X)^\prime$ via the
pairing
\begin{equation*}
\langle f,\psi\rangle=\int_X f(x)\psi(x)\,d\mu\quad(f\in L_1^{\text{\textup{loc}}}(X,\mu),\,\psi\in
C_c(X)).
\end{equation*}
This applies in particular to $L_p(X,d\mu)\subset L_1^{\text{\textup{loc}}}(X,\mu)$ for $1\leq
p\leq\infty$, so that $L_p(X,\mu)\hookrightarrow C_c(X)^\prime$. One verifies
easily that the induced weak topology of $C_c(X)^\prime$ is weaker than the
original topology of $L_p(X,\mu)$, so that the second possibility in Theorem
\ref{thm:spacedist} applies to $L_p(X,\mu)$ for $1\leq p\leq\infty$.

This being said, we remark that our method obviously requires supplementary information on the interplay between
$\test$ and functions on $X$. This is provided by the following three results,
which are based on the structure equation $X=\Omega\cap C$.

\begin{lemma}\label{lem:restrictionsurjective}
Taking restrictions from $\Omega$ to $X$ maps $C_c(\Omega)$ continuously into
and onto $C_c(X)$.
\end{lemma}

\begin{proof}
Let $\widetilde{\psi}\in C_c(\Omega)$. If $x\in X$ and $\widetilde\psi(x)\neq 0$,
then $x\in X\cap\supp\widetilde\psi$, which is closed in $X$. Thus,
$\supp\widetilde\psi\vert_X\subset X\cap\supp\psi=C\cap\supp\widetilde\psi$.
Since the latter set is compact the restriction map is into. It is evidently
continuous. To show surjectivity, let $\psi\in C_c(X)$ and note that
$X=\Omega\cap C$ is closed in the metrizable (hence normal) space $\Omega$. By
the Tietze extension theorem \cite[IX.4.2]{Bourbaki2} there exists
$\widetilde\psi\in C(\Omega)$ extending $\psi$. By a version of Urysohn's lemma
 \cite[2.12]{RudinRCA}) there exists $\widetilde f\in C_c(\Omega$) which is identically $1$ on $\supp\psi$. Then $\widetilde
f\widetilde\psi\in C_c(\Omega)$ extends $\psi$.
\end{proof}

\begin{lemma}\label{lem:testfunctionsannihilatorsX}
Let $V$ denote $C_c(X)$, $C(X)$, or $C_0(X)$. Then, for $T\in V^\prime$ and $g\in
V$, the following are equivalent:
\begin{enumerate}
\item $\langle T,\test\cdot g\rangle=0$.
\item $\langle T, C_c(X)\cdot g\rangle=0$.
\item $\langle T,V\cdot g\rangle=0$.
\item $\supp T\subset Z(g)$.
\end{enumerate}
\end{lemma}

\begin{proof}
We give the proof for $V=C(X)$; the other cases are similar. Assume (1) and
consider the continuous linear functional $\psi\mapsto\langle T,\psi\cdot
g\rangle$ on $C(\Omega)$. The assumption implies that in particular $\langle
T,\testcOmega\cdot g\rangle=0$; therefore $\langle T,C(\Omega)\cdot g\rangle=0$
by the density of $\testcOmega$ in $C(\Omega)$. In particular $\langle
T,C_c(\Omega)\cdot g\rangle=0$, which by the previous lemma is equivalent to
$\langle T,C_c(X)\cdot g\rangle=0$. Hence (1) implies (2). Assuming (2),
consider the continuous linear functional $\psi\mapsto\langle T,\psi\cdot
g\rangle$ on $C(X)$. By density of $C_c(X)$ in $C(X)$ we conclude that (2)
implies (3). It is trivial that (3) implies (1). Analogously to the proof of
Lemma~\ref{lem:testfunctionsannihilators} one sees that (2) implies (4). There
exists a complex Borel measure $\mu$ on X with compact support such that
\begin{equation}\label{eq:representation}
\langle T,f\rangle=\int_X f(x)\,d\mu(x)\quad(f\in V),
\end{equation}
which makes it clear that (4) implies, e.g., (2).
\end{proof}

\begin{corollary}\label{cor:testfunctionsXideals}
Let $V$ denote $C_c(X)$, $C(X)$, or $C_0(X)$ and suppose $\gen\subset V$ is a generating set. Then
$\overline{\sum_{g\in\gen}\test\cdot g}$ is equal to the closed ideal generated
in $V$ by $\gen$.
\end{corollary}

The following result now becomes clear from the combination of
Corollary~\ref{cor:nogrowth}, Lemma~\ref{lem:testfunctionsannihilatorsX},
\eqref{eq:representation} (and its local analogue for $C_c(X)^\prime$) and the previous corollary.

\begin{theorem}\label{thm:testfunctionsX}
Endow $C(X)$ with its Fr\'echet topology \textup{(}resp., $C_c(X)$ with its LF-topology\textup{)}
and suppose $\gen\subset C(X)$ \textup{(}resp., $\gen\subset C_c(X)$\textup{)}. For each
$g\in\gen$ choose an admissible space $L_g$ for a holomorphic weight and
consider the corresponding subspace $L=\sum_{g\in\gen} L_g\cdot g$.

Then the annihilator of $L$ consists of all complex Borel measures $\mu$ with
compact support \textup{(}resp., all Radon measures $T$\textup{)} such that the support of $\mu$
\textup{(}resp., $T$\textup{)} is contained in $\bigcap_{g\in\gen}Z(g)$.

The closure of $L$ in $C(X)$ \textup{(}resp., $C_c(X)$\textup{)} is equal to the closed ideal
generated by $\gen$ in $C(X)$ \textup{(}resp., $C_c(X)$\textup{)} and consists of all functions in
$C(X)$ \textup{(}resp.\ $C_c(X)$\textup{)} vanishing on $\bigcap_{g\in\gen}Z(g)$. In particular,
$L$ is dense in $C(X)$ \textup{(}resp., $C_c(X)$\textup{)} if and only if the elements of $\gen$
have no common zero.

In the case of $C_c(X)$ the sequential closure of $L$ is independent of the
choice of the $L_g$ if $\gen$ is finite.
\end{theorem}

\begin{corollary}
Suppose $\gen\subset C(X)$ \textup{(}resp., $\gen\subset C_c(X)$\textup{)}. Then the module generated by $\gen$ over $\poly$ in $C(X)$ \textup{(}resp., $C_c(X)$\textup{)} is dense
in $C(X)$ \textup{(}resp., $C_c(X)$\textup{)} if and only if the elements of $\gen$ have no common
zero. The analogous statement holds for the module generated over $\galg$.
\end{corollary}

\begin{remark}
\begin{enumerate}
\item In particular $\poly$ is dense in $C(X)$ since the constant function $1$ has no zeros.
For compact $X$ we thus retrieve Weierstrass' result.
\item For compact $X$ the explicit description of the closed ideal $\overline L$
of the Banach algebra $C(X)$ also follows from \cite[Corollary 8.3.1]{Larsen}.
\end{enumerate}
\end{remark}

The dual theorem follows from Corollary~\ref{cor:nogrowth}, Lemma
\ref{lem:testfunctionsannihilatorsX}, and \eqref{eq:representation} (and its
local analogue for $C_c(X)^\prime$). The last statement in the theorem is
demonstrated by considering $C(X)$, analogously to the proof of Theorem
\ref{thm:distributions}.

\begin{theorem}\label{thm:Radonmeasures}
Endow $C(X)^\prime$ \textup{(}resp., $C_c(X)^\prime$\textup{)} with the weak topology. Suppose
$\gen\subset C(X)^\prime$ \textup{(}resp., $\gen\subset C_c(X)^\prime$\textup{)}. For each
$T\in\gen$ choose an admissible space $L_T$ for a holomorphic weight and
consider the corresponding subspace $L=\sum_{T\in\gen} L_T\cdot T$ of compactly
supported complex Borel measures \textup{(}resp., Radon measures\textup{)}.

Then the annihilator of $L$ consists of all $\psi\in C(X)$ \textup{(}resp., all $\psi\in
C_c(X)$\textup{)} such that $\psi$ vanishes on $\overline{\bigcup_{T\in\gen}\supp T}$,
where closure is taken in $X$.

The closure of $L$ in $C(X)^\prime$ \textup{(}resp., $C_c(X)^\prime$\textup{)} consists of
all
compactly supported complex Borel measures \textup{(}resp., consists of all Radon measures\textup{)} $S$ in
$C(X)^\prime$ \textup{(}resp., $C_c(X)^\prime$\textup{)} such that $\supp
S\subset\overline{\bigcup_{T\in\gen}\supp T}$. The subspace $L$ is dense in
$C(X)^\prime$ \textup{(}resp., $C_c(X)^\prime$\textup{)} if and only if
$\overline{\bigcup_{T\in\gen}\supp T}=X$.

The sequential closure of $L$ is in both cases independent of the choice of the
$L_T$ if $\gen$ is finite. In the case of $C_c(X)^\prime$, if the set $\gen$ of
arbitrary cardinality consists of continuous functions on $X$, then the
sequential closure of $L$ is independent of the choice of the $L_T$.
\end{theorem}

\begin{corollary}
Suppose $\gen\subset C(X)^\prime$ \textup{(}resp., $\gen\subset C_c(X)^\prime$\textup{)}. Then the
module generated by $\gen$ over $\poly$ is weakly dense in $C(X)^\prime$ \textup{(}resp.,
$C_c(X)^\prime$\textup{)} if and only if $\overline{\bigcup_{T\in\gen}\supp T}=X$, where
closure is taken in $X$. The analogous statement holds for the module generated
over $\galg$.
\end{corollary}

We now consider $C_0(X)$ and discuss the partial filtration of Subsection
\ref{subsec:filtration}. Evidently, $C_0(X)\subset C_0(X)_\weights$ by
considering the weight $w=1$. In fact, $C_0(X)=C_0(X)_0$. To see this, let $g\in
C_0(X)$. For $x\in\Rn$ let $w(x)=\vert g\vert^{1/2}$ if $x\in X$ and $w(x)=0$ if
$x\notin X$. Then $w\in \weights_0$ and $\opw\cdot g\subset C_0(X)$. Furthermore
if $g\in C_0(X)$ and if there exists a weight $w$ in $\weights_d$ for some
$0<d\leq\infty$ (resp., a weight $w$ in $\weightsqa$, resp., a weight $w$ in
$\weightshol$) which is strictly positive on $X$ and such that $w^{-1}\cdot g$
vanishes at infinity (with respect to $X$), then $g\in C_0(X)_d$ (resp., $g\in C_0(X)_{\textup{\text{qa}}}$, resp., $g\in C_0(X)_{\textup{\text{hol}}}$) since $\opw\cdot f\subset C_0(X)$.

For $X=\Rn$ we note the following explicit example: if $g$ is continuous on $\Rn$ and if $\vert
g\vert$ is equal to one of the majorants in the first part of Proposition
\ref{prop:examplesqa} for sufficiently large $\Vert x\Vert$, then $g\in
C_0(\Rn)_{\textup{\text{qa}}}$. This follows from considering the quasi-analytic weight $w$
given by $w=\vert g\vert^{1/2}$ for large $\Vert x\Vert$ and defined to be equal
to $1$ at the remaining part of $\Rn$.

The following theorem follows from Theorem~\ref{thm:spacefunctions}, Lemma
\ref{lem:testfunctionsannihilatorsX}, Corollary~\ref{cor:testfunctionsXideals},
and the obvious analogue of \eqref{eq:representation}. Here, as with $C(X)$, it
is also possible to resort to Banach algebra theory for the explicit description
of the closure once this closure has been identified as an ideal.

\begin{theorem}\label{thm:vanishinfinity}
Let $\gen\subset C_0(X)$. For each $g\in\gen$ choose an admissible subspace
$L_g$ and consider the corresponding subspace $L=\sum_{g\in\gen} L_g\cdot g$.

Then the annihilator of $L$ consists of all complex Borel measures $\mu$ such
that the support of $\mu$ \textup{(}resp., $T$\textup{)} is contained in $\bigcap_{g\in\gen}Z(g)$.

The closure of $L$ is equal to the closed ideal generated by $\gen$ in $C_0(X)$
and consists of all functions in $C_0(X)$ vanishing on $\bigcap_{g\in\gen}Z(g)$.
The subspace $L$ is dense in $C_0(X)$ if and only if the elements of $\gen$
have no common zero.
\end{theorem}

\begin{corollary}
Let $\gen\subset C_0(X)$. Then the module generated by $\gen$ over $\galg$ is
dense in $C_0(X)$ if and only if the elements of $\gen$ have no common zero. If
$\gen\subset C_0(X)_{\textup{\text{qa}}}$, then the analogous statement holds for the module
generated over $\poly$.
\end{corollary}

We turn to the spaces $L_p(X,\mu)$ ($1\leq p\leq \infty$) for a Borel measure
$\mu$ on $X$. As with $C_0(X)$ we have $L_p(X,\mu)\subset L_p(X,\mu)_\weights$.
If $g\in L_p(X,\mu)$ and if there exists a weight $w$ in $\weights_d$ for some
$0\leq d\leq\infty$ (resp., a weight $w$ in $\weightsqa$, resp., a weight $w$ in
$\weightshol$) such that $w^{-1}\cdot g\in L_p(X,\mu)$, then $g\in L_p(X,\mu)_d$
(resp., $g\in L_p(X,\mu)_{\textup{\text{qa}}}$, resp., $g\in L_p(X,\mu)_{\textup{\text{hol}}}$).

For $X=\Rn$ we note the following explicit example: if $g$ is Borel measurable, if $\vert
g\vert$ is equal to one of the majorants in the first part of Proposition
\ref{prop:examplesqa} for sufficiently large $\Vert x\Vert$, and if $\vert
g\vert^\nu\in L_p(\Rn,\mu)$ for some $0\leq\nu<1$, then $g\in
L_p(\Rn,\mu)_{\textup{\text{qa}}}$. This follows from considering the quasi-analytic weight $w$
given by $w(x)=\vert g\vert^{1-\nu}$ for sufficiently large $\Vert x\Vert$ and
defined to be equal to $1$ at the remaining part of $\Rn$.

\begin{theorem}\label{thm:Lpspaces}
Let $\mu$ be a Borel measure on $X$ and suppose $\gen\subset L_p(X,\mu)$ \textup{(}$1\leq
p<\infty$\textup{)} is countable. For each $g\in\gen$ choose an admissible subspace $L_g$
and consider the corresponding subspace $L=\sum_{g\in\gen} L_g\cdot g$.

Then the annihilator of $L$ consists of all $f\in L_q(X,\mu)$ vanishing
$\mu$-almost everywhere on the complement of the common zero locus of the
$g\in\gen$.
Here $q$ is the conjugate exponent of $p$.

The closure of $L$ consists of all $g\in L_p(X,\mu)$ vanishing $\mu$-almost
everywhere on the common zero locus of the $g\in\gen$. The subspace $L$ is
dense in $L_p(X,\mu)$ if and only if the $g\in\gen$ have $\mu$-almost no common
zeros.
\end{theorem}

\begin{proof}
In view of Theorem~\ref{thm:spacedist} the annihilator of $L$ can be identified
with those $f \in L_q(X,\mu)$ such that
\begin{equation}\label{eq:annihilator}
\int_X \psi(x) f(x)g(x) \,d\mu(x)=0
\end{equation}
for all $g\in\gen$ and $\psi\in\test$. Choose such $f$ and fix $g\in\gen$. Then
the left-hand side in \eqref{eq:annihilator} defines a continuous functional on
$C_c(\Omega)$ which vanishes on $\testcOmega$. By density we conclude that
\eqref{eq:annihilator} holds for all $\psi\in C_c(\Omega)$ or, as is
equivalent by Lemma~\ref{lem:restrictionsurjective}, for all $\psi\in C_c(X)$.
Consider the complex Borel measure $\nu$ on $X$ defined by
\begin{equation*}
\nu (E)=\int_E f(x)g(x)\,d\mu(x)
\end{equation*}
for Borel sets $E$. From the Radon--Nikodym theorem and \eqref{eq:annihilator} we have
\begin{equation*}
\int_X \psi(x) \,d\nu(x)=0\quad(\psi\in C_c(X)).
\end{equation*}
By the Riesz representation theorem this implies that the total variation of $\nu$ is zero, i.e.,
\begin{equation*}
\int_X \vert f(x)g(x)\vert \,d\mu(x)=0.
\end{equation*}
Therefore $fg=0$ almost everywhere ($\mu$). Since $g\in\gen$ was arbitrary we conclude that
$f$ vanishes $\mu$-almost everywhere on the complement of the common zero locus
of the $g\in\gen$. The converse statement for the annihilator is obvious. As to
$\overline L$, the necessity of the condition on the zero locus follows from the
fact that a convergent sequence in an $L_p$-space has a subsequence that
converges almost everywhere. The sufficiency is evident from the description of
$L^\perp$.
\end{proof}

\begin{remark}
Let $\mu^*$ be the completion of $\mu$ as obtained by including the subsets of
sets of measure zero in the domain of definition of the set function. By \cite
[p.~56, Theorem C]{Halmos} this completion $\mu^*$ can alternatively be defined as
the result of the Carath\'eodory extension procedure. Then the theorem also holds
for $\mu^*$ since $L_p(X,\mu)$ and $L_p(X,\mu^*)$ are canonically isomorphic
as a consequence of \cite[p.~143, Lemma 1]{RudinRCA}.
\end{remark}

\begin{corollary}\label{cor:Hamburger}
Let $\mu$ be a Borel measure on $X$ or the completion of such a measure and let
$1\leq p<\infty$. Let $\gen\subset L_p(X,\mu)_{\textup{\text{qa}}}$ be countable. Then the
module generated by $\gen$ over $\poly$ is dense in $L_p(X,\mu)$ if and only if
the $g\in\gen$ have $\mu$-almost no common zeros. The analogous statement holds
for the module generated over $\galg$ by a countable subset $\gen$ of $L_p(X,\mu)_0$.
\end{corollary}

\begin{corollary}\label{cor:polydenseLp}
Let $\mu$ be a Borel measure on $X$ or the completion of such a measure. Suppose

\begin{equation*}
\int_X w(x)^{-1}\,d\mu<\infty
\end{equation*}
for some quasi-analytic weight $w$ which is measurable on $X$. Then $\poly$ is
dense in $L_p(X,\mu)$ \textup{(}$1\leq p<\infty$\textup{)} and the same holds for the span of the
trigonometric functions $e_{i\lambda}$ corresponding to spectral parameters in a
somewhere dense set.
\end{corollary}

\begin{proof}
Replacing $w$ by a majorant as in Theorem \ref{thm:improvements} we may assume
that $w$ is strictly positive. Note that then $\op_{w^{1/p}}\cdot 1\subset
L_p(X,\mu)$. Since $w^{1/p}\in\weights_{\textup{\text{qa}}}$ we conclude that $1\in
L_p(X,\mu)_{\textup{\text{qa}}}$. The constant function $1$ has no zeros and the density
statements therefore follow from Theorem~\ref{thm:Lpspaces}.
\end{proof}

For $X=\Rn$, $p=2$ and the choice of a holomorphic weight $\exp (-\epsilon \Vert x\Vert)$ for some $\epsilon>0$ the polynomial part of
the corollary specializes to the assertion that the polynomials are dense in
$L_2(\Rn,\mu)$ if
\begin{equation}\label{eq:intcondition}
\int_\Rn e^{\epsilon \Vert x\Vert}\,d\mu<\infty
\end{equation}
for some $\epsilon>0$. This is the multidimensional analogue \cite[Theorem
3.1.18]{DunklXu} of Hamburger's theorem \cite{Hamburger}. According to the
corollary, the result can be strengthened, not only with respect to the
underlying point set, the trigonometric functions, and other values than $p=2$,
but most notably the integrability condition can be relaxed. The integrand in
\eqref{eq:intcondition} can, e.g., be replaced by any positive function which at
infinity is equal to $\exp (\epsilon \Vert x\Vert/\log a \Vert x\Vert)$ for some
$\epsilon,\,a>0$. The other examples of quasi-analytic weights in
Subsection~\ref{subsec:examplesqa} provide even more lenient integrability
conditions.

\subsection{Determinate multidimensional measures}

As an application of the preceding section we have the following result:

\begin{theorem}\label{thm:determinate}
Let $\mu$ be a probability measure on $\Rn$ such that
\begin{equation*}
\int_\Rn w(x)^{-1}\,d\mu<\infty
\end{equation*}
for some measurable quasi-analytic weight. Then $\mu$ is determinate.
Furthermore, $\poly$ is then dense in $L_p(\Rn,\mu)$ for all $1\leq p<\infty$ and
the same holds for the span of the trigonometric functions $e_{i\lambda}$
corresponding to spectral parameters in a somewhere dense set.
\end{theorem}

\begin{proof}
The density is a special case of Corollary \ref{cor:polydenseLp}. We see that
there exists $p_0>2$ such that $\poly$ is dense in $L_{p_0}(\Rn,\mu)$. It is
known that $\mu$ is then determinate \cite{Fuglede}.
\end{proof}

Note that the criterion in the above theorem is directly in terms of the measure
and circumvents, e.g., the computation of moments. For the choice of a
holomorphic weight $\exp (-\epsilon \Vert x\Vert)$ for some $\epsilon>0$ the
condition specializes to \eqref{eq:intcondition} and the corresponding result is
known \cite[Theorem 3.1.17]{DunklXu} but, as with the density of $\poly$ in
Corollary~\ref{cor:polydenseLp}, we see that the classical condition
\eqref{eq:intcondition} can be relaxed.

The determinacy of the measure can also be proved directly as an application of an extended
Carleman theorem in arbitrary dimension, using the classification of quasi-analytic weights on
$\Rn$. We will report separately on this in \cite{deJeu}.

\section{Interpretation in terms of Lie groups}\label{sec:interpretation}
We will now review the approach of the density results in Section
\ref{sec:density} in a formalism that indicates how the setup can be generalized
to certain classes of connected Lie groups. It will become apparent that such a
generalization yields results situated on the unitary dual of a Lie group rather
than on the group itself. It also shows that the results on $\Rn$ as we have
derived them are actually results on ${\widehat\R}^n$.

Given the indicative character of this section we are not precise about the
various topologies that one has to impose in order to validate the line of
reasoning below. The case of $\Rn$ in the present paper does not provide too
much information on this aspect, since various notions as well as topologies
coincide in that situation. We therefore concentrate on the formal structure and
leave a possible further concretization on the basis of a more thorough
investigation to the interested reader.

For background on the notions from Lie theory the reader is referred to, e.g.,
\cite{Knapp,Warner}.

\subsection{Formalism}\label{subsec:formalism}
Let $G$ be a connected Lie group with unit element $e$, universal enveloping
algebra $\env(G)$, and space of compactly supported smooth functions $\test(G)$.
Suppose $\{(\pi_\lambda,H_\lambda)\}_{\lambda\in\Lambda}$ is a collection of
unitary representations $\pi_\lambda$ of $G$ on Hilbert spaces $H_\lambda$. Let
$\Lambda^\infty=\coprod_{\lambda\in\Lambda} H_\lambda^\infty$ denote the
disjoint union of the smooth vectors in the various representations and put
$\End \Lambda^\infty=\coprod_{\lambda\in\Lambda}\End H_\lambda^\infty$, the
disjoint union of the linear transformations of these spaces of smooth vectors.
One has the naturally associated vector space $\Gamma\Lambda^\infty$ of sections
of $\Lambda^\infty$ and the algebra $\Gamma\End\Lambda^\infty$ of sections of
$\End\Lambda^\infty$. We let $\Gamma\End\Lambda^\infty$ act on
$\Gamma\Lambda^\infty$ by the obvious action in each fiber.

In $\Gamma\End\Lambda^\infty$ we have three natural subspaces---in fact
subalgebras---as images of a Fourier transform:
\begin{enumerate}
\item For $g\in G$ define its Fourier transform $\widehat
g\in\Gamma\End\Lambda^\infty$ by $\widehat g (\lambda)=\pi_\lambda (g)$. The
Fourier transform extends to the group algebra $\C[G]$.
\item For $D\in\env(G)$ define its Fourier transform $\widehat
D\in\Gamma\End\Lambda^\infty$ by $\widehat D(\lambda)=\pi_\lambda(D)$, where we
let $\pi_\lambda$ also denote the associated representation of $\env(G)$ on
$H_\lambda^\infty$.
\item For $f\in\test(G)$ define its Fourier transform $\widehat f\in\
\Gamma\End\Lambda^\infty$ by the fiberwise integral $\widehat f(\lambda)=\int_G
f(g)\pi_\lambda(g)\,dg$.
\end{enumerate}
Suppose that the three subalgebras $\widehat\env(G)$, $\widehat\C[G]$ and
$\widehat\test(G)$ of $\Gamma\End\Lambda^\infty$ are all contained in some
subspace $\op$ of $\Gamma\End\Lambda^\infty$. We assume that $\op$ carries a
Fr\'echet topology and aim at proving that these subspaces of $\op$ (or only the
second and the third to start with) have equal closure in $\op$. We do this by
showing that they are all dense. The crucial tool for this, as we will see, is
the Fourier transform $\Fourier:G\mapsto\op$ defined by $\Fourier(g)=\widehat
g$.

Before outlining the possible proof we illustrate how this setup corresponds to
the case of $\Rn$ in Section~\ref{sec:density}. In the case of $\Rn$ we choose
$\Lambda=\Rn$ (actually one should write $\Lambda={\widehat\R}^n$) with
$\pi_\lambda$ equal to the one-dimensional irreducible unitary representation
with spectral parameter $\lambda$. In this situation, matters are simplified---or
perhaps obscured---since both $\Gamma\Lambda^\infty$ and
$\Gamma\End\Lambda^\infty$ can be identified with functions on $\Rn$, the action
of the latter on the former being pointwise multiplication. Under this
identification pointwise multiplication by elements of $\galg$ now corresponds
to $\widehat\C[G]$, but with an important difference: the role of the variables
has been reversed compared to Section~\ref{sec:density}. Namely, in the present
picture the variable in the underlying point set $\Lambda$ is the parameter for
the unitary dual of $\Rn$, i.e., $\lambda$. Similarly, the image of the
universal enveloping algebra of $\Rn$ under the present Fourier transform now
corresponds to pointwise multiplication by elements of $\poly$ (with $\lambda$
as variable). The counterpart of $\op$ is pointwise multiplication by the
elements of one of our spaces $\opw$ (which in our approach could as well be
assumed to be Fr\'echet, as we have seen). The definition of $\Fourier$ in the
general context corresponds in the case of $\Rn$ to the map
$x\mapsto\{\lambda\rightsquigarrow\exp i(\lambda,x)\}$. This is the map $\psi_1$
of Section~\ref{sec:density}, but again there is a reversal of roles for
$\lambda$ and $x$.

The role of the space $\test(\Rn)$ as the basic dense subspace of $\opw$ is now
most naturally taken over by $\widehat\test(G)$. In the case of $\Rn$ the latter
space does not coincide with $\test(\Rn)$, but $\widehat\test(\Rn)$ is dense in
$\opw$ just as well (here the range result from Fourier analysis is essential!),
so we could equivalently have worked with this space from the very beginning.
Conceptually it is, as will become apparent, actually the natural space to work
with, but in case of $\Rn$ we preferred to use $\test(\Rn)$ as the basic
subspace since that is the most practical space for the applications on $\Rn$.

As a parallel to $\Rn$, in the general case we make the fundamental assumption that we
have already been able to show that $\widehat\test(G)$ is dense in $\op$. This
can analogously to $\Rn$ be expected to be basically a range result from
harmonic analysis for $G$.

Let us now indicate how one could prove in this general setup that
$\overline{\widehat\C[G]}=\overline{\widehat\test(G)}$ under the assumption that
$\Fourier:G\mapsto\op$ is continuous. For $\Rn$ this corresponds to proving
that $\overline\galg=\overline{\widehat{\test}}$ in $\opw$ if $w\in\weights_0$,
since under the latter condition on the weight we established in Proposition
\ref{prop:mapiscontinuous} that $\psi_1:\Rn\mapsto\opw$ is continuous. Consider for
$f\in\test(G)$ the weak integral
\begin{equation*}
I_f=\int_G f(g) \widehat g\,dg,
\end{equation*}
which exist in the Fr\'echet space $\op$. Suppose $T\in\op^\prime$ vanishes on
$\widehat\C[G]$. Then $T$ vanishes on all weak integrals $I_f$. On the other
hand, we have at least formally that $I_f(\lambda)=\widehat f(\lambda)$ for all
$\lambda\in\Lambda$, i.e., $I_f=\widehat f$. Therefore such $T$ vanishes on
$\widehat f$ for all $f\in\test(G)$, i.e., on a subspace which we had assumed to
be dense in $\op$. Hence $T=0$ as required.

Suppose next that we want to prove that
$\overline{\widehat\env(G)}=\overline{\widehat\C[G]}=\overline{\widehat\test(G)}$.
For $\Rn$ this corresponds to
$\overline\poly=\overline\galg=\overline{\widehat\test}$ in $\opw$, which we
were able to demonstrate for $w\in\weightsqa$. Let us assume that
$\Fourier:G\mapsto\op$ is smooth, which we concluded in the case of $\Rn$ for
$\psi_1$ in Corollary~\ref{cor:emapsmooth}, provided $w\in\weights_\infty$. Then
we already have $\overline{\widehat\C[G]}=\overline{\widehat\test(G)}$ by the
previous argument. In addition one verifies that formally the following relation
holds:
\begin{align}\label{eq:formalrelations}
D(T\circ\Fourier)(g)=\langle T,\widehat g\circ \widehat D\rangle\quad(g\in
G,\,D\in\env(G),\,T\in\op^\prime).
\end{align}
The necessary additional ingredient is now a theorem asserting that
$T\circ\Fourier$ has the quasi-analytic property at $e\in G$ for all
$T\in\op^\prime$. Indeed, if in that case $T\in\op^\prime$ vanishes on
$\widehat\env(G)\subset\op$ we conclude from
\eqref{eq:formalrelations} that $D(T\circ\Fourier)(e)=0$ for all $D\in\env(G)$.
The hypothesized quasi-analytic property then yields that $T\circ\Fourier=0$,
implying that $T$ vanishes on $\widehat\C[G]$. We had already concluded that the
latter space is dense in $\op$; therefore $T=0$ as required.

For the sake of completeness, let us note that for $T\in\op^\prime$ the map
$g\mapsto\langle T,\widehat g\rangle$ is a function on $G$ which (up to
normalization and a minus sign) was identified in Section~\ref{sec:density} as
the Fourier transform of the tempered distribution $T$ on $\Rn$. In the general
setup $T$ is now understood as a tempered distribution on ${\widehat\R}^n$. In
the general case the identification of the corresponding function appears to
have no natural analogue, reflecting our observation in the case of $\Rn$ (cf.\
Remark~\ref{rem:generalizations}) that the identification of the function
$\Fouriertyp$ as a Fourier transform is not essential, but that a range result
from harmonic analysis cannot be dispensed with.

The above framework constitutes the formal generalization to a connected Lie
group $G$ of our approach to the spaces $\opw$. Once this outline has been
completed for a particular group, one can exploit results on equal closure
obtained in such spaces $\op\subset\Gamma\End\Lambda^\infty$ to obtain closure
results in topological vector spaces $V\subset\Gamma\Lambda^\infty$ by applying
Theorem~\ref{thm:closurecont} to the fiberwise action of the elements of $\op$
on the elements of $V$, provided that the necessary continuity hypotheses are
satisfied by this action. Evidently, $V$ can also consist of equivalence classes
of sections in $\Gamma\Lambda^\infty$ (the analogue of, e.g., $L_p$-spaces) or
the action can be transposed to dual objects (the analogue of spaces of
distributions and measures).

The formalism as indicated explains at a conceptual level the intimate connection between the
spaces $\galg$ and $\poly$ in parts of classical approximation theory: the former space is the
group algebra of $\Rn$ and the latter is the universal enveloping algebra of $\Rn$, both acting
on the smooth vectors in the various irreducible unitary representations of the group. It is now
also apparent why these and other of the subspaces of $\opw$ that we encountered are actually
algebras. We did not actually use the multiplicative structure of such spaces---it merely seemed
to allow a convenient terminology in terms of modules---but it is now evident that these
spaces are the image (or the closure of an image) of naturally occurring algebras in Lie theory
under a Fourier transform. Therefore these subspaces of $\opw$ had to be algebras themselves.

The outline above makes it obvious that generalizations of our results on $\Rn$
can only be expected to be situated on sets of representation of groups. Since
abelian groups can be viewed as unitary duals, such groups form an exceptional
class in the sense that one can obtain results situated on these groups
themselves. From the present point of view this is a coincidence. The case of
the self-dual $\Rn$ is in a sense particularly deceptive since a maximal number
of various different structures corresponding to the general framework then
coincide.

\subsection{Perspectives}\label{subsec:perspectives}
We now make some tentative comments on the practical feasibility and context of
the above outline.

First of all, a technical point of some importance is the quasi-analytic
property on Lie groups: What conditions on a smooth function $f$ on a connected
Lie group $G$ are sufficient to ensure $f$ vanishes identically if $Df(e)=0$ for
all $D\in\env(G)$? If the exponential map is surjective then one can in
principle resort to $\Rn$ for knowledge on quasi-analytic classes. In the case
of arbitrary $G$ one can hope to be able to use the local character of theorems
in the vein of Theorem~\ref{thm:DCtheorem}. A theorem on quasi-analytic classes
on general Lie groups however does not seem to be available at present.

As to the practical context of the general scheme, it stands to reason that
additional conditions on $\Lambda$ and the various topologies are necessary to
obtain meaningful applications. The case of $\Rn$ provides the following clues.

First of all, our definition of the spaces $\opw$ for $\Rn$ involved derivatives
(now with respect to $\lambda$) and the behavior at infinity. This suggests
that the representations corresponding to the elements of $\Lambda$ should be
related to each other and to a topology on $\Lambda$---as was only to be
expected.

More interesting is the following trivial application of our method for $\Rn$.
Let $c_1,\,c_2\in\C$ and $\lambda_0\in\Rn$ be such that $c_1\neq 0$. Choose
$f_1\in\schw_{\textup{\text{qa}}}$ (see Subsection~\ref{subsec:schwartz} for notation) without
zeros and such that $f_1(\lambda_0)=c_1$. This is possible since $c_1\neq 0$.
Choose $f_2\in\schw$ such that $f_2(\lambda_0)=c_2$. Then $\poly\cdot f_1$ is
dense in $\schw$ as a consequence of the results in Subsection
\ref{subsec:schwartz}. In particular, if $\epsilon>0$ is given, then there exists
a polynomial $P_\epsilon$ such that $\Vert f_1 - P_\epsilon\cdot
f_2\Vert_\infty<\epsilon$. Specializing to $\lambda_0$ yields $\vert
c_2-P_\epsilon(\lambda_0)\cdot c_1\vert<\epsilon$. One can obviously replace the
Schwartz space by other spaces of functions and obtain the same result in a
similar fashion. This result, though trivial for $\Rn$, has an interpretation in
the formalism for the general case. By virtue of its type of proof, which will
have an obvious analogue if meaningful applications as in the case of $\Rn$ are
possible, it can now be interpreted as the density in $H_{\lambda_0}^\infty$ of
the orbit of any nonzero smooth vector under the universal enveloping algebra,
and this for arbitrary $\lambda_0\in\Lambda$. Since $H_{\lambda_0}^\infty$ is in
turn dense in $H_{\lambda_0}$ and since the above density of the orbit is in
particular valid for analytic vectors, one is led to the conclusion that the
analytic vectors are dense in each $H_{\lambda_0}$, i.e., that all $\pi_\lambda$
($\lambda\in\Lambda$) have to be irreducible. This locates the spaces $\op$ on
the unitary dual of $G$. Perhaps it is already sufficient to assume that the
$\pi_\lambda$ are only generically irreducible, thus admitting a somewhat wider
context.

In view of all the above we tentatively expect that a generalization of our
method can be established for (a subset consisting of irreducible
representations in) a series of (possibly induced) unitary representations for
reductive groups. The duals of connected nilpotent groups appear to be another
group of candidates. Once established, relatively easy applications as in
Section~\ref{sec:applications} are within reach, each situated on subsets of the
unitary dual or perhaps (if the representations are only generically
irreducible) on slightly larger sets. We leave it to the reader to determine whether the relative ease of the possible
applications merits the effort of establishing the generalization in particular
cases of interest.

To conclude we note that the above formalism in principle still admits a further
generalization along the same lines to, e.g., nonunitary Hilbert representations
or Fr\'echet representations and their smooth vectors.

\subsection{The $n$-torus}
The general outline can obviously be completed for the $n$-torus $\Tn$, yielding
results situated on $\widehat \T^n\cong\Zn$. On the other hand, since $\Zn\subset\Rn$, nontrivial
results on $\Zn$ can also be obtained immediately from the results for the spaces $\opw$ on $\Rn$, simply by restricting the elements of $\opw$ to $\Z^n$. This evidently reflects the canonical inclusion $\widehat \T^n\subset\widehat \R^n$. As the reader may verify, stronger statements than these statements by restriction could be obtained by actually completing the outline, provided that one had results on generalized quasi-analyticity (in the sense of \cite{Mandelbrojt}) for smooth functions on $\Tn$, which are substantially stronger than those for smooth functions on $\Rn$. This would lead to M\"untz--Szasz-type theorems, stating that the closure of the span of ``sufficiently many'' monomials in a function space $\op_w$ on $\widehat T^n\cong\Zn$ is equal to that of the closure of the subspace of finitely supported elements of $\opw$ (since these spaces are in fact both dense). Such substantially stronger results on generalized quasi-analyticity for periodic smooth functions on $\Rn$ however do not appear to be available. In view of this, the results on $\widehat\T^n$ obtained by restriction from $\widehat\R^n$ as above appear to be optimal within the framework of our method, at least for practical purposes.

\appendix

\section{Equal closure and continuous maps}\label{app:closurecont}

This appendix is concerned with the relation between continuous maps and the
property of having equal closure. It supplements Lemma
\ref{lem:fundamentallemma}.

\begin{proposition}\label{prop:sequentialpreserved}
If $X$ and $Y$ are topological vector spaces \textup{(}not necessarily Hausdorff\textup{)} and
$0\in X$ has a countable neighborhood base, a linear map $\phi:X\mapsto Y$ is
continuous if and only if, for all subsets $A$ and $B$ of $X$ having equal
closure in $X$, the images $\phi(A)$ and $\phi(B)$ again have equal sequential
closure in $Y$.
\end{proposition}

Note that closure and sequential closure coincide in $X$.

\begin{proof}
As to the ``only if'' part, assume that $\phi$ is continuous and $\overline
A=\overline B$. It is sufficient to show that $y\in{\overline{\phi(B)}}^s$ for
$y\in{\overline{\phi(A)}}^s$. Let $\{a_n\}_{n=1}^\infty$ be a sequence in $A$
such that $\lim_{n\rightarrow\infty}\phi(a_n)=y$. Let $\{U_n\}_{n=1}^\infty$ be
a neighborhood base of $0$ in $X$. We may assume that $U_1\supset U_2\supset
U_3\supset\ldots\,$. Since $a_n\in {\overline{A}}={\overline{B}}$ for all $n$,
there exists for all $n$ an element $b_n\in B$ such that $a_n-b_n\in U_n$. Then
$\lim_{n\rightarrow\infty}\phi(b_n)=y$. To see this, let $V_0$ be an arbitrary
neighborhood of $0$ in $Y$. Choose a neighborhood $V_1$ of $0$ in $Y$ such
that $V_1+V_1\subset V_0$. There exists $N_1$ such that $y-\phi(a_n)\in V_1$ for
all $n\geq N_1$. By the continuity of $\phi$, there exists $N_2$ such that
$\phi(U_{N_2})\subset V_1$. Then for all $n\geq\max(N_1,N_2)$ we have
$y-\phi(b_n)=(y-\phi(a_n))+\phi(a_n-b_n)\in V_1 +
\phi(U_n)\subset V_1 +
\phi(U_{N_2})\subset V_1+V_1\subset V_0$.

As to the ``if'' part, let $\phi$ preserve the property of having equal
(sequential) closure. Suppose that $\phi$ is not continuous at $0$. Then there
exists an open neighborhood $V$ of $0$ in $Y$ such that $\phi^{-1}(V)$ is not a
neighborhood of $0$ in $X$. Using a decreasing neighborhood base at $0\in X$
as above, one easily constructs a sequence $\{x_n\}_{n=1}^\infty\subset X$ such
that $\phi(x_n)\in Y-V$ for all $n$ and such that $x_n\rightarrow 0$. Then
$\overline {\bigcup_{n=1}^\infty\{x_n\}\cup\{0\}}=\overline
{\bigcup_{n=1}^\infty\{x_n\}}$, so $0=\phi(0)\in
{\overline{\bigcup_{n=1}^\infty\{\phi(x_n)\}\cup
\{\phi(0)\}}}^s = {\overline{\bigcup_{n=1}^\infty\{\phi(x_n)\}}}^s\subset
\overline{\bigcup_{n=1}^\infty\{\phi(x_n)\}}\subset \overline{Y-V}=Y-V$. But $0\notin Y-V$
by construction. Contradiction.
\end{proof}

\begin{theorem}\label{thm:closurecont}
Let $Y$ be a topological vector space. Suppose that $I$ is a nonempty index set
and assume that for each $i\in I$ a topological vector space $X_i$ and a
continuous linear map $\phi_i:X_i\mapsto Y$ are given.

Let $\choice$ be a choice function, assigning to each $i\in I$ a dense subspace
$\choice_i$ of $X_i$, and consider the corresponding algebraic span
$L_\choice=\sum_{i\in I}\phi_i(\choice_i)$ of the images $\phi_i(\choice_i)$.
Then:
\begin{enumerate}
\item $L_\choice^\perp$ is independent of $\choice$. In fact, if $\choice^\prime$ is any other such choice function,
then $L_\choice^\perp=\bigcap_{i\in I}\phi_i(\choice^\prime_i)^\perp$.
\item If $Y$ is locally convex, then the closure of $L_\choice$ is independent of $\choice$.
\item If $I$ is finite, then the closure of $L_\choice$ is independent of $\choice$.
\item If $I$ is finite and the $X_i$ all have a countable neighborhood base at $0$, then the sequential closure of $L_\choice$ is independent of $\choice$.
\end{enumerate}
\end{theorem}

\begin{proof}
As to (1), invoking Lemma~\ref{lem:fundamentallemma} yields that
\begin{equation*} L_\choice^\perp=\bigcap_{i\in I}\phi_i(\choice_i)^\perp=
\bigcap_{i\in I}\left\{\overline{\phi_i(\choice_i)}\right\}^\perp=\bigcap_{i\in I}\left\{\overline{\phi_i(\choice_i^\prime)}\right\}^\perp=\bigcap_{i\in I}\phi_i(\choice_i^\prime)^\perp.
\end{equation*}
This proves (1) and (2) is then immediate. For the remaining statements assume
$I=\{1,\ldots,n\}$. Consider $\widehat{X}=\bigoplus_{i=1}^n X_i$ in its product
topology and define the continuous linear map $\phi:\widehat{X}\mapsto Y$ by
$\phi(x_1,\ldots,x_n)=\sum_{i=1}^n\phi_i(x_i)$. Consider the subspace
$\widehat{M}_\choice=\bigoplus_{i=1}^n \choice_i$ of $\widehat X$. Since
${\overline{\widehat M}}_\choice=\widehat X$, Lemma~\ref{lem:fundamentallemma}
implies (3). If the $X_i$ are first countable, then so is $\widehat{X}$ and (4)
therefore follows from Proposition~\ref{prop:sequentialpreserved}.
\end{proof}

\section{A Denjoy--Carleman-type theorem in several variables}\label{app:DCtheorem}

In this appendix we prove a theorem on quasi-analytic classes in several
variables. If the $M(j,m)$ in the formulation of the theorem are all strictly
positive, then with the aid of \cite[Theorem 1.8.VII]{Mandelbrojt} the result
could also be inferred from \cite{Hryptun}. In our application of the theorem
below in the proof of Theorem~\ref{thm:Fouriertransforms} the $M(j,m)$ are,
however, only known to be nonnegative. Therefore we provide an independent
proof, which is perhaps somewhat simpler than
 the proof in \cite{Hryptun} where
a proof is given by the repeated use of Bang's formulas. We reduce the theorem
to the one-dimensional case by induction.

\begin{theorem}\label{thm:DCtheorem}
For $j=1,\ldots,n$, let $\{M(j,m)\}_{m=0}^\infty$ be a sequence of nonnegative
real numbers. Put $\mu(j,m)=\inf_{k\geq m} M(j,k)^{1/k}$
$(j=1,\ldots,n;\,m=1,2,\ldots)$ and suppose
\begin{equation*}
\sum_{m=1}^\infty \frac{1}{\mu(j,m)}=\infty\quad(j=1,\ldots,n).
\end{equation*}
Let $R>0$ and suppose $f:(-R,R)^n\mapsto\C$ is of class $C^\infty$. Assume that
there exist nonnegative constants $C$ and $r$ such that
\begin{equation*}
\vert \parta f(x)\vert\leq Cr^{\vert \alpha\vert}\prod_{j=1}^n M(j,\alpha_j)
\end{equation*}
for all $\alpha\in\Nn$ and all $x\in (-R,R)^n$. Then, if $\parta f(0)=0$ for all
$\alpha\in\Nn$, $f$ is actually identically zero on $(-R,R)^n$.
\end{theorem}

\begin{proof}
By induction. For $n=1$ the result is the Denjoy--Carleman theorem if the
$M(1,m)$ are all strictly positive \cite{Carleman,Mandelbrojt}. If $M(1,m_0)=0$
for some $m_0>0$, then $f$ is a polynomial of degree at most $m_0-1$ and the result
is obvious. If $M(1,0)=0$ then the statement is trivial, completing the case $n=1$.

Assuming the theorem for $n-1$, we define for each multi-index
$(\alpha_1,\ldots,\alpha_{n-1})\in\N^{n-1}$ the function
$\phi_{\alpha_1,\ldots,\alpha_{n-1}}:(-R,R)\mapsto\C$ by
$$\phi_{\alpha_1,\ldots,\alpha_{n-1}}(t)=\partial_1^{\alpha_1}\cdots\partial_{n-1}^{\alpha_{n-1}}f(0,\ldots,0,t),$$
where we have used the shorthand notation $\partial_j=\partial/\partial x_j$.
Then all derivatives of $\phi_{\alpha_1,\ldots,\alpha_{n-1}}$ vanish at
$0\in\R$. Since for $m=0,1,\ldots$ and $t\in(-R,R)$ we have
\begin{equation*}
\vert (d/dt)^m\phi_{\alpha_1,\ldots,\alpha_{n-1}}(t)\vert\leq
\left(Cr^{\sum_{j=1}^{n-1}\alpha_j}\prod_{j=1}^{n-1}M(j,\alpha_j)\right)r^m M(n,m),
\end{equation*}
the result as established for $n=1$ implies that
$\phi_{\alpha_1,\ldots,\alpha_{n-1}}$ is identically zero on $(-R,R)$, for
arbitrary $(\alpha_1,\ldots,\alpha_{n-1})\in\N^{n-1}$.

Next, for each $t\in (-R,R)$ define $\psi_t:(-R,R)^{n-1}\mapsto\C$ by
$\psi_t(x_1,\ldots,x_{n-1})=f(x_1,\ldots,x_{n-1},t)$. Since
$\partial_1^{\alpha_1}\cdots\partial_{n-1}^{\alpha_{n-1}}\psi_t(0,\ldots,0)=\phi_{\alpha_1,\ldots,\alpha_{n-1}}(t)$,
the vanishing of all $\phi_{\alpha_1,\ldots,\alpha_{n-1}}$ on $(-R,R)$ implies
that all derivatives of $\psi_t$ vanish at $0\in\R^{n-1}$, for arbitrary
$t\in(-R,R)$. Additionally, we have, for all
$(x_1,\ldots,x_{n-1})\in(-R,R)^{n-1}$ and all $t\in (-R,R)$,
\begin{equation*}
\vert \partial_1^{\alpha_1}\cdots\partial_{n-1}^{\alpha_{n-1}}\psi_t(x_1,\ldots,x_{n-1})\vert\leq
\left(C\, M(n,0)\right)\cdot r^{\sum_{j=1}^{n-1}\alpha_j}\prod_{j=1}^{n-1}M(j,\alpha_j).
\end{equation*}
For each $t\in(-R,R)$ the induction hypothesis then implies that $\psi_t$ is
identically zero on $(-R,R)^{n-1}$, i.e., $f$ vanishes on $(-R,R)^n$. This
completes the induction step.
\end{proof}





\begin{thebibliography}{99}

\bibitem{Bernstein}
S.N.~Bernstein, On weight functions (Russian), {\it Dokl.\ Akad.\ Nauk SSSR}
{\bf 77} (1951), 549--552.

\bibitem{Bourbaki1}
N.~Bourbaki, ``Topologie G\'en\'erale I,'' Hermann, Paris, 1971.

\bibitem{Bourbaki2}
N.~Bourbaki, ``Topologie G\'en\'erale II,'' Hermann, Paris, 1971.

\bibitem{BourbakiTVS}
N.~Bourbaki, ``Topological Vector Spaces, Chapters 1-5,'' Springer Verlag,
Berlin-Heidelberg-New York-London-Paris-Tokyo, 1987.

\bibitem{Carleman}
T.~Carleman, ``Les fonctions quasi-analytiques,'' Gauthiers-Villars, Paris,
1926.

\bibitem{Carleson}
L.~Carleson, On Berstein's approximation problem, {\it Proc.\ Amer.\ Math.\ Soc.\ }{\bf 2} (1951), 953--961.

\bibitem{Chirka}
E.M.Chirka, ``Complex analytic sets,'' Kluwer Academic Publishers,
Dordrecht-Boston-London, 1989.

\bibitem{DunklXu}
C.F.~Dunkl and Y.~Xu, ``Orthogonal polynomials of several variables,'' Cambridge
University Press, Cambridge-New York-Oakleigh-Madrid-Cape Town, 2001.

\bibitem{Fuglede}
B.~Fuglede, The multidimensional moment problem, {\it Expo.\ Math.\ }{\bf 1} (1983), 47--65.

\bibitem{Hall}
T.~Hall, Sur l'approximation polyn\^omiale des fonctions continues d'une
variable r\'eelle, in ''9e Congr\`es Math.\ Scand., Helsingfors 1938,''
Mercator, Helsingfors, 1939.

\bibitem{Halmos}
P.R.~Halmos, ``Measure Theory,'' Springer Verlag, New York-Heidelberg-Berlin,
1974.

\bibitem{Hamburger}
H.~Hamburger, Beitr\"age zur Konvergenztheorie der Stieltjesschen Kettenbr\"uche, {\it Math.~Z.\
}{\bf 4} (1919), 186--222.

\bibitem{Horvath1}
J.~Horv\'ath, L'approximation polynomiale sur un ensemble non compact, {\it Math.\
Scand.\ }{\bf 2} (1954), 83--90.

\bibitem{Horvath2}
J.~Horv\'ath, ``Topological Vector Spaces and Distributions,'' Addison-Wesley,
Reading-Palo Alto-London-Don Mills, 1966.


\bibitem{Hryptun}
V.G.~Hryptun, An addition to a theorem of S.~Mandelbrojt, {\it Ukrain. Mat.\ Z.\ }{\bf 28}
(1976), 841--844. English translation: {\it Ukrainian Math.\ J.\ }{\bf 28} (1976), 655--658.


\bibitem{deJeu}
M.F.E.\ de Jeu, Determinate multidimensional measures, the extended Carleman theorem and
quasi-analytic weights, to appear in {\it Ann.\ Probab.}

\bibitem{Koosis}
P.~Koosis, ``The logarithmic integral I,'' Cambridge University Press,
Cambridge-New York-New Rochelle-Melbourne-Sydney, 1988.

\bibitem{Knapp}
A.W.~Knapp, ``Representation theory of semisimple groups, an overview based on
examples,'' Princeton University Press, Princeton, 1986.

\bibitem{Larsen}
R.~Larsen, ``Banach algebras, an introduction,'' Marcel Dekker, New York, 1973.

\bibitem{Lin}
G.D.~Lin, On the moment problems, {\it Statistics \& Probability Letters} {\bf
35} (1997), 85--90.

\bibitem{Mandelbrojt}
S.~Mandelbrojt, ``S\'eries adh\'erentes, r\'egularisation des suites,
applications,'' Gauthiers-Villars, Paris, 1952.

\bibitem{MandelbrojtRice}
S.~Mandelbrojt, ``General theorems of closure,'' The Rice Institute Pamphlet
Special Issue, The Rice Institute, Houston, 1951.

\bibitem{Mergelyan}
S.N.~Mergelyan, Weighted approximation by polynomials, {\it Uspekhi Mat.\ Nauk
}{\bf 11} (1956), 107--152. English translation: {\it Amer.\ Math.\ Soc.\
Transl.\ (2) }{\bf 10} (1958), 59--106.

\bibitem{Nachbin1}
L.~Nachbin, ``Elements of approximation theory,'' Van Nostrand, New
York-Chicago-San Francisco, 1967.

\bibitem{Nachbin2}
L.~Nachbin, On the weighted approximation of continuous differentiable
functions, {\it Proc.\ Amer.\ Math.\ Soc.\ }{\bf 111}, (1991), no. 2, 481--485.

\bibitem{Nachbin3}
L.~Nachbin, Sur les alg\`ebres denses de fonctions diff\'erentiables sur une
vari\'et\'e, {\it C.~R.\ Acad.\ Sci.\ Paris }{\bf 228} (1949), 1549--1551.

\bibitem{Ostrowski}
A.~Ostrowski, \"Uber quasianalytische Funktionen und Bestimmtheit aymptotischer Entwicklungen, {\it Acta Math.\ }{\bf 53} (1929), 181--266.


\bibitem{RudinFA}
W.~Rudin, ``Functional Analysis,'' 2nd edition, McGraw-Hill, Singapore, 1991.

\bibitem{RudinRCA}
W.~Rudin, ``Real and Complex Analysis,'' McGraw-Hill, New York, 1985.

\bibitem{Schwartz}
L.~Schwartz, ``Th\'eorie des distributions,'' Hermann, Paris, 1978.

\bibitem{Treves}
F.~Treves, ``Topological Vector Spaces, Distributions and Kernels,'' Academic
Press, New York-London, 1967.

\bibitem{Warner}
G.~ Warner, ``Harmonic analysis on semi-simple Lie groups I, II,'' Springer
Verlag, Berlin-Heidelberg-New York, 1972.


\bibitem{Whitney}
H.~Whitney, On ideals of differentiable functions, {\it Amer.\ J.\ Math.\ }{\bf
70} (1948), 635--658.

\bibitem{Zapata}
G.~Zapata, Bernstein approximation problem for differentiable functions and quasi-analytic
weights, {\it Trans.\ Amer.\ Math.\ Soc.\ }{\bf 182} (1973), 503--509.

\end{thebibliography}
\end{document}